\documentclass[12pt]{article}
\usepackage{blindtext}
\usepackage{tikz-cd}
\usepackage{tikz}

\usepackage{amsmath}
\usepackage{amscd}
\usepackage{amssymb}

\numberwithin{equation}{section}

\newtheorem{thm}{Theorem}[section]

\newtheorem{lem}[thm]{Lemma}
\newtheorem{prop}[thm]{Proposition}

\newcommand{\thmref}[1]{Theorem~\ref{#1}}
\newcommand{\propref}[1]{Proposition~\ref{#1}}

\newcommand{\lemref}[1]{Lemma~\ref{#1}}

\newsavebox{\SmallMathBox}

\def\pdo{\psi{\rm do}}

\def\Ci{C^\infty}

\def\dd{\partial}

\def\Di{D\kern -.65em /}
\def\Dii{D\kern -.45em /}
\def\di{{\dd}\kern -.55em /}
\def\dii{{\dd}\kern -.40em /}

\def\na{\nabla}

\def\ol{\overline}

\def\to{\rightarrow}
\def\too{\longrightarrow}
\def\mto{\mapsto}

\def\x{\times}

\def\re{{\rm Re}}

\def\NN{{\bf N}}

\def\QQ{{\bf Q}}

\def\Aa{{\mathcal A}}

\def\Cc{{\mathcal C}}
\def\Dd{{\mathcal D}}
\def\Ee{{\mathcal E}}
\def\Ff{{\mathcal F}}

\def\Ii{{\mathcal I}}

\def\Pp{{\mathcal P}}

\def\Sss{{\mathcal S}}

\def\Vv{{\mathcal V}}
\def\Ww{{\mathcal W}}
\def\Zz{{\mathcal Z}}

\def\={\cong}
\def\>{\supset}
\def\<{\subset}
\def\ii{^{-1}}
\def\si{^{-s}}

\def\12{\frac{1}{2}}

\def\6{\cup}
\def\ox{\otimes}

\def\dl{\mathchar'26\mkern-12mu d}

\def\uR{\underline{\mathbb R}}

\def\a{\alpha}

\def\R{\mathbb{R}}
\def\C{\mathbb{C}}
\def\Z{\mathbb{Z}}
\def\N{\mathbb{N}}

\def\d{\delta}
\def\D{\Delta}

\def\g{\gamma}
\def\G{\Gamma}

\def\la{\lambda}
\def\La{\Lambda}

\def\N{\NN}
\def\o{\infty}

\def\Q{\QQ}

\def\Si{\Sigma}

\def\ch{\mbox{\rm ch}}

\def\End{\mbox{\rm End}}

\def\Hom{\mbox{\rm Hom}}

\def\index{\mbox{\rm index\,}}

\def\Ker{\mbox{\rm Ker}}
\def\ker{\mbox{\rm ker}}

\def\Si{S\kern -.65em /}

\def\str{{\rm Str}}

\def\Str{{\rm Str}}

\def\Tr{{\rm Tr}}

\def\Strmb{{\rm Str}_{\mbox{\tiny M/B}}}

\def\Trmb{{\rm Tr}_{\mbox{\tiny M/B}}}

\def\tr{\mbox{\rm tr\,}}

\def\Af{\mathbb{A}}
\def\Bf{\mathbb{B}}
\def\Cf{\mathbb{C}}

\def\If{\mathbb{I}}

\def\Kf{\mathbb{K}}

\def\Mf{\mathbb{M}}
\def\Nf{\mathbb{N}}

\def\Ps{\mathbb{P}}

\def\Zf{\mathbb{Z}}

\def\tm{{\mbox{\tiny{M}}}}
\def\tn{{\mbox{\tiny{N}}}}
\def\tb{{\mbox{\tiny{B}}}}
\def\tmb{{\mbox{\tiny{M/B}}}}
\def\smb{{\mbox{\small{M/B}}}}
\def\tpb{{\mbox{\tiny{P/B}}}}
\def\tnb{{\mbox{\tiny{N/B}}}}
\def\tmnb{{\mbox{\tiny{$(M\times_B N)/B$}}}}

\def\ord{\mbox{\rm ord}}

\def\bs{\backslash}

\def\cc{\textsf{c}}

\def\As{\textsf{A}}
\def\Bs{\textsf{B}}
\def\Cs{\textsf{C}}
\def\Ds{\textsf{D}}
\def\Es{\textsf{E}}
\def\Fs{\textsf{F}}
\def\Gs{\textsf{G}}
\def\Hs{\textsf{H}}
\def\Is{\textsf{I}}
\def\Js{\textsf{J}}
\def\Ks{\textsf{K}}
\def\Ls{\textsf{L}}

\def\Ps{\textsf{P}}
\def\Qs{\textsf{Q}}
\def\Rs{\textsf{R}}

\def\Ts{\textsf{T}}
\def\Us{\textsf{U}}
\def\Vs{\textsf{V}}
\def\Ws{\textsf{W}}

\renewcommand{\r}{\right}
\renewcommand{\l}{\left}

\newcommand{\ww}{{\bf w}}

\newcommand{\sslash}{\mathbin{/\mkern-6mu/}}

\begin{document}

\title{The Poisson summation formula for \\[2mm] fibrations of Riemannian manifolds}

\author{S.G.Scott}

\date{}

\maketitle

\section{Introduction}

Let $(X,g)$ be a closed Riemannian manifold and let  $\D$ be the metric Laplacian on $X$ with eigenvalues $\la$. The wave trace
\begin{equation}\label{wave tr}
\Tr(e^{it\sqrt{\D}}) = \sum_{\la} e^{it\sqrt{\la}}   
\end{equation} 
is a tempered distribution in $t$ whose singular support is contained in the set of lengths $L$ of closed geodesics on $(X,g)$.
The wave trace Poisson summation formula on $X$ due to Chazarain \cite{C} and  Duistermaat-Guillemin \cite{DG}  computes the singularity structure of  \eqref{wave tr} in terms of distributions in $t$ of  increasing homogeneity, and  regularity. The formula  says that there is an asymptotic singularity expansion 
\begin{eqnarray}\label{wave Poisson summation}
\Tr(e^{it\sqrt{\D}}) & \sim & a_{\,0,-n+1}(t+i0)^{-n+1} + a_{\,0,-n+2}(t+i0)^{-n+2} + \cdots  \nonumber \\
 &+&  \sum_{L>0} \Biggl( a_{\,L,-1}(t-(L+i0))^{-1}\biggr.  \\
 &&    \hskip 10mm + \ \ \left.\sum_{j\geq 0}  a_{\,L,j}(t-(L+i0))^j \log(t-(L+i0))\right).  \nonumber
\end{eqnarray}
The coefficients $a_{\,0,-n+j}$ at the $t=0$ singularity correspond to the coefficients of the heat trace expansion and are relatively straightforward to compute, in terms of integrals of curvature polynomials, while the wave invariants $a_{\,L,j}$ for the singularity at  $t=L>0$ are more subtle 
dynamical curvature invariants, with detailed formulae  given in fundamental works of Guillemin \cite{G} and Zelditch \cite{Z}. \eqref{wave Poisson summation} was extended to manifolds with boundary by Guillemin and Melrose \cite{GM}.\\

Invariants associated to a manifold (of a given class)  often admit a natural extension to corresponding fibrewise invariants for a fibration of manifolds
 $$\pi: M\to B$$
(with fibres of that given class). Here, $M$ and $B$ are compact boundaryless manifolds and $\pi$ a smooth submersion; equivalently,  $\pi$ is a fibre bundle with smooth diffeomorphic fibres $M_b := \pi\ii(b)$. \footnote{For compact manifolds there is no distinction between a `fibration' and a `fibre bundle'.}
Higher invariants of this type, defined vertically along the fibres, 
reside on the parameter space $B$  as functions, differential forms or distributions, or currents, or as cohomology classes. Higher spectral-geometric invariants are the class of higher invariants associated to a smooth family of  Fourier integral operators (FIOs) along the fibres
$$\As = \As_{[0]} + \cdots + \As_{[n]}$$
with $\As_{[k]}$ the component which raises form degree on $B$ by $k$, and $n=\dim B$; in local coordinates $(b_1,\ldots, b_n)$ on $B$ 
\begin{equation}
\As =  \sum _I A_{I,b} \,db_I
\end{equation}
with $db_I  = db_{k_1}\wedge  \cdots \wedge db_{k_r}\in \Aa^{|I|}(B)$ a basis of differential forms for $B$ running over ordered  multi-indices $I = (k_1,\ldots, k_r)\in\mathbb{N}^r$ while $A_{I,b}$ is a FIO  on $M_b$ in the classical sense depending smoothly on $b\in B$, the summands with $|I| = k$ corresponding to $\As_{[k]}$. Here, $\Aa^*(N)$ is the de Rham algebra of differential forms on a manifold $N$\\

Let $\pi: M\to B$ be a smooth fibration of Riemannian manifolds of fibre dimension $q$. 
Let $\Qs$ be a smooth family of classical pseudo-differential operators of non-negative integer order along the fibres of $\pi$ with 
\begin{equation}\label{fibre Lap}
    \Qs_{[0]} = \D
\end{equation} 
(the smooth family of metric Laplacians defined by the Riemannian metric on each fibre $(M_b,g_b)$).  We assume smooth maps $L_\kappa : B\to\R^1$ such that 
 $\{L_\kappa(b) \ | \ \kappa\in\Z\}$ is the discrete set of periods of closed geodesics of length $L_\kappa(b)$ for the Riemannian metric $g_b$  on the fibre $M_b$. The higher wave-trace Poisson summation formula is then as follows.
 
\begin{thm}\label{thm1}
 There is an asymptotic singularity expansion   
\begin{equation}\label{fibrewise Poisson summation}
    \Tr_\tmb\l(e^{-it\sqrt{\Qs}}\r)  \ \sim \ \sum_{d=0}^{\dim B}\sum_{\kappa\in\Z}  \sum_{j\geq 0}  
{\bf a}_{\kappa,j,[d]}\,(t- L_\kappa + i0)^{-q + j} \hskip 25mm 
\end{equation}
$$\hskip 35mm + \ \sum_{d=0}^{\dim B}\sum_{\kappa\in\Z}  \sum_{l\geq 0}  
{\bf a}^{\,\prime}_{\kappa,l,[d]}\,(t- L_\kappa + i0)^{-q + l} \log(t- L_\kappa + i0)$$
with differential form wave invariants 
$$  {\bf a}_{\kappa,j, [d]}, \  {\bf a}^{\,\prime}_{\kappa,l, [d]}\ \in\Aa^d(B)$$
and $\Trmb$ the vertical trace along the fibres. The  coefficient of  $(t- L_\kappa + i0)^{-q + j}$
$$  \sum_d {\bf a}_{\kappa,j, [d]} \ = \  \int_{\Zz/B}  \Omega_{\kappa,j} \in\Aa^*(B)$$
is a differential form of mixed form degree for $\Zz\to B$ a subfibration of the cosphere bundle $S^*_\tmb $ along the fibres (the fibration over $B$ with fibre $S^*M_b$ at $b\in B$) with $\Omega_{\kappa,j}\in\Aa^*(\Zz)$ determined by the vertical Lagrangian distribution kernel of the v-FIO $e^{{-it\sqrt{\Qs}}}$, and similarly for the coefficients  ${\bf a}^{\,\prime}_{\kappa,l,[d]}$.\\

The vertical wave invariants $\textsf{a}_{\,0,j, [d]}\in\Aa^*(B)$  at the singularity at $L_0:=0$ coincide with the differential form coefficients in the asymptotic expansion of the vertical heat trace  $\Trmb(e^{-t\Qs})$ as $t\to 0+$. \\

Since $\D$ is a vertical differential operator then the sum over the $L_0=0$ terms has no non-zero log terms and the $d=0$ component of \eqref{fibrewise Poisson summation} coincides pointwise with \eqref{wave Poisson summation}; if $\Qs_{[0]} = \Ps$ is a positive self-adjoint elliptic pseudodifferential operator of integer order $m>0$ then the corresponding vertical wave trace  $\Tr_\tmb\l(e^{-it\Qs^{1/m}}\r)$ has a similar expansion to \eqref{fibrewise Poisson summation} and which may have non-zero log terms in the $ L_0=0$ sum, and the fibre integrals are over corresponding subfibrations of $S^*_\tmb$.
\end{thm}

Here, the fibrewise (half) wave operator $e^{it\sqrt{\Qs}}$ is a smooth family of FIOs along the fibres of $\pi$, a vertical FIO (or v-FIO), which is seen to coincide pointwise at form degree 0  with the classical wave operators \eqref{wave tr} $$(e^{it\sqrt{\Qs}})_{[0]} = e^{it\sqrt{\D}}.$$ Ellipticity properties of a v-FIO $\As$ are determined by those of the order zero component $\As_{[0]}$ of $\As$, consequent on the Duhamel principle combined with the nilpotence of the terms  $\As_{[k]}$ with $k>0$. The square root $\sqrt{\Qs}$ is defined via holomorphic functional calculus.  \\

The vertical trace $\Trmb$ is integration along the fibre diagonal of the fibrewise Schwartz kernels of the family of FIOs. For this, it is assumed that the vertical tangent bundle $T_\tmb$ along the fibres of $\pi:M\to B$ is oriented, that the fibres of $\pi$ are smoothly oriented manifolds. There is then a push-forward Gysin map in  cohomology $\pi_! : H^r(M, R)\to H^{r-q}(B, R)$  over a ring $R$ with $q$ the dimension of the fibre, which for $R=\R$ is realised at the form level by integration over the fibre 
\begin{equation}\label{iof form}
    \int_{M/B} : \Omega^r(M)\to \Omega^{r-q}(B)
\end{equation}
characterised by the  identity
\begin{equation}\label{int fibre identity}
    \int_B \l(\int_{M/B} \alpha\r) \wedge \omega = \int_M \alpha\wedge\pi^*\omega
\end{equation}
and by $d_B\int_{M/B} \alpha = \int_M d_M \pi^*\alpha$, corresponding to  $\pi_! = \int_{M/B}$ in real cohomology. Fibre integration $\int_{M/B}$ provides the basic tool for vertical Schwartz distribution theory in which the vertical Schwartz kernel theorem identifies a smooth family of FIOs $\As$ with  corresponding vertical oscillatory integral kernels $K_b(x,y)$ with $(x,y)\in M_b\x M_b$. If the singular support of the kernels is disjoint from the fibrewise diagonal then one has a vertical trace map defining a mixed degree differential form on $B$, via
 $$\Trmb(\As) = \int_{M/B} K_{|{\rm diag}}, \hskip 10mm \Trmb(\As)(b) = \int_{M_b} K_b(x,x).$$
$K$ has coefficients in differential forms on $M$, 
 giving meaning to the integrals. \\

The operator $\Qs$ in \thmref{thm1} is required to be of Laplace-type in the weak sense that its form degree zero component is pointwise a classical generalised Laplace operator ($\sqrt{\Qs}$ could be replaced by $\Ps^{1/m}$ for any positive self-adjoint elliptic $\pdo$ $\Ps^{1/m}_{[0]}$ of order $m$ as in \cite{DG}), while the coefficients of the components of $\Qs$ of non-zero form degree  can be essentially arbitrary positive integer order classical $\pdo$s. The vertical FIO $e^{it\sqrt{\Qs}}$ satisfies a corresponding wave equation but in general not one of geometric origin. There is, however, a class of natural geometric $\Qs$ which arise in this more general framework as the curvature $\Bf^2$ of a choice of Quillen-Bismut superconnection $\Bf$ on the fibration acting on $\Z_2$-graded spaces of sections, with $\Bf^2_{[0]} = \Ds^2$ a smooth family of Dirac Laplacians. It is convenient to rescale this to $\Bf_t$ by multiplying $\Bf_{[d]}$ by $t^{-d+1}$.
\begin{thm}\label{thm2}
Let $\Strmb$ be the vertical super trace ($\mathbb{Z}_2$-graded integration over the fibre diagonal) on a fibration $\pi:M\to B$ of Riemannian spin manifolds. Then the asymptotic singularity expansion of the distributional wave trace  
\begin{eqnarray}\label{wave trace chern character}
\Strmb\l(e^{-i\,{\rm sgn}(t)\sqrt{\Bf^2_t}}\r) & \sim &
\sum_{d=0}^{\dim B}\sum_{j\geq 0} \ww_{-q+j,d} \, (t+i0)^{-q + j - d}  \notag \\[2mm]
  &+&   \sum_{d=0}^{\dim B}\sum_{\kappa\neq 0} \sum_{j\geq 0}     \ww_{\kappa,j,d}(t-(L_\kappa+i0))^{-q+j -d}\log(t-(L_\kappa +i0))  \nonumber
 \end{eqnarray}
 is such that:
 \begin{itemize}
     \item The vertical wave trace invariants  $\ww_{-q+j,d}$ at $t=0$ are closed differential forms, representing cohomology classes in $H^d(B)$. 
     \item If $\Bf$ is chosen to be the Bismut superconnection then 
     \begin{equation}\label{sconn wave form = A-hat genuus}
         \frac{1}{2}\ww_{d,d} = \l((2\pi)^{-q/2}\int_{M/B} \widehat{A}(R^{\mbox{\tiny{M/B}}})\r)_{[d]},
     \end{equation}
     where the right-hand side is the $d$-form component of the vertical $\widehat{A}$-genus form. For any test form $\phi$, the limit $\lim_{t\to 0} (e^{-i\,{\rm sgn}(t)\sqrt{\Bf^2_t}}(\phi))$ then exists in $\Aa^*(B)$.
     \item If the function $L_\kappa : B\to \R$ is constant then the coefficients $\ww_{\kappa,i,d}$ are closed differential forms.
 \end{itemize} 
\end{thm}
Generically, then, the  coefficients $\ww_{\kappa,i,d}$ are not closed forms, not cohomological.
The $\widehat{A}$-genus $\int_{M/B} \widehat{A}(R^{\mbox{\tiny{M/B}}}) \in \Aa^*(B)$ 
is obtained by integrating over the fibre the mixed degree local families index form $\widehat{A}(R^{\mbox{\tiny{M/B}}}) = {\rm det}\left(\frac{R^{\mbox{\tiny{M/B}}}/2}{\sinh(R^{\mbox{\tiny{M/B}}}/2)}\right)^{1/2} \in \Aa^*(M)$ with $R^{\mbox{\tiny{M/B}}}$ the curvature of the induced metric connection on the tangent bundle along the fibres \cite{B,BGV}. \\

Comments on  \thmref{thm1} and classical Poisson summation for families are provided in the Appendix. 

\subsubsection*{Acknowledgements:} I am grateful to Alan Greenleaf and Richard Melrose for helpful comments, to King's College London for research leave, and to the Departamento de Matemáticas, Universidad de los Andes, Bogotá, for their hospitality. 

\section{Vertical Lagrangian distributions and FIOs}

The vertical structures, along the fibres, needed to discuss smooth families of FIOs on a fibration $\pi: M\to B$ are outlined in this Section. 

\subsection{Vertical distributions}

Let $\pi:M\to B$ be a fibration of smooth manifolds and let $E$ be a vector bundle over $M$. Consider the formal infinite-dimensional vector bundle $\pi_*(E) \to B$ whose fibre at $b\in B$ is the space of smooth sections 
$\Ci(M_b, E_b \otimes |\wedge_b|^{1/2})$
over the fibre $M_b = \pi\ii(b)$ of $M$ at $b\in B$, where $|\wedge_b|^{1/2}$ is the 1/2-density bundle along  $M_b$ and $E_b := E_{|M_b} = i_b^*(E)$ is the restriction of $E$ to the fibre for the inclusion map $i_b:M_b\hookrightarrow M$. In order for a continuous linear operator 
\begin{equation}\label{AonB}
    A: \Ci(B,\pi_*(E)) \to \Ci(B,\pi_*(E))
\end{equation} 
to correspond to a smooth family of such operators $$A_b : \Ci(M_b, E_b \otimes |\wedge_b|^{1/2}) \to \Ci(M_b, E_b \otimes |\wedge_b|^{1/2})$$ along the fibres of $\pi$ it must act as a tensor:  for $s\in \Ci(B,\pi_*(E))$ and $f\in \Ci(B)$
\begin{equation}\label{Atensor}
    A(fs) = f As.
\end{equation}
If \eqref{Atensor} holds, then  $A\in \Ci(B, \End(\,\pi_*(E)\,)).$
This allows $A$ to act non-tensorially along fibre directions of $M$ (for example as a differential operator along the fibres) but not in horizontal $B$ directions. Rigorously, any use of infinite dimensional bundles is avoided by identifying $\Ci(B,\pi_*(E))$ with the space of sections of the finite-rank bundle $E\otimes |\wedge_\pi|^{1/2}$ over $M$, where $|\wedge_\pi|^{1/2}$ is the bundle of 1/2-densities along the fibres. The operator in \eqref{AonB} is then replaced by an operator 
\begin{equation}\label{AonM}
    A:\Ci_0(M, E\otimes |\wedge_\pi|^{1/2}) \to  \Ci(M, E\otimes |\wedge_\pi|^{1/2})
\end{equation}
and \eqref{Atensor} by the equivariance requirement for $f\in\Ci(B)$, $\phi\in \Ci_0(M, E\otimes |\wedge_\pi|^{1/2})$
\begin{equation}\label{Averticaltensor}
    A((\pi^*f) \phi) = (\pi^*f) A\phi.
\end{equation} 
 The analysis of such operators takes place within the space of vertical distributions on $\pi$, an element of which is a continuous linear map
 \begin{equation}\label{Averticaldistrn}
\begin{tikzcd}
\Ci_0(M, E\otimes |\wedge_\pi|^{1/2}) \ \   \arrow[d,  "\nu_{\mbox{\tiny{M/B}}}"' ]   \\
\Ci(B)
\end{tikzcd}
\end{equation}
 which is horizontally equivariant: for $f\in\Ci(B)$ 
 \begin{equation}\label{verticalequiv}
    \pi^*f\cdot \nu_{\mbox{\tiny{M/B}}} = f \nu_{\mbox{\tiny{M/B}}},
\end{equation} 
where 
$\cdot$ is  distributional  whilst the multiplication on the right it is in the $\Ci$ sense; thus, for a test section $\phi$ on $M$ \eqref{verticalequiv} requires
\begin{equation}\label{verticalequiv2}
     \nu_{\mbox{\tiny{M/B}}}((\pi^*f)\,\phi) = f\,\nu_{\mbox{\tiny{M/B}}}(\phi)
\end{equation}
in $\Ci(B)$. Denote the space of such distributions by  $\Dd^{\,\prime}(M/B, E)$. There is a canonical inclusion map\footnote{\eqref{caninc} can more generically be considered as a map from $C^0(M, E^* \otimes |\wedge_\pi|^s)$ with $s\in [0,1]$  to vertical distributions \eqref{Averticaldistrn} on $\Ci_0(M, E\otimes |\wedge_\pi|^{1-s}).$} defined by integration over the fibre
\begin{equation}\label{caninc}
   C^0(M, E^* \otimes |\wedge_\pi|^{1/2}) \hookrightarrow \Dd^{\,\prime}(M/B, E),
\end{equation}
$$\psi\longmapsto \nu_{\mbox{\tiny{M/B}}}^\psi, \hskip 10mm \nu_{\mbox{\tiny{M/B}}}^\psi(\phi)  = \int_{\mbox{\tiny{M/B}}}\psi(\phi), $$
where $\psi(\phi)$ is the dual pairing (the fibrewise trace)
$$(E^*\otimes |\wedge_\pi|^{1/2}  )\otimes (E\otimes |\wedge_\pi|^{1/2})\to   |\wedge_\pi|.$$ 
\eqref{caninc} is injective and dense, and satisfies \eqref{verticalequiv}, \eqref{verticalequiv2}
$$\nu_{\mbox{\tiny{M/B}}}^\psi((\pi^*f)\,\phi)(b)   =  \int_{M_b}\psi(f(b)\phi(m))  =  f(b) \int_{M_b}\psi(\phi(m)) 
    = \l(f\nu_{\mbox{\tiny{M/B}}}^\psi(\phi)\r)(b).$$
Matters extend without difficulty to form valued distributions and operators.  Let $\Aa(B)$ be the de Rham algebra of differential forms on the manifold $B$ and set 
\begin{equation}\label{Aak}
    \Aa^k(M/B, E) := \Ci_0(M, \pi^*(\wedge^k T^*_B)  \otimes E\otimes |\wedge_\pi|^{1/2} )
\end{equation}
and $\Aa(M/B, E) = \bigoplus_k\Aa^k(B, E)$; the notation $ \Aa^*(B, \pi_* (E))$ may also be used. A vertical form-valued distribution, or vertical current, is a continuous linear map
\begin{equation}\label{Form Averticaldistrn}
\begin{tikzcd}
 \Aa(M/B, E):= \Ci_0(M, \pi^*(\wedge T^*_B) \otimes E\otimes |\wedge_\pi|^{1/2})   \arrow[d,  "\mu_{\mbox{\tiny{M/B}}}"' ]   \\
 \Aa(B),
\end{tikzcd}
\end{equation}
 which is horizontally equivariant: for $\omega\in\Aa(B)$ 
 \begin{equation}\label{formverticalequiv}
    \pi^*\omega\cdot \mu_{\mbox{\tiny{M/B}}} =  (-1)^{|\omega|} \omega\, \mu_{\mbox{\tiny{M/B}}}.
\end{equation}
Multiplication on de Rham algebras is wedge product and $\cdot$ is distributional; so, for a form $\omega \in\Aa(B)$ and a test form $\beta$ on $M$ \eqref{formverticalequiv} requires
\begin{equation}\label{formverticalequiv2}
     \mu_{\mbox{\tiny{M/B}}}(\pi^*\omega\wedge \beta) =  (-1)^{|\omega|}\omega\wedge\mu_{\mbox{\tiny{M/B}}}(\beta)
\end{equation}
in $\Aa(B)$. We denote the space of such distributions by $\Af^{\,\prime}(M/B,E)$. 
The dense inclusion
\begin{equation}\label{formcaninc}
   \Aa(M/B, \pi_* (E)) = C_0^0(M, \pi^*(\wedge T^*_\tb) \otimes E^* \otimes |\wedge_\pi|^{1/2}) \hookrightarrow \Af^{\,\prime}(M/B,E)
\end{equation}
with $\pi^*\beta\otimes\theta\otimes\delta^{1/2}_\pi\longmapsto \nu_{\mbox{\tiny{M/B}}}^{\pi^*\beta\otimes\theta\otimes\delta^{1/2}_\pi}$ is given by
\begin{eqnarray}
    \nu_{\mbox{\tiny{M/B}}}^{\pi^*\beta\otimes\theta\otimes\delta^{1/2}_\pi}(\pi^*\alpha\otimes \phi\otimes\tau^{1/2}_\pi) & =&  \int_{M/B}\pi^*(\beta\wedge\alpha)  
 \theta(\phi)\,\delta^{1/2}_\pi\otimes\tau^{1/2}_\pi \notag \\
    & =&  \beta\wedge\alpha\cdot\int_{M/B}\theta(\phi)\,\delta^{1/2}_\pi\otimes\tau^{1/2}_\pi \label{form caninc}
\end{eqnarray}
A vertical current $\mu_{\mbox{\tiny{M/B}}}\in \Af^{\,\prime}(M/B,E)$ ranging in $\Aa(B)$ composed with a classical current  $\mu_B\in\Af^\prime(B)$ 
yields a classical current 
 \begin{equation}\label{de Rham caninc}
     \upsilon_M  = \mu_B\circ \mu_{\mbox{\tiny{M/B}}} \in\Af^\prime(M).
 \end{equation}
 For the fundamental case \eqref{formcaninc} this is the property \eqref{int fibre identity} of fibre integration:  $\nu_{\mbox{\tiny{M/B}}}^{\pi^*\beta\otimes\theta\otimes\delta^{1/2}_\pi}$ composes with $\mu_B = \int_B$ to give $\upsilon_M = \int_M$, that is,
$$ \int_B  \beta\wedge\alpha\cdot\int_{M/B}
\theta(\phi)\,\delta^{1/2}_\pi\otimes\tau^{1/2}_\pi  
= \int_M  \pi^*(\beta\wedge\alpha)
\theta(\phi)\,\delta^{1/2}_\pi\otimes\tau^{1/2}_\pi. $$
This is \eqref{int fibre identity} adapted to densities rather than forms, via the line bundle isomorphism  (of trivialisable line bundles) $|\wedge_\pi| \cong \wedge^{{\rm max}}T^*_{\mbox{\tiny{M/B}}}$ (since $\pi$ has orientable fibres). \\

{\small \noindent Comment: {\it Indeed,  $|\wedge_\pi|$ can be replaced throughout with the vertical exterior cotangent bundle $\wedge T^*_\tmb$. With a choice of vector bundle isomorphism 
\begin{equation}\label{cotangent splitting}
    T^*_{\mbox{\tiny{M}}} \cong \pi^*(T^*_{\mbox{\tiny{B}}}) + T^*_{\mbox{\tiny{M/B}}},
\end{equation}
and hence $\wedge T^*_\tm \cong \pi^*(\wedge T^*_\tb) \otimes \wedge T^*_\tmb$ for $\wedge := \oplus_{k\geq 0} \wedge^k$, then
$\Aa^*(M/B, E)$ in \eqref{Aak} is replaced by
$\Aa^*(M,E) := \Ci(M, \wedge T^*_\tm \otimes E) =  \Ci(M, \pi^*(\wedge^k T^*_B)  \otimes E\otimes \wedge T^*_\tmb)$ 
and development proceeds in a parallel manner.
This may be more natural and canonical in some geometric applications. For $\omega\in\Aa(B)$  and  $\alpha\in\Aa(M)$ integration over the fibre 
\begin{equation}\label{form int over fibre}
 \int_\tmb (\pi^*\omega \wedge \alpha)  =  \omega \wedge \int_\tmb \alpha   
\end{equation}
 is a linear map as in \eqref{iof form} from $r$ forms on $M$ to $r-q$ forms on $B$, whereas on $\Aa^*(M/B, E)$  \eqref{Aak} $q$-forms along the fibre are identified with densities along the (oriented) fibre and  integration over the fibre in this convention preserves form degree
$$\int_\tmb: \Aa^r(M/B) \to \Aa^r(B)$$ satisying \eqref{form int over fibre}.}}\\

The vertical Schwartz kernel theorem states that a continuous linear map 
\begin{equation}\label{FSKT4}
\As:\Ci_0(N, \sigma^*(\wedge T^*_\tb)\otimes F \otimes |\wedge_\sigma|^{1/2}) \to \Af^{\,\prime}(M/B,E)
\end{equation}
which commutes with the action of $\Aa(B)$ is equivalent to a vertical distribution 
\begin{equation}\label{Fvertkernel}
\Ks_\As \in \Af^{\,\prime}((M\times_B N)/B, E^*\boxtimes F), 
\end{equation}
that is, to a horizontally equivariant map
$$\Ks_\As: \Ci_0(M\times_B N, \rho^*(\wedge T^*_\tb) \otimes (E\otimes |\wedge_\pi|^{1/2}) \boxtimes (F^*\otimes|\wedge_\sigma|^{1/2}))) \to \Aa(B),$$
via the correspondence $(\As\psi)(\phi) = \Ks_\As(\phi\otimes\psi)$.\\

$\As$ commutes with the $\Aa(B)$ action  means for $\psi\in\Ci_0(N, \sigma^*(\wedge T^*_\tb)\otimes F \otimes |\wedge_\sigma|^{1/2})$  
that $(\sigma^*\beta)\cdot\As_{[k]}(\psi) = (-1)^{k|\beta|}\beta\As_{[k]}(\psi) $ in $\Af^{\,\prime}(M/B,E)$, 
i.e.
$\As(\sigma^*\beta \wedge \psi) = (-1)^{k|\beta|}\beta \wedge \As(\psi)$.\\

$\As$  is thus equivalent to giving a vertical distribution $k_\As$ on the fibre product fibration $\rho:M\x_B N \to B$ with fibre $M_b\x N_b$.  Here $M\x_B N = (\pi_M\x \pi_N)\ii({\rm diag}(B))$ for ${\rm diag}(B)$ the diagonal in $B\x B$; since $\pi_m$ and $\pi_N$ are submersions  $\pi_M\x \pi_N:  M\x N \to B\x B$ is transverse to ${\rm diag}(B)$ and so  $M\x_B N$  is a smooth manifold of dimension
$\dim M + \dim N - \dim B$. Vector bundles $\Vv\to M$, $\Ww\to N$ determine a vector bundle $\Vv\boxtimes\Ww :=
p_1^*(\Vv)\otimes p_2^*(\Ww)$ over $M\x_X N$, where $p_1,p_2$  are the projection maps  $M\times_{\pi}N$ to $M,N$.  The vertical distribution $k_\As$ in turn corresponds to a family of classical  Schwartz kernels $k_{A_b}$ on the fibres $M_b\x N_b$, and hence to a family of operators $\{A_b \ | \ b\in B\}$ along the fibres.  The movement here between a vertical  (horizontally equivariant) distribution and a smooth family of classical distributions along the fibres of $\pi:M\to B$ is:
\begin{lem}\label{vertical distrn b}
A vertical distribution $\nu_{\mbox{\tiny{M/B}}}\in \mathbb{A}^{\,\prime}(M/B,E)$  restricts on each fibre $M_b$ of $\pi$ to a classical distribution current $\nu_b\in \mathbb{A}^{\,\prime}(M_b, E_b)$ which depends smoothly on $b\in B$ (i.e $b\mto\nu_b(\phi_b)$ is smooth differential form on $B$ for each $\phi_b\in\Aa_0(M_b, E_b)$). Conversely, a smooth family $\{\nu_b\}$  of such distributions determines uniquely a vertical distribution.
\end{lem}

The mechanics of \eqref{FSKT4} are that if $A$ ranges in $C^0(M, \pi^*(\wedge T^*_\tb)\otimes E \otimes |\wedge_\pi|^{1/2})$  then  
for $s\in \Ci_0(N, \sigma^*(\wedge T^*_\tb)\otimes F \otimes |\wedge_\sigma|^{1/2})$ and $m \in M$
\begin{eqnarray*}
    (\As s)(m) &=& \int_{N/B} \Ks_\As (m,n) s(n) \\ 
    &=& \int_{N_b} \Ks_b(m_b,n_b) s(n_b), 
\end{eqnarray*}
where $$b = \pi(m =m_b) \ \ \ {\rm and} \ \ \ n =n_b \in N_b,$$ while as a distribution $\As s\in \Af^{\,\prime}(M/B,E)$ acts on a test section  $\phi\in \Ci_0(M, \pi^*(\wedge T^*_\tb)\otimes E \otimes |\wedge_\pi|^{1/2})$ by 
\begin{eqnarray*}
    \As s(\phi) &=& \int_{M/B} \As s(m)\phi(m) \\ 
   &=& \int_{M/B} \int_{N/B} \Ks_\As (m,n) s(n) \phi(m) = \int_{(M\x_B N)/B} \Ks_\As (m,n) s(n)\phi(m) 
\end{eqnarray*}
defining the differential form on $B$  $$b\mto \As s(\phi) (b) := \int_{M_b} \int_{N_b} \Ks_\As (m_b,n_b) s(n_b) \phi(m_b).$$
For example, if $\pi:M\to B$ is a Riemannian fibration and $\Ls$ a family of differential operators with form degree zero component $\Ls_{[0]} = \{\D_{g_b}\}$ a family of Laplace-type operators defined by the Riemannian metric $g_b$ on $M_b$, then for $t>0$ the heat operator $e^{-t\As}$ is well-defined as a smooth family of smoothing operators with vertical heat kernel 
$ \Ks_{e^{-t\Ls}} \in \Ci(M\times_B M, \pi^*(\wedge T^*_\tb)\otimes (E \otimes |\wedge_\pi|^{1/2})\boxtimes (E^* \otimes |\wedge_\pi|^{1/2}) ),$ the form degree zero component coinciding pointwise with the usual heat kernel of $\D_{g_b}$. See \cite{B,BGV} and \S3 here for details of the construction. 



\subsection{Vertical Lagrangian distributions and FIOs}

A vertical (fibrewise vertically-embedded) conic Lagrangian fibration in $T^*_M\bs 0$ is a closed conic Lagrangian submanifold $\Lambda$  of $T^*_M\bs 0$ which is the total space of  a  conic subfibration $\la: \Lambda \to B$ of the cotangent bundle along the fibres $T^*_{M/B}\bs 0$ (minus the zero section). Thus,   $\la: \Lambda \to B$ has fibre 
at $b\in B$ a conic Lagrangian submanifold 
$$\la\ii(b) = \Lambda_b  \< T^*_{M_b}\bs 0.$$ 
One has
$$\dim T^*_M\bs 0 = 2q + 2\beta$$
and by the Lagrangian condition   $$\dim \Lambda = q + \beta,$$ where $$q := \dim M_b = \dim \Lambda_b =  \frac{1}{2}\dim T^*_{M_b}, \ \ \ \ \beta = \dim B.$$ 
$T^*_{M/B}\bs 0$ (a subbundle of $T^*_M\bs 0$) has dimension $2q + \beta$ and is considered here as a fibration over $B$ with fibre at $b\in B$ the cotangent bundle $T^*_{M_b}\bs 0$ of dimension $2q$ of $M_b$. $\Lambda$ is not (in general) a Lagrangian submanifold of the (in general) non-symplectic manifold $T^*_{M/B}\bs 0$, although the fibres of $\Lambda$ are Lagrangian submanifolds of the fibres of $T^*_{M/B}\bs 0$. $\Lambda$ is, on the other hand, assumed to be a Lagrangian subfibration of $T^*_M\bs 0$ when the latter is considered as a fibration over $B$ with fibre $(T^*_M\bs 0)_{|M_b} = T^*_{M_b}\bs 0 + \pi^* T^*_b B$ over $b\in B$.\\

Any vertical conic Lagrangian fibrations is parametrised locally by non-degenerate vertical phase functions.  Using local triviality of  $\pi$ and the corresponding classical fact for conic Lagrangian submanifolds it follows that a $q + \beta$ dimensional subfibration $L\< T^*_{M/B}\bs 0 \< T^*_M\bs 0$ is a a vertical Lagrangian fibration if and only if every $m\in L$ has a fibred conic neighbourhood  $\Gamma = \cup_b \Gamma_b$ such that  $\Gamma_b\cap L = \Lambda_{\phi_b}$ some nondegenerate vertical phase function $\phi$.
For such local considerations one may use fibrewise charts defined relative to a local trivialization $M_{|W} \cong W \x M_b$ of $\pi$ over a chart 
$W = (W, (b_1, \ldots, b_\beta)) \< B$ and fibre $M_b = \pi\ii(b)$ with $b\in W$ and local chart $U = (U, (x_1, \ldots, x_q)) \< M_b$. In these coordinates $m = (b,x) = (b_1, \ldots, b_\beta, x_1, \ldots, x_q)$, $\phi \in \Ci(W\x U \x \R^N\bs 0)$ is a local vertical phase function for $\pi$ if $\phi(b,x,\xi)$ is homogeneous of degree 1 in $\xi$ and for $\xi\neq 0$ and for each fixed $b_0$ 
\begin{equation}\label{fibre phase}
    d_{x,\xi}\phi(b_0,x,\xi)  \neq 0.
\end{equation}
 Thus, $\phi$ is in the classical sense a local phase function on each fibre $M_b$ with $N$ phase variables which as $b$ varies  extends by default to a local phase function $M$, free of critical points, parametrising locally $\Lambda$.  $\phi$ is vertically non-degenerate on a fibred cone $\Gamma \< W\x U \x\R^N$ if $\Gamma \to W$ is fibred by cones $\Gamma_b < \{b\}\x U \x\R^N$ such that for each fixed $b_0$ setting $\phi_{b_0}(x,\xi) = \phi(b_0,x,\xi)$ the differentials $d_{x,\xi}\l(\frac{\partial \phi_{b_0}}{\partial \xi_1}\r), \ldots,d_{x,\xi}\l(\frac{\partial \phi_{b_0}}{\partial \xi}\r)$ are linearly independent on 
$$C_\phi^{M_{b_0}} = \{ (b_0,x,\xi) \in \Gamma_{b_0} \ | \ d_\xi \phi_{b_0}(x,\xi) =0\}.$$
It follows from the Implicit Function theorem that these form a closed conic submanifold $C_\phi^M$ of $W\x U \x\R^N$ of dimension $q+\beta$ fibred by closed conic submanifolds $C_\phi^{M_{b_0}}$ of $\{b_0\}\x U \x\R^N$ of dimension $q$. The fibrewise non-degeneracy is equivalently characterised as the non-degeneracy of the Hessians  $H\phi_{b_0}$ on the normal bundles to  $C_\phi^{M_{b_0}}$.  In this case, it follows from standard symplectic geometry that the maps
\begin{equation}\label{hphiMb}
    h_{\phi_{b_0}} : C_\phi^{M_{b_0}}  \to T^*_{M_{b_0}}\bs 0, \ \ \ (b_0, x,\xi) \mapsto (x,d_x\phi_{b_0}), 
\end{equation}
are respectively immersions onto conic Lagrangian submanifolds $\Lambda_{\phi_{b_0}} \< T^*_{M_{b_0}}\bs 0.$ $\phi$ is not in general non-degenerate in horizontal directions, and so not in general a non-degenerate phase function for $\Lambda$.\\

An important class of vertical conic Lagrangians occur as the  conormal bundle along the fibres  $\Lambda = i^*(T^*_{M/B}\bs 0) / (T^*_{Z/B}\bs 0)$  of a smooth subfibration $i: Z\hookrightarrow M$  of $\pi: M\to B$. $\Lambda$ has fibre at $b\in B$ the conormal bundle $$\Lambda_b  =  (i_b^* T^*_{M_b}\bs 0) / (T^*_{Z_b}\bs 0) \to Z_b$$  to the submanifold $i_b: Z_b \hookrightarrow M_b$ of the fibre $M_b$ of $\pi$. $\Lambda_b$ defines a Lagrangian submanifold of $T^*_{M_b}\bs 0$. On the other hand, there is a vector bundle isomorphism over $M$ of $\Lambda$ with the conormal bundle of $Z$ as a smooth submanifold of $M$ - for, from \eqref{cotangent splitting} 
\begin{equation}\label{L=fibreL}
 \frac{i^*(T^*_M\bs 0)}{T^*_Z\bs 0} = \frac{ i^*T^*_{\mbox{\tiny{M/B}}}\bs 0 + i^* \pi^* T^*_{\mbox{\tiny{B}}}\bs 0}{T^*_{\mbox{\tiny{Z/B}}}\bs 0 + (\pi_{|Z})^* T^*_{\mbox{\tiny{B}}}\bs 0} 
= \frac{ i^*T^*_{\mbox{\tiny{M/B}}}\bs 0}{T^*_{\mbox{\tiny{Z/B}}}\bs 0} 
\end{equation}
since $\pi_{|Z}:= \pi\circ i$. Thus, $\Lambda \cong i^* (T^*_M\bs 0) / (T^*_Z\bs 0)\to Z$ 
 is diffeomorphic to a Lagrangian submanifold of $T^*_M\bs 0$. Locally, $\Lambda$ is parametrised by a vertical phase function which is linear with respect to the $\xi$ variables, so $ \phi(b,x,\xi) = \langle \Phi(b,x), \xi\rangle$
 with $\Phi: W\x U\< M\to \R^N$. Since $\phi$ for $b=b_0\in W$ fixed must be a non-degenerate phase function in the classical sense, then  $\langle \Phi^\prime_{x_j} (b_0,x), \xi\rangle \neq 0$ for some $j$ when $d_\theta\phi_{b_0} = \Phi(b_0, \ )=0$, and hence $N  \leq q $ and $\Phi^\prime_x(b_0, \ )$ has rank $N$ at those $x\in M_{b_0}$ where  $\Phi(b_0,x)=0$, this equation therefore defining a submanifold 
$$Z_{b_0}= \{ (b_0,x) \in M_{b_0} \ | \ \phi_1(b_0,x) = \cdots = \phi_N(b_0,x) = 0\}$$ of $M_b$ of codimension $N\leq q$, where $\Phi=(\phi_1, \ldots, \phi_N)$. Next as we allow $b=b_0$ to vary within $W\<B$ the resulting fibre bundle $Z$ still has the same codimension $N$ in $M$ (since $M_b$ is also varying with $b\in W$). Thus, for $\phi(b,x,\xi)$ to define also a phase function on $M$, $\Phi^\prime_{b,x}$ also has rank $N$ at those $(b,x)$ in $M$ where  $\Phi(b,x)=0$, so we require additionally to maintain verticality that $\Phi^\prime_{b}$ has rank $0$ at those $(b,x)$ in $M$ where  $\Phi(b,x)=0$, so that
 $$Z = \{ (b,x) \in M \ | \ \phi_1(b,x) = \cdots = \phi_k(b,x) = 0, \ {\rm rk}(\Phi^\prime_{b}(b,x))=0\}$$
 fibres over $B$ with fibre $Z_{b_0}$. 
 The latter conditions then ensure that $\phi(b,x,\xi) = \sum_j \xi_j\,\phi_j(b,x)$ is a local linear phase function parametrising each Lagrangian fibre $\Lambda_b$ as
$$\Lambda_{\phi_b} = \{(b,x,\theta^\prime)\in T^*_{M_b}\bs 0  \ | \ (b,x)\in Z_b, \,\theta^\prime = \sum_{j=1}^{N} \xi_j\,d_x\phi_j(b,x) \}.$$ 
To additionally   be non-degenerate for the conormal bundle $\Lambda$ requires that ${\rm rk}(\Phi^\prime_{b}(b,x))=0$ at those $(b,x)$ for which $\Phi(b,x)=0$ giving $\Lambda_\phi =\cup_b \Lambda_{\phi_b}$ corresponding to the equality \eqref{L=fibreL}. This holds for the basic case where $\Phi(b,x)$ is the local quadratic form $\sum_{i,j} g^{i,j}(b) b_i x_j$ with $g^{i,j}$ the dual Riemannian metric on $M_b$. \\[3mm]

If $\la: \Lambda \to B$ is a vertical conic Lagrangian fibration in $T^*_M\bs 0$ and let $m\in \R$, then the space 
\begin{equation}\label{vert Lagrangian}
    \If^m(M/B, \Lambda)_{[k]}
\end{equation}
of vertical Lagrangian distributions  of (vertical) order $m\in\R$ and form degree $k\in\{0,1,\ldots, q+\beta\}$ associated to $\Lambda$ consists of those vertical distributions $u^\tmb\in\Af^{\,\prime}(M/B)$ which may be written microlocally as a finite sum of vertical Fourier integral distributions acting on a test form $\omega\in \Aa_0(M,  |\wedge_\pi|^{1/2})$ on $M$ by
\begin{equation}\label{vert osc int}
    u_{[k]}^\tmb(\omega) =  \int_\tmb \int_{\Ww^\Lambda/M} e^{i\phi(m,\xi)} \, a(m,\xi)_{[k]}\wedge \omega(m)
\end{equation}
for  $\Ww^\Lambda \to M$ a fibration of  rank $N$ vector spaces and $\phi : \Ww^\Lambda \to \R$ a vertical non-degenerate phase function parametrising $\Lambda$, with $\xi$ an element of the fibre $ \Ww^\Lambda_m\cong\R^N$, and 
\begin{equation}\label{vert symbol}
 a_{[k]}\in \Sss^{m+\frac{q}{4} - \frac{N}{2}}(\Ww^\Lambda)_{[k]}
\end{equation}
a form-valued classical symbol of order $m+\frac{q}{4} - \frac{N}{2}$ and form degree $k$, meaning a smooth map
$$a_{[k]}: \Ww^\Lambda \to|\wedge_\pi|^{1/2}\otimes\pi^* (\wedge^k T^*_B) \ox \End E$$ 
 satisfying  estimates 
  $|\partial^\alpha_b\partial^\beta_x\partial^\gamma_\xi a_{[k]}(b,x,\xi)| < C_{\alpha,\beta,\gamma, K}(1 + |\xi|)^{m+\frac{q}{4} - \frac{N}{2} - |\gamma|}$ 
in local coordinates on a
compact subset $K$ of $M$ uniformly in $\xi$. As a map on $B$
\begin{equation}\label{vert osc int b}
    u_{[k]}^\tmb(\omega)(b_0) =  \sum_{|I|=k} \int_{M_{b_0}} \int_{\R^N} e^{i\phi_{b_0}(x,\xi)} \, a^I_{b_0}(x,\xi)\, \omega^I_{b_0}(x)\ dx\,d\xi \ db_I(b_0)
\end{equation}
is a sum of oscillatory integrals in the classical sense but taking values in differential forms on $B$, where $i^*_b (a \wedge \omega)(b,x)  = \sum_{|I|=k} a^I_b(x,\xi)\, \omega^I_b(x)\ dx\,d\xi \,db^I(b_0)$ is the restriction to $M_{b_0}$, with $x\in M_{b_0}$ and $db_I\in\Aa^k(B)$ a differential $k$-form on $B$.  In view of \eqref{fibre phase} this is well-defined as an oscillatory integral. A general element of $$ \If^m(M/B, \Lambda) = \bigoplus_k\If^m(M/B, \Lambda)_{[k]} $$ can be written as a finite sum of terms of the form $ u^\tmb =\sum_k u_{[k]}^\tmb$ plus a smooth family of smoothing operators. Equivalently, 
a vertical distribution is said to belong to $ \If^m(M/B, \Lambda)$ if for every point $m=(b,x)\in M$ there exists a fibrewise coordinate patch $W\x U$, centred at $m$, such that $u$ is equal to a vertical oscillatory integral \eqref{fibre phase} (with integral over $M_b$ restricted to $U$). It is easy to see that a fibrewise (bundle) diffeomorphism transforms such local oscillatory integrals into integrals of the same type.  \\

The dimension $N$ in \eqref{vert osc int b} may be different for different summands for a generic element of $\If^m(M/B, \Lambda)$, though phase variables can be added in the usual way to bring all local representatives up to same rank without changing, up to equivalence, the corresponding Fourier integral. Conversely, the number $N$ of vertical phase variables in \eqref{vert osc int b} can always be reduced to $k:= \dim T(\Lambda_b) \cap T_{M_b}$ using vertical stationary phase approximation (v-SPA), with an equivalent phase function $\tilde\phi(x, \tilde\xi)$ and symbol 
$\tilde a_{[k]}\in \Sss^{m+\frac{q}{4} - \frac{k}{2}}$ - equivalent insofar as the corresponding Fourier integral  defines the same vertical distribution as \eqref{vert osc int b}.
For v-SPA, consider a fibration $\rho:Z\to B$ and a vertical oscillatory integral 
$$\int_{Z/B} e^{ir\phi(z)} q(z)\omega,$$
where $r$ is a real variable, $q$ is a compactly supported section of $\pi^*(\wedge^* T^*_B)$, $\omega \in \Ci(Z, |\wedge^{{\rm max}}T^*_{Z/B} |\bs 0)$ a vertical density (along the fibres), and where $\phi : Z \to \R$ is a vertical Bott-Morse function, meaning that $\phi$ restricts on each fibre $Z_b$ to a Bott-Morse function $\phi_b (y) := \phi(b,y)$ (in the usual sense) such that the union of the sets of critical points form a fibration over $B$. More precisely, one requires that:  
\begin{enumerate}
    \item The fibrewise critical set  
    $${\rm Crit}_{\phi_b} = \{z=(b,y)\in M_b \ |\ d\phi_b = d_y\phi(b,y)=0 \}$$
is a submanifold of $M_b$  for each $b\in B$, forming in the subspace topology a fibration $$\rho_\phi :  \Ww := \cup_b {\rm Crit}_{\phi_b}\to B$$ of $\rho$, so $\rho_\phi = \rho_{|\Ww}$, comprising a possible disjoint union of connected subfibrations $\rho_{\phi, i} : \Ww_i\to B$  of $\rho_\phi$, so $\Ww = \bigsqcup_i \Ww_i$ with $\Ww_i$   a connected submanifold of $M$ but which may have disconnected fibres $\Ww_{i,b} = \bigsqcup_k W_{i,b,k}$. More invariantly, avoiding local fibrewise coordinates, a choice of of Riemannian metric on $B$ and a choice of metric on $T_{Z/B}$ determine an orthogonal projection map $P_{Z/B}: T^*Z \to T^*_{Z/B}$, then
$$ \Ww = \{z\in M \ |\ P_z(d_z\phi) =0 \}$$
with $P_z$ the action of $P_{Z/B}$ on $T^*_z Z$ and $ \Ww_b = \{z\in M_b\ |\ P_z(d_z\phi) =0 \}$.
     \item For each $m\in \Ww$ the fibrewise Hessian $h\phi_b = d^2\phi_b$ is non-degenerate on the normal bundle  to the fibres $N_{Z,\Ww} \to \Ww$ with respect to $T_{Z/B}$; that is, where
     \begin{equation}\label{normaltoWw}
       (T_{Z/B})_{|\Ww} = T_{\Ww/B}+ N_{Z,\Ww}.
     \end{equation}
\end{enumerate}
It follows that $\phi_b$ is a constant $\phi_b(W_{i,b,k})$ on each connected fibre $W_{i,b,k}$. The v-SPA formula is as follows.\\

\begin{thm}\label{vspa0} As $r\to \infty$ there is an asymptotic formula  
    \begin{equation}\label{v-SPA0}
    \int_{Z/B} e^{ir\phi(z)} q(z)\omega \sim \sum_i \left(\frac{2\pi}{r}\right)^{\frac{n_i}{2}}\, e^{i\frac{\pi}{4}{\rm sgn}(\Ww_i)}\, e^{ir\phi(\Ww_i)} \sum_{j\geq 0} \int_{\Ww_i/B} \alpha_{j,i} \ r^{-j}
\end{equation}
where $n_i$ is the (fibre) codimension of $\Ww_i$ in $Z$, ${\rm sgn}(\Ww_i)$ is the fibrewise signature of the Hessian on $\Ww_i$, and where
\begin{equation}\label{alpha0+1}
    \alpha_{0,i}(z) = q(z)|\det(h_z\phi_b)|^{-1/2}\omega
\end{equation}
while  $\alpha_{j,i}(z)$ is a vertical derivative, along the fibre, of this density. 
\end{thm}
 {\it Proof:}
The non-degeneracy of the Hessian $H\phi_b$ means that $|\det(H\phi_b)|^{1/2}$ is a trivialising section of the  line bundle $|\wedge^{{\rm max}}N_{Z,\Ww} |^{1/2}$ and hence from \eqref{normaltoWw} that the quotient 
$$|\det(H\phi_b)|^{-1/2} \omega \in \Ci(\Ww, |\wedge_{\rho_\phi}|)$$
is a vertical density, which can be integrated along the fibres of $\rho_\phi:\Ww\to B$, and likewise on each connected subfibration $\Ww_i\to B$, so the integrals $\int_{\Ww_i/B} $ are well-defined - given \eqref{alpha0+1}. The expressions on either side of \eqref{v-SPA0} are forms on the base manifold $B$ and evaluating at $b\in B$  \eqref{v-SPA0} is the statement that 
\begin{equation*}
    \int_{Z_b} e^{ir\phi_b(z)} q(z)\omega_b \sim \sum_{i,k} \left(\frac{2\pi}{r}\right)^{\frac{n_{i,b,k}}{2}}\, e^{i\frac{\pi}{4}{\rm sgn}((H\phi_b)_{| W_{i,b,k}})}\, e^{ir\phi( W_{i,b,k})} \sum_{j\geq 0} \int_{W_{i,b,k}} (\alpha_{j,i})_{| W_{i,b,k}} \ r^{-j},
\end{equation*}
where $n_{i,b,k} $ is the fibre codimension of $W_{i,b,k}$.
But this is just the usual SPA formula for manifolds of \cite{C}, \cite{deV}, differing only here in taking values in the de Rham algebra $\Aa^*(B)$ on $B$. Thus, \eqref{v-SPA0} holds.  \hfill $\Box$ \\[3mm]

For example, if ${\rm Crit}_{\phi_b}$ is a discrete subset of the fibre $Z_b$ the fibration $\Ww\to B$ is a covering  and the sum of integrals over the $W_{i,b,k}$ is a discrete sum over the  critical points of $\phi_b$. An interesting example is to take a fibrewise Bott-Morse function on the tautological line bundle (the open Mobius strip) fibred over $\R P^1$.\\

A vertical (horizontally equivariant)  Fourier integral operator $$\As:\Aa_0(N, F) \to \Af^{\,\prime}(M/B, E)$$ of order $m\in\R$  between fibrations $\sigma: N\to X$ and   $\pi:M\to B$ 
is one whose vertical distribution Schwartz kernel 
 $\Ks_\As \in \Af^{\,\prime}((M\times_B N)/B, E^*\boxtimes F)$  is in the space 
\begin{equation}\label{vert FIO Lagrangian}
    \If^m((M\times_B N)/B, E^*\boxtimes F, \Lambda^\tmnb)
\end{equation}
of vertical Lagrangian distributions defined with respect to a vertical conic Lagrangian subfibration  
$\la: \Lambda^\tmnb \to B$ of the cotangent bundle $T^*_{(M\x_BN)/ B}\bs 0$  along the fibres of the fibre product $M\x_B N$ (where $T^*_{(M\x_BN)/ B}\bs 0$ is considered as a fibration over $B$ with fibre the cotangent bundle $T^*_{M_b\x N_b}\bs 0$ over $M_b\x N_b$ at $b\in B$), so the fibre of $\Lambda^\tmnb$ is a Lagrangian submanifold 
 $$\Lambda_b  \< T^*_{M_b\x N_b}\bs 0.$$ 
 of dimension $\dim(M_b) + \dim(N_b)$.\\
 
Thus, $\Ks_\As$ defines a continuous bilinear form from $\Aa_0(M, E) \x\Aa_0(N, F)$ to $\Aa(B)$ and so a continuous linear map from $\Aa_0(N, F)$ to continuous linear forms on $\Aa_0(M, E)$ with values in $\Aa(B)$, and requiring that $\Ks_\As$ is an element of \eqref{vert FIO Lagrangian} means that locally 
it has $k$-form component which is a finite sum of terms of the form
\begin{equation}\label{vert kernel}
    \Ks_\As (m,n)_{[k]} =   \int_{\R^p} e^{i\phi_b(m,n,\xi)} \, a(m,n,\xi)_{[k]}
\end{equation}
for 
$$(m,n) = (b,x,y) \in M_b \x N_b,$$
  $\phi : \R^p\to \R$ a vertical non-degenerate phase function parametrising $\Lambda_b^\tmnb$  and 
\begin{equation*}
 a_{[k]}\in \Sss^{m+\frac{q}{4} - \frac{N}{2}}(B\x M_b \x N_b\x \R^p, \Hom(F,E))_{[k]}
\end{equation*}
a classical symbol of order $m+\frac{q}{4} - \frac{N}{2}$ of form degree $k$, defining 
\begin{equation*}
    \As_{[k]}(u) =  \int_\tmnb \int_{\Ww/(M\times_B N)} e^{i\phi(m,n,\xi)} \, a(m,n,\xi)_{[k]} u(n)
\end{equation*}
for $\Ww \to M\times_B N$ a fibration of  rank $p$ vector spaces and $\phi : \Ww \to \R$ a vertical non-degenerate phase function parametrising $\Lambda$, with $\xi\in\Ww_{(m,n)}\cong\R^p$. \\ 

To discuss compositions we refine matters to vertical FIOs associated to a (vertical) canonical relation. For this one needs the following identifications for fibrations. 
\begin{prop}\label{tmxn}
There are fibrewise-symplectomorphisms between the following fibrations over $B$ with fibre the cotangent bundles along the corresponding fibres:
\begin{equation}\label{fibre symp isom}
T^*_\tmnb = \mu^*T^*_\tmb + \nu^*T^*_\tnb  = T^*_\tmb \x_B T^*_\tnb,  
\end{equation}
where $\mu: M\x_BN\to M$, $\nu: M\x_BN\to N$ are the projections maps. The first identification is a vector bundle isomorphism (over $M\x_BN$), the second a  fibrewise symplectomorphism of fibrations over $B$. 
\end{prop}
The proof is contained in the Appendix.  It follows that  a vertical conic subfibration  
$\la: \Lambda^\tmnb \to B$ of the cotangent bundle $T^*_{(M\x_BN)/ B}\bs 0$ which is vertically Lagrangian for $\omega_\tmb + \omega_\tnb$, so the fibre of $\Lambda^\tmnb$ is a Lagrangian submanifold $\Lambda_b  \< T^*_{M_b\x N_b}\bs 0$ for $\omega_{M_b} + \omega_{N_b}$,  determines  
a vertical conic subfibration  $\Cs^\tmnb \to B$ of the fibration $T^*_\tmb \x_B T^*_\tnb\bs 0$ which is vertically Lagrangian for $\omega_\tmb - \omega_\tnb$, so the fibre of $\Cs^\tmnb$ is a Lagrangian submanifold $\Cs_b  \< T^*_{M_b}\x T^*_{N_b}\bs 0$ for $\omega_{M_b} - \omega_{N_b}$. Such a vertical conic subfibration  $\Cs^\tmnb$ is called a vertical homogeneous canonical relation from $T^*_\tnb$ to  $T^*_\tmb$.  The space $\If^m(M/B\times_B N, \Lambda^\tmnb)$ of vertical FIOs $\As$ associated to $\Lambda^\tmnb$ (or rather to the kernel $\Ks_\As$ of $\As$)  will hence also be viewed as the class of vertical FIOs defined by by the corresponding vertical canonical relation  $\Cs^\tmnb$, and denoted $\If^m(M/B\times_B N, \Cs^\tmnb).$ \\

A vertical canonical relation $\Cs : = \Cs^\tmnb \< T^*_\tmb \x_B T^*_\tnb\bs 0$ maps subfibrations $\Es\to B$ of $T^*_\tnb$ to possibly singular subfibrations of $T^*_\tmb$, via $\Cs(\Es)_b =\{\rho\in  T^*_{M_b}\bs 0 \ |\ \exists \eta \in \Es_b \ (\rho , \eta) \in \Cs_b \} $, while vertical canonical relations $\Cs^1 \< T^*_\tmb \x_B T^*_\tnb\bs 0$ and 
$\Cs^2 \< T^*_\tnb \x_B T^*_\tpb\bs 0$ (for a third fibration $\chi:P\to B$)  can be composed to give the vertical canonical relation fibration $\Cs^1\circ \Cs^2 \< T^*_\tmb \x_B T^*_\tpb\bs 0$ whose fibre at $b\in B$ is  $\Cs^2_b\circ \Cs^1_b = \{(\alpha, \beta) \in T^*_{M_b} \x T^*_{P_b} \ | \ \exists \gamma\in T^*_{N_b}, \  (\alpha, \gamma) \in \Cs^1_b, (\gamma, \beta) \in \Cs^2_b \}$, which  if $\As_1\in \If^{m_1}(M/B\times_B N, \Cs^1)$ and  $\As_2\in \If^{m_2}((N\times_B P)/B, \Cs^2)$ sit within the vertical transverse intersection calculus is then mirrored by the composed vertical FIO 
\begin{equation}\label{comp calculus}
    \As_1\circ\As_2\in \If^{m_1 + m_2}((M\times_B P)/B, \Cs^1\circ \Cs^2).
\end{equation}
For, noting that  
$$\Cs_1\circ \Cs_2   = (\Pi_M \x \Pi_P)((\Cs^1\x\Cs^2) \cap_B T^*_\tmb \x \Delta(T^*_\tnb) \x  T^*_\tpb)$$
where $\Pi_M$ is the fibrewise projection to $T^*_\tmb$ from $T^*_\tmb \x T^*_\tnb \x T^*_\tnb \x T^*_\tpb$, and where $\Delta(T^*_\tnb)$ is the fibration with fibre the diagonal in $T^*_{N_b}\x T^*_{N_b}$, and $\cap_B$ is the fibrewise intersection, then Hormander´s condition for classical FIOs extends naturally to the requirement for vertical FIOs that if each fibre of $(\Cs^1\x\Cs^2)$ intersects each fibre of $T^*_\tmb \x \Delta(T^*_\tnb) \x  T^*_\tpb$ transversally then $\Cs^1\circ \Cs^2$ is a smooth fibration and \eqref{comp calculus} holds. More generally, the less restrictive {\it clean} intersection calculus likewise holds for vertical FIOs. For Theorems 1 and 2 conveniently one has  $\Cs \< T^*_\tmb \x_B T^*_\tnb\bs 0$ is the graph $\Cs= {\rm gr}_\tmb(G)$ of a vector bundle map $G\in \Hom(T^*_\tnb, T^*_\tmb)$. A fundamental special case of this situation is for the algebra 
$$\Psi^*(M/B,E) := \If^*((M\times_B M)/B, {\rm gr}_\tmb(I))$$
of vertical pseudodifferential operators (v-$\pdo$s) defined for the identity bundle map $I: T^*_\tmb\to T^*_\tmb$. In general, without further constraints, $\If^*((M\times_B M)/B, \Cs)$ does not form an algebra, but it is always a left and right module for composition with $\Psi^*(M/B,E)$. This is immediate when $\Cs= {\rm gr}_\tmb(G)$. More generally, for v-FIOs $\Fs_i\in \If^{q_i} ((M \times_B M)/B, \Gs(\Ts_i))_{[d_i]}$ for $i=1,2$, associated to vertical canonical relations $\Cs_i = \Gs(\Ts_i)$ defined by graphs of vertical (fibrewise) canonical transformation maps $\Ts_i : T_\tmb^* \to T_\tmb^*$, that 
\begin{equation}\label{comp canon rels}
\Fs_2\circ \Fs_1\in \If^{q_1 + q_2} ((M \times_B M)/B, \Gs(\Ts_2\circ \Ts_1))_{[d_1 + d_2]}.
\end{equation}
 
\section{Proof of \thmref{thm1} and \thmref{thm2}}

Let $\pi: M \to B$ be a compact Riemannian fibration. This requires that $M$ and $B$ are closed manifolds with a choice of metric $g_\tmb$ on the tangent bundle along fibres $T_\tmb$, along with a choice of Riemannian metric $g_B$ on $T_B$; the choice of $g_B$ is not important but facilitates constructions. The metric $g_B$ can be lifted to $\pi^*T_B$, which combined with    $g_\tmb$ yields a Riemannian metric 
\begin{equation}\label{split metric}
g_M := g_\tmb + \pi^*g_B  \ \  {\rm on} \ \ T_M  = T_{M/B} + \pi^*T_B 
\end{equation} 
by assuming $T_{M/B}$ and $\pi^*T_B$ to be orthogonal. Associated to the vertical metric $g_\tmb =\{g_b \in \Ci(M_b, S^2T^*_{M_b}) \ |\ b\in B\}$ one has a smooth family of Laplace-Beltrami operators $\Delta_\tmb = \{\Delta_b  \ |\ b\in B\}$ with $\Delta_b$ the usual Laplace operator on the fibre $(M_b,g_b)$, given locally by $\Delta_{g_b} = \sqrt{\det g_b}\ii\partial_i \sqrt{\det g_b} \,g^{ij}_b\partial_j$. $\Delta_\tmb\in \Psi^2(M/B)_{[0]}$  defines an element of order 2 in the space of vertical $\pdo$s of form degree 0, and to this one may fibrewise construct the family of positive square root operators $\sqrt{\Delta_\tmb} = \{\sqrt{\Delta_{g_b}}  \ |\ b\in B\} \in \Psi^1(M/B)_{[0]}$ occurring - pointwise - in the distributional wave trace \eqref{wave tr}. The more general class of vertical $\pdo$ contemplated in Theorems (1.1) and (1.2) is of mixed differential-form degree and mixed order v-$\pdo$ 
$$ \Qs = \Qs_{[0]} + \Qs_{[1]} \cdots +  \Qs_{[\dim B]},$$
 where  
$$\Qs_{[j]} \in \Psi^*(M/B)_{[j]}$$
is a v-$\pdo$ of differential form degree $j$ on $B$ while the form degree 0 component is
\begin{equation}\label{P order 2}
\Qs_{[0]} = \Delta_\tmb.
\end{equation}
We assume $\Qs_{[j]}$  to have non-negative integer order 
$$\nu_j\in \N.$$
It is important to slightly broaden  this class  to  where the form degree zero part
\begin{equation}\label{P order m}
\Qs_{[0]} =  \Ps \in \Psi^m(M/B)_{[0]}
\end{equation}
is a smooth family $\{P_b  \ |\ b\in B\}$ of positive elliptic self-adjoint $\pdo$s of integer order $m>0$. $\Qs$ has constant order $(\nu_0, \nu_1, , \ldots, \nu_{\dim B})\in \mathbb{N}^{\dim B +1}$,  where $\nu_i ={\rm ord}(\Qs_{[i]})$ - so $\nu_0=2$ for \eqref{P order 2} and  $\nu_0=m$ for \eqref{P order m}. Whilst the non-negative integer order v-$\pdo$s  $\Qs_{[j]}$ for $j>0$ can be arbitrary, it is essential that  $\Qs_{[0]}$  be a smooth family of elliptic operators for the following constructions to be possible. \\

A geometrically natural example of such a v-$\pdo$ is the curvature $\Qs =\Bs^2$ of a Quillen-Bismut superconnection $\Bs$ \cite{Q,B,BGV} on $\pi_*(E^+ + E^-)$, adapted to a smooth family of formally self-adjoint elliptic $\pdo$s $\Ds =
\begin{bmatrix}\label{D pm}
  0 & \Ds^- \\
  \Ds^+ & 0 \\
\end{bmatrix}$. $\Bs$ is a
$\pdo$ on $\Aa(M/B, \pi_*(E^+ + E^-))  = \Ci(M,\pi^*(\wedge T^*_B)\otimes (E^+ + E^-)\otimes|\wedge_{\pi}|^{1/2})$ of odd-parity for the $\Zf_2$-grading with $\Bs( \omega\,\psi) = d\omega \, \psi + (-1)^{|\omega|}\omega\,\Bs(\psi)$  for $\omega\in\Aa(B)$ and $\psi\in \Aa(B,\pi_*(E^+ + E^-))$ with $\Bs_{[0]} = \Ds \ ,$ where $\Bs = \sum_{i=0}^{\dim
B}\Bs_{[i]}$ and $\Bs_{[i]}$  raises form
degree by $i$. It follows  that
$\Bs_{[1]}$ is a covariant derivative
while 
$\Bs_{[i]}\in\Aa^i(M/B,\Psi(E^+ + E^-))$  is a v-$\pdo$. 
In a local fibrewise
trivialization  $\Bs_{|W} = d_{B|W} + \sum_I P_{b,I} db_I,$ with $P_{b,I}$ a classical integer order $\pdo$ along the fibre. The curvature is the v-$\pdo$
$\Qs:= \Bs^2 \in \Aa(B,\pi_*(\Ee))$
with form degree zero part an even-parity vertical 2nd order Laplace-type elliptic operator 
$\Qs_{[0]} = \Ds^2 = \begin{bmatrix}
  \Ds^-\Ds^+ & 0 \\
  0 & \Ds^+\Ds^- \\
\end{bmatrix} \in \Aa^0(B,\pi_*(E^+ + E^-)).$
A spin structure on $T_\tmb$ and a bundle of Clifford
modules $E\to M$ equipped with a connection $\nabla^E$  restricting to a Clifford
connection on $E_{|M_z}$ defines an associated family of
compatible Dirac operators $\Ds_\tmb$. Let $\nabla^{\pi_*(\Ee)}$ be the
canonical Hermitian connection induced on $\pi_*(\Ee)$ defined by $\nabla^{\pi_*(E)}_X s = \nabla^E_{X^H}s + ws$ for $X^H$ the lifting of $X\in \Ci(B,T_B)$ to   $\Ci(M,\pi^*(T_B))$ specified by the splitting \eqref{split metric} and $w$ the divergence form of $g_\tmb$ with respect to $X^H$.  Finally, let $\cc(T)$ be the element of 
$\Aa^2(B,\End(\pi_*(E)))$  defined by Clifford multiplication by the
curvature tensor of the fibration associated to \eqref{split metric}. The Bismut superconnection $\Bf = \Bf_{[2]} + \Bf_{[1]} + \Bf_{[0]}$ on $\pi_*(E)$ is then defined by 
\begin{equation}\label{B-superconnection}
    \Bf = \cc(T) + \nabla^{\pi_*(E)}+ \Ds,
\end{equation} 
defining the v-$\pdo$ $$\Qs =\Bf^2  = \Bf^2_{[4]} + \Bf^2_{[2]}  + \Bf^2_{[1]} + \Bf^2_{[0]}$$
of order $(0,2,1,2)$; here, $\Bf^2_{[4]}$ a bundle endomorphism (of order 0), $\Bf^2_{[2]}$ a vertical differential operator of order 2, $\Bf^2_{[1]}$ is a vertical differential operator of order 1, and $$\Bf^2_{[0]}=\Delta_\tmb:=\Ds^2_\tmb$$ an elliptic vertical differential operator of order 2. We consider the vertical wave trace of $\Bf^2$ further below. \\

\subsection{The vertical $\pdo$ $\Qs^{1/m}\in \Psi^*(M/B)$:}

Let $\Qs=\sum_{k=0}^{\dim B}\Qs_{[k]}\in\Psi^*(M/B)$ be a
v-$\pdo$ of constant order
$$(\nu_0=m,\nu_1,\ldots,\nu_{\dim B}) \in\Nf^{\,\dim B +1},$$ 
with $\Ps = \Qs_{[0]}$ as above. 
Construction of the vertical wave operator requires knowledge of the microlocal structure of the v-$\pdo$ $\Qs^{1/m}$ on the fibration $\pi:M\to B$. $\Qs^{1/m}$ is defined by holomorphic functional calculus with form degree zero component
\begin{equation}\label{mth root}
    \Qs^{1/m}_{[0]} =\Ps^{1/m} \  \in \Psi^1(M/B)
\end{equation}
the smooth positive elliptic family $\{P^{1/m}_b  \ |\ b\in B\} $ in the classical single operator sense (for each $b$); for \eqref{P order 2}  this is
\begin{equation}\label{sq root}
    \Qs^{1/2}_{[0]} = \sqrt{\Delta_\tmb} :=  \{\sqrt{\Delta_{g_b}}  \ |\ b\in B\} \in \Psi^1(M/B).
\end{equation}
The principal symbol of $\Ps = \Qs_{[0]}\in \Psi^1(M/B)_{[0]}$, constructed fibrewise, has an invariant realization 
\begin{equation}\label{p0}
    {\bf p}_{0}\in \Ci(T^*_\tmb,\varphi^*(\End E))
\end{equation}
for the vertical tangent bundle $\varphi : T^*_\tmb \to M$.
 $\Qs$ is vertically-elliptic with
principal angle $\pi$ insofar as $\Qs_{[0]}:=\Ps
$ is parameter elliptic with principal angle
$\pi$ for each $b\in B$.    Thus,
 $p_{0} - \la I\in \Ci(T_\tmb\backslash\{0\},p^*(\End(E)))$ 
is an invertible bundle map for $\la\notin \R^1_+$ for the identity bundle operator $I$.  The  resolvent v-$\pdo$ $(\Qs-\la I)\ii$ for $\la\notin \R_+$ is constructed from the v-resolvent $(\Ps - \la I)\ii \in \Aa^0(B,\Psi^{\nu_0}(\Ee))$, ellipticity of $\Qs_{[0]}$ thus implying that of $\Qs$. For this, a vertical parametrix for $\Ps-\la I$ is constructed using a standard microlocal procedure: with respect to trivializations, inductively define vertical symbols $r_j[\la](b,x,y,\xi)$ quasihomogeneous
in $(\la,\xi)$ of degree $-\nu -j$  by
 $r_{-2}[\la](b,x,\xi) = ({\bf
p}_{0}(b,x,\xi) - \la I)\ii, $ and
$$r_{-2-j}[\la](b,x,\xi) = \ r_{-2}[\la](b,x,\xi) 
\sum_{\stackrel{|\a| + k + l = j}{k<j}}\frac{(-i)^{\a}}{\a
!}\partial^{\a}_{\xi}{\bf p}_k (b,x,\xi)
\partial^{\a}_{x}r _{-2-l} [\mu](b,x,\xi),$$
so that
$\left(\sum  p_k(b,x,\xi) - \la I\right)\circ\left(\sum r_{-2-j}[\la]\right) \sim {\bf I}.$
The vertical
poly-quasihomogeneous symbol $r[\la] \sim \sum r_{-2-j}[\mu]$ so defined yields a local parametrix which can be patched together to define $\Rs = {\rm
OP}(r)\in\Ci(M/B,\Psi^{-m}(E))$ with $(\Ps-\la\Is)\Rs = I - \Js_\infty$ where $\Js_\infty\in
\Ci(M/B,\Psi^{-\o}(\Ee))$ is a vertical smoothing operator. The $L^2$ operator norm of $\Js_\infty(b)$ is
$O(|\la|\ii)$ locally uniformly on $B$  and so $I-\Js_\infty$ is invertible in
$\Ci(M/B,\Psi(E))$ for sufficiently large $|\la|$ with
$(\Ps-\la\Is)\ii = R + R\sum_{j\geq 1}\Js_\infty^j.$ On the other hand, if $\As \in
\Aa^{>0}(B,\Psi(\Ee))$  so $\As$ (strictly) raises form degree in $\Aa(M/B,\pi_*(E))$,
then  $(I - \As)\ii = I + \As + \ldots + \As^{\dim B}$ in $\Aa(M/B,\Psi(E))$. 
We may write $$\Qs = 
\Ps + \Qs_+$$  with $$\Qs_+  = \Qs_{[1]} + \cdots \Qs_{[\dim B]} \in
\Aa^{>0}(M/B,\Psi(E))$$ 
and $$\Ps = \Qs_{[0]}\in
\Aa^0(B,\Psi(\Ee),$$ 
where  we assume each $\Qs_{[k]}$ has (constant) integer order. It follows for large $\la$ in an open sector $\G_\pi\subset \Cf \bs\R_+$ that  $(\Qs - \la \Is)\ii \in \Aa(B,\Psi(\Ee))$ 
with degree zero component equal pointwise to the usual
$\pdo$ resolvent for $\Ps := \Qs_{[0]}$ and 
\begin{eqnarray}
   (\Qs - \la \Is)\ii &=& (\Ps - \la \Is)\ii \left(\Is +
\Qs_+(\Ps - \la \Is)\ii\right)\ii\notag\\ 
& = & (\Ps - \la \Is)\ii  + \sum_{k=1}^{\dim B}
(-1)^{k}(\Ps - \la \Is)\ii \left(\Qs_+(\Ps - \la \Is)\ii\right)^k \label{res exp}\\
& = &  (\Ps - \la \Is)\ii  + \sum_{k=1}^{\dim B} (-1)^k(\Ps - \la \Is)\ii
\left(\sum_{i=0}^{\dim B} \Qs_{[i]}(\Ps - \la \Is)\ii\right)^k 
\end{eqnarray}
So as a sum over form degree on $B$
\begin{equation}\label{res exp d}
    (\Qs - \la \Is)\ii =  \sum_{d=0}^{\dim B} \l((\Qs - \la \Is)\ii\r)_{[d]}
\end{equation}
\begin{equation*}
=   (\Ps -
\la \Is)\ii  + \sum_{d=1}^{\dim B}\sum_{p_1+\ldots +p_k = d}  (-1)^k (\Ps -
\la \Is)\ii \Qs_{[p_1]}(\Ps - \la \Is)\ii \ldots \Qs_{[p_k]}(\Ps -
\la \Is)\ii,  
\end{equation*}
where  the inner sum is over integer partitions $p_1, \ldots, p_k$ of $d$ with $p_j>0$. 
The term 
\begin{equation*}
\Rs_{p_1, \ldots, p_k | d}[\la] := (\Ps -
\la \Is)\ii \Qs_{[p_1]}(\Ps - \la \Is)\ii \ldots \Qs_{[p_k]}(\Ps -
\la \Is)\ii
\end{equation*}
in \eqref{res exp d} is a v-$\pdo$ of order
\begin{equation}\label{e:sum}
  \nu_1 + \ldots + \nu_k - (k+1)m.
\end{equation}
Modulo vertical smoothing operators, $R[\la] = {\rm OP}(r_\D[\la])$ is the resolvent  parametrix for
$\Ps-\la\Is$ so constructed, $\Qs_{p_j} = {\rm OP}({\bf q}_{[p_j]})$, for weakly polyhomogeneous symbols
\begin{equation}\label{b w symbols}
r_\D[\la](b,x,\xi) \sim \sum_{l\geq 0}r_{-m-l}[\la](b,x,\xi), \ \ 
{\bf q}_{[p_j]}(b,x,\xi) \sim \omega_{[p_j]}\otimes \sum_{
l\geq 0}{\bf q}_{\nu_j-l}(b,x,\xi).
\end{equation}
Then $(\Qs - \la \Is)\ii$ has v-symbol in fibrewise coordinates  
$$ {\bf r}[\la](b,x,\xi)  = r_\D[\la](b,x,\xi) +  \sum_{d=1}^{\dim B}\sum_{p_1+\ldots +p_k = d} {\bf r}_{p_1,\ldots, p_k|d}[\la](b,x,\xi),$$
\begin{equation}\label{d parametrix symbols}
 {\bf r}_{p_1,\ldots, p_k|d}[\la](b,x,\xi) = r_\D[\la](b,x,\xi)\circ {\bf q}_{[p_1]}(b,x,\xi)\circ\ldots \circ{\bf
q}_{[p_k]}(b,x,\xi)\circ r_\D[\la](b,x,\xi),
\end{equation}
\vskip 2mm
Here, for v-symbols
 ${\bf a} = \omega_{[i]}\otimes a, {\bf b} =
\sigma_{[j]}\otimes b$ then
${\bf a}\circ {\bf b} = \omega_{[i]}\wedge \sigma_{[j]}\otimes
(a\circ b)$ 
for the graded tensor product $\otimes$,
and $a\circ b \sim \sum_j (a\circ b)_j$ with
$(a\circ b)_j = \sum_{|\a| + k + l=j}\frac{(-i)^{\a}}{a
!}\partial^{\a}_{\xi}(a)_k
\partial^{\a}_{x}(b)_l,$. In \eqref{b w symbols} 
${\bf q}_{\nu_j-l}(b,x,\xi)$ is homogeneous in $\xi$ of degree
$\nu_j-l$ and $r_{-m-l}[\la](b,x,\xi)$ is quasi-homogeneous in
$(\xi,\la^{1/m})$ of degree $-m-l$, i.e. $r_{-m-l}[t^m\la](b,x,t\xi)= t^{-m-l}r_{-m-l}[\la](b,x,\xi)$ 
$t>1$,  so modulo vertical smoothing operators $\Rs_{p_1, \ldots, p_k | d}[\la] ={\rm OP}({\bf r}_{p_1,\ldots, p_k|d}[\la])$ has poly-quasi-homogeneous symbol of degree $\nu_1 + \ldots + \nu_k - (k+1)m$
\begin{equation}\label{r}
    {\bf r}_{p_1,\ldots, p_k | d}[\la](b,x,\xi) \sim \sum_{l\geq 0}{\bf r}_{p_1,\ldots, p_k | d, l}[\la](b,x,\xi)_{[d]} 
\end{equation}
for fixed $(p_1,\ldots, p_k)$ of fixed differential form order $d=p_1 + \ldots + p_k$; so for $t>1$
\begin{equation}\label{res d symbol homog}
{\bf r}_{p_1,\ldots, p_k | d, l}[t^m\la](b,x,t\xi)_{[d]} \ \ = \hskip 5cm
\end{equation}
$$   \hskip 4cm  t^{\nu_1 + \ldots + \nu_k - (k+1)m-l}\,{\bf r}_{p_1,\ldots, p_k | d, l}[\la](b,x,\xi)_{[d]}. $$
Assume, as holds for the Bismut superconnection curvature, that
\begin{equation}\label{e:orderlessthanr}
 v_k = \ord(\Qs_{[k]}) \leq m, \hskip 10mm k=1,\ldots,\dim B.
\end{equation}
Then \eqref{res exp}, \eqref{res exp d}  imply an operator norm estimate in
$\Aa(B)$ as $\la\to\infty$ in $\G_\pi$
\begin{equation}\label{e:Olambdainverse}
\| (\Qs - \la \Is)\ii\|^{(l)}_{M/Z} = O(|\la|\ii) \ ,
\end{equation}
where 
$\| \ . \ \|^{(l)}_{M/B} : \Aa(M,\Psi^{0}(\Ee)) \too \Aa(B)$ 
is the vertical operator Sobolev $l$ norm associated to the vertical
metric  $ |
\pi^*(\a)\otimes\Psi\otimes\upsilon |_{M/B} = \a
\int_{M/B}|\psi|^2\,\upsilon^2$
for elements 
$\psi\otimes\pi^*(\a)\otimes\upsilon\in \G(M,\pi^*(\wedge T^*
_B)\otimes\Ee\otimes |\wedge|^{1/2}_{\pi})$ independently of the form of the tensor
product.  In \eqref{e:Olambdainverse},  $a(\la) = O(f(\la))$ for $a:\Cf\to
\Aa(B)$  if for each $l$ 
 there is a constant $C(l)$
such that $\|a(\la)\|_l \leq C(l)f(\la)$ for the  $C^l$ norm $\|\ . \ \|_l$ on $B$. The proof of \eqref{e:Olambdainverse} follows
 \cite{Se} using the vertical parametrix $\Bs$. Given \eqref{e:orderlessthanr}, and setting $$\dl\la := (i/2\pi)d\la,$$  the complex power $\Qs_\pi\si$ is defined for $\re(s)>0$ by
\begin{equation}\label{e:powers 3}
     \Qs_\pi\si = \int_{C}\la_\pi\si\, (\Qs -
\la\Is)\ii \, \dl\la,
\end{equation}
which by \eqref{res exp d} has the structure over form degree on $B$
\begin{eqnarray*}
&& \int_{C}\la_\pi\si\, (\Ps -\la \Is)\ii \, \dl\la\\ 
&& + \sum_{d=1}^{\dim B}\sum_{p_1+\ldots +p_k = d}  (-1)^k \int_{C}\la_\pi\si\,  (\Ps -
\la \Is)\ii \Qs_{[p_1]}(\Ps - \la \Is)\ii \ldots \Qs_{[p_k]}(\Ps -
\la \Is)\ii\, \dl\la, \notag 
\end{eqnarray*}
that is,
\begin{equation}\label{powerstructure}
\Qs_\pi\si =  \Ps_\pi\si   + \sum_{d=1}^{\dim B}\sum_{p_1+\ldots +p_k = d}  (-1)^k\int_{C}\la_\pi\si\,  (\Ps -
\la \Is)\ii \Qs_{[p_1]}(\Ps - \la \Is)\ii \ldots \Qs_{[p_k]}(\Ps -
\la \Is)\ii\, \dl\la,
\end{equation}
$$\hskip 20mm = \left.\Ps_\pi\si   \hskip 65mm \right\}= (\Qs_\pi\si)_{[0]}\hskip 45mm$$ 
$$ \hskip 27mm  + \left.\int_{C}\la_\pi\si\,  (\Ps -
\la \Is)\ii \Qs_{[1]}(\Ps - \la \Is)\ii \, \dl\la\hskip 20mm \right\} = (\Qs_\pi\si)_{[1]}\hskip20mm$$
$$\hskip 35mm\left.\begin{array}{l}
    \ + \int_{C}\la_\pi\si\,  (\Ps -
\la \Is)\ii \Qs_{[2]}(\Ps - \la \Is)\ii \, \dl\la \\[3mm]
\hskip 2mm + \int_{C}\la_\pi\si\,  (\Ps -
\la \Is)\ii \Qs_{[1]}(\Ps - \la \Is)\ii \Qs_{[1]}(\Ps - \la \Is)\ii  \, \dl\la
\end{array}\right\}  = (\Qs_\pi\si)_{[2]}$$

$$+ \cdots,$$

where $\la_\pi\si$ is the branch of $\la\si$ defined by cutting the complex plane along the negative real axis $\R_-$, and where $C$ is a negatively oriented contour encircling $\R_-$ and the origin. One notes that the form degree 0 component of $\Qs_\pi\si$ coincides with $ \Ps_\pi\si = \{P_b\si \ | \ b\in B \}$ the smooth family  of pointwise complex powers in the classical sense \cite{Se}. 
So defined, $\Qs_\pi\si$
is a v-$\pdo$ with from  \eqref{res exp},\eqref{res exp d} vertical symbol 
\begin{eqnarray}
    {\bf q}_\pi\si(b,x,\xi) &=& \int_{C}\la_{\pi}\si\, {\bf r}[\la](b,x,\xi) \ \dl\la \label{cx symbol}\\
 &=& \sum_{d=0}^{\dim B}  \sum_{p_1+\ldots +p_k = d} 
\int_{C}\la_{\pi}\si\, (-1)^k {\bf r}_{p_1,\ldots, p_k | d}[\la](b,x,\xi) \ \dl\la \notag 
\end{eqnarray}
\begin{equation}\label{e:powersexpansion} 
    \sim \sum_{d=0}^{\dim B}  \sum_{p_1+\ldots +p_k = d} \sum_{l\geq 0}(-1)^k
\underbrace{\int_{C}\la_{\pi}\si\,  {\bf r}_{p_1,\ldots, p_k | d, l}[\la](b,x,\xi)_{[d]} \ \dl\la}_{=\ {\bf q}\si_{p_1,\ldots, p_k |d,l}(b,x,\xi)},
\end{equation}  
\vskip 3mm
\noindent a sum of  classical polyhomogeneous v-symbols, with ${\bf q}\si_{p_1,\ldots, p_k |d,l}(b,x,\xi)$  from \eqref{res d symbol homog} homogeneous of degree 
 \begin{equation}\label{s homog}
      (- ms -mk +\nu_1 + \ldots + \nu_k)-l,
 \end{equation}
 that is, 
$${\bf q}\si_{p_1,\ldots, p_k |d,l}(b,x,t\xi) = t^{(- ms -mk+\nu_1 + \ldots + \nu_k) - l}\,{\bf q}\si_{p_1,\ldots, p_k |d,l}(b,x,\xi),$$ 
where $\nu_j$ is the order of the classical v-$\pdo$ $\Qs_{[p_j]}$.\\

The complex powers form a semi-group $\Qs_\pi\si\Q_\pi^{-w} = \Q_\pi^{-s-w}$  for $\re(s), \re(w) >0$ by the standard argument given in \cite{Se}, and are hence  extended to any $s\in\C$ by   choosing a
positive integer $N$ with $\re(s) + N > 0$ and defining (independently of the choice of  $N$  by the semigroup property) $\Qs_\pi\si = \Qs^N \Qs_\pi^{-s-N} \in \Psi^*(M/B,E).$
Specifically, 
\begin{equation}\label{1/m}
    \Qs_\pi^{1/m} := \Qs \Qs_\pi^{-1 +\frac{1}{m}} \in \Psi^*(M/B,E).
\end{equation}
One may retain the same formula \eqref{e:powersexpansion} for the symbol of $ \Qs_\pi^{1/m}$ provided $b$ is restricted to a sufficiently small bounded open set $O\< B$, for local fibrewise coordinates, and the contour $C$ is replaced by a simple closed contour $C_0$ which encloses the (finite number of) eigenvalues of the self-adjoint principal symbol ${\bf p}_0$ of $\Ps$  for each $b\in O$, as guaranteed by vertical ellipticity. Then,  $\Qs_\pi^{1/m} $ has local vertical symbol 
 \begin{equation}\label{1/m powersexpansion}
 {\bf q}_\pi^{1/m}(b,x,\xi) =\sum_{d=0}^{\dim B}  \sum_{p_1+\ldots +p_k = d}  (-1)^k 
\underbrace{\int_{C_0}\la_{\pi}^{1/m}\,{\bf r}_{p_1,\ldots, p_k | d}[\la](b,x,\xi) \ \dl\la}_{=\ {\bf q}^{1/m}_{p_1,\ldots, p_k |d,l}(b,x,\xi)},
\end{equation}
with ${\bf q}^{1/m}_{p_1,\ldots, p_k |d,l}(b,x,\xi)$   homogeneous of degree 
 \begin{equation}\label{1/m homog}
      (1 -mk +\nu_1 + \ldots + \nu_k) - l,
 \end{equation}
 so for $t>1$
$${\bf q}^{1/m}_{p_1,\ldots, p_k |d,l}(b,x,t\xi) = t^{(1 -mk +\nu_1 + \ldots + \nu_k) -l}\,{\bf q}^{1/m}_{p_1,\ldots, p_k |d,l}(b,x,\xi).$$ 
The composition on right-hand side of \eqref{1/m} is well defined since the canonical relation defining $\Qs$ is just the identity bundle map. Note also that
\begin{equation}\label{1/m d=0}
    (\Qs_\pi^{1/m})_{[0]} = \Ps_\pi^{1/m} \ \  \in \Psi^m(M/B,E)_{[0]}
\end{equation}
 is the smooth family  of pointwise $1/m$  powers $\{P_b^{1/m} \ | \ b\in B \}$   in the standard classical $\pdo$ sense. \\

The  assumption \eqref{e:orderlessthanr} is removed  using 
$ \| \la^{n-s}\, \dd_{\la}^{n-1}(\Qs - \la \Is)\ii\|^{(l)}_{M/Z} =
O(|\la|\ii)$ for integer $n>0$ as $\la\to\infty$ along $C$. For then, integrating by parts in
\eqref{powerstructure} gives
\begin{equation}\label{e:F-by-parts}
  \Fs\si  =  \frac{1}{(s-1)\ldots (s-n)}.\frac{i}{2\pi}\int_C\la^{n-s}\,
  \dd_{\la}^{n}(\Fs - \la\Is)\ii \, d\la,
\end{equation}
since
\begin{equation*}
\dd_{\la}^{n}(\Qs - \la\Is)\ii = \sum_{\stackrel{k=0}{n_0 + \ldots
+ n_k = n}}^{\dim B} (-1)^k\dd_{\la}^{n_0}(\Ps - \la \Is)\ii
\Qs\,\dd_{\la}^{n_2}(\Ps - \la \Is)\ii \ldots \Qs \,
\dd_{\la}^{n_k}(\Ps - \la \Is)\ii\,d\la
\end{equation*}
and $\dd_{\la}^{n_i}(\Ps - \la \Is)\ii \in
\Aa(B,\Psi^{-mn_i-m}(\Ee))$, then by taking $n\geq N$ for
sufficiently large $N$ we may ensure an estimate $\|
\dd_{\la}^{n}(\Qs - \la \Is)\ii\|^{(l)}_{M/Z} = O(|\la|\ii) $
without assuming \eqref{e:orderlessthanr}. For the general case,
one has for $\re(s) > n \geq N$ via
\eqref{e:F-by-parts}
$$
\Fs\si = \Ps\si + \sum_{k=0}^{\dim B}\frac{1}{(s-1)\ldots
(s-n)}.\frac{i}{2\pi}\int_{C}\la^{n-s}\, \dd_{\la}^n\left(\,(\Ps
- \la\Is)\ii (\Qs (\Ps - \la\Is)\ii )^k\,\right) \, d\la \ , $$
and then matters proceed as above, with the same symbol formulae. 
Thus, one has:
\begin{prop}\label{Qmthroot}
    \begin{equation*}
        \Qs_\pi^{1/m} := \Qs \Qs_\pi^{-1 +1/m} \in \Psi^*(M/B,E)
    \end{equation*}
is a vertical $\pdo$ associated to the fibration $\pi: M\to B$ with 
 \begin{equation}\label{1/m 0}
    (\Qs_\pi^{1/m})_{[0]} = \Ps^{1/m} \in \Psi^1(M/B,E)
\end{equation}
a positive elliptic family of classical $\pdo$s  and  
\begin{equation}\label{1/m d}
    (\Qs_\pi^{1/m})_{[d]} =  \sum_{p_1+\ldots +p_k = d} \Qs^{1/m}_{p_1,\ldots, p_k |d},
\end{equation} 
where
\begin{equation}\label{1/m d 2}
   \Qs^{1/m}_{p_1,\ldots, p_k |d} \in \Psi^{1 -mk +\nu_1 + \ldots + \nu_k}(M/B,E)_{[d]}
\end{equation}
with local vertical polyhomogeneous symbol of differential form degree $d$
\begin{equation}\label{symbol 1/m d 3}
{\bf q}^{1/m}_{p_1,\ldots, p_k |d}(b,x,\xi) \sim \sum_{l\geq 0} {\bf q}^{1/m}_{p_1,\ldots, p_k |d,l}(b,x,\xi),
\end{equation}
with ${\bf q}^{1/m}_{p_1,\ldots, p_k |d,l}(b,x,\xi)$   homogeneous in $\xi$ of degree 
 $(1 -mk +\nu_1 + \ldots + \nu_k) -l$.
\end{prop}

\subsection{The vertical wave operator $e^{it\Qs^{1/m}}$:}

$e^{it\Qs^{1/m}}$ is a vertical FIO whose vertical kernel is characterised as the kernel of the operator assigning to $u_0\in \Aa^*(M/B, \pi^*(E))$ the solution 
$$u\in \Aa^*((\uR\x_BM)/B, \pi^*(E))  := \Ci_0(\uR\x_BM, \pi^*(\wedge T^*_B) \otimes E\otimes |\wedge_\pi|^{1/2}), $$ 
$$u(t,m) = e^{it\Qs^{1/m}}u_0(m),$$  to the Cauchy problem 
\begin{equation}\label{Cauchy}
\frac{\partial u}{\partial t} + i\Qs_\pi^{1/m} u =0, \hskip 1cm u(0,m) = u_0(m),
\end{equation}
where $\uR \to B$ is the trivial fibration (line bundle) with fibre $\R^1$ at $b\in B$. In computing the fibrewise wave trace it is convenient to treat the  parameter $t$ as a parameter, so that for each $t\in\R^1$  the wave kernel $ \Us^t$ of $e^{it\Qs^{1/m}}$ is
 defined by a vertical Lagrangian distribution 
\begin{equation*}
    \Us^t \in   \If^*((M\times_B M)/B, E^*\boxtimes E, \Cs_t^\Us)
\end{equation*}
(as is done for the heat trace).  $\Cs_t^\Us$ is a vertical canonical relation equal to the graph  of the fibrewise flow $\Phi^t$ of the vertical Hamiltonian vector field of the form degree zero component 
\begin{equation}\label{wave kernel 0}
    \Us_{[0]}^t  = e^{it\Ps^{1/m}}  \in  \If^0((M\times_B M)/B, E^*\boxtimes E, \Cs_t^\Us)_{[0]}.
\end{equation}
Initially, however,  $t$ will be  treated as a variable.  Then  
\begin{equation}\label{wave kernel t var}
    \Us \in   \If^*(((\uR\x_BM)\times_B M)/B, E^*\boxtimes E, \tilde{\Cs}^\Us)
\end{equation}  
modifying by $-1/4$ the orders of the vertical FIOs, specifically 
\begin{equation}\label{wave kernel 0 t var}
    \Us_{[0]}    \in  \If^{-1/4}(((\uR\x_b M)\times_B M))/B, E^*\boxtimes E, \tilde\Cs^\Us)_{[0]},
\end{equation}
where $\tilde\Cs^\Us$ has from \cite{DG}  fibre at $b\in B$
\begin{equation}\label{CUb0t}
    \tilde C^{U_b} =\{(t,\tau_b) \in T_b^*\R\bs 0, (x,\xi), (y,\xi) \in T^*_{M_{b_0}}\bs 0 \ |  \ \tau_b + p_b^{1/m}(x,\xi)=0,  (x,\xi) = \Phi_b^t(y,\xi) \}.
\end{equation}
To compute the trace of the kernel \eqref{wave kernel t var} as integration over fibre of the diagonal in $M\times_B M$ a one expects a map (for vertical kernels smooth along the diagonal)
\begin{equation}\label{vertical wave trace}
    \Tr_\tmb: \If^*(((\uR\x_B M)\times_B M)/B, E^*\boxtimes E, \Cs^\Us) \to \Aa(B, \Sss^\prime(\R^1)):=\Aa(B)\ox_{\Ci(B)}\Sss^\prime(\R^1),
\end{equation}
to tempered distributions on $\R^1$ with values in forms on $B$; the objective is to compute the singularity structure of the trace (modulo smooth functions on $\R^1$) of the differential form coefficients. (In the following we may neglect reference to the bundle $E$.) The first step in constructing the solution to \eqref{Cauchy} is to do so for the form degree zero component $\Qs_{[0]}^{1/m}=\Ps^{1/m}$ of $\Qs$. Doing so requires nothing more than classical wave trace constructions and essentially nothing more than this is needed to solve the full vertical Cauchy problem \eqref{Cauchy},  just as $\Qs$ and the form degree zero part $\Qs_{[0]}^{1/m}= \Ps^{1/m} $ are enough to determine $\Qs^{1/m}$. 
\begin{prop}\label{0 form wave eqn}
Let $u_0\in \Aa^k(M/B, \pi^*(E))$. Then the Cauchy problem 
\begin{equation}\label{Cauchy 0}
\frac{\partial u}{\partial t} + i\Ps_\pi^{1/m} u =0, \hskip 1cm u(0,m) = u_0(m),
\end{equation}
has unique solution $u(t,m) = e^{-it\Ps^{1/m}}u_0(m)$, $u \in \Aa^k((\uR\x_B M)/B, \pi^*(E))$,  defined by a vertical FIO  of order 0 (and form degree 0) denoted $e^{-it\Ps^{1/m}}$ with Schwartz kernel a vertical Lagrangian distribution 
\begin{equation}\label{wave kernel}
    \Us_0 \in   \If^{-1/4}(((\uR\x_B M)\times_B M)/B, E^*\boxtimes E, \Cs^\Us)_{[0]}
\end{equation}
in which $\Cs^\Us$ is a vertical canonical relation equal to the graph  of the fibrewise flow $\Phi^t$ of the vertical Hamiltonian vector field associated to the vertical leading symbol of $\Ps^{1/m}.$  
\end{prop}
{\it Proof:} It is enough to construct $\Us_0$ locally around  each $b_0\in B$ and then patch together local kernels to obtain $ \Us_0$ via a partition of unity. Since $\pi:M\to B$ is locally trivial there is a neighbourhood $V_{b_0}\<B$ with 
$M_{|V_{b_0}} = \pi\ii(V_{b_0})$   identified fibrewise diffeomorphically with $M_{b_0} \x V_{b_0}$ and a vector bundle isomorphism $E_{|\pi\ii(V_{b_0})} \cong E_{|M_{b_0}}\x V_{b_0}$.  This defines a trivialisation $\pi_*E_{b_0}\cong V_{b_0}\times \Ci_0(M_{b_0}, E_{b_0}\otimes |\wedge_{M_{b_0}}|^{1/2})$  of $\pi_*E$ over $V_{b_0}$ in which $u \in \Aa^k((\uR\x_B M)/B, \pi^*(E))$ becomes
\begin{equation}\label{u deg k}
u_{|\pi\ii(V_{b_0})} =  \sum_{|I|=k} u_{b,I} \ db_I,
\end{equation}
where $u_{b,I}\in \Ci(\uR\x M_{b_0}, E_{b_0}\otimes |\wedge_{M_{b_0}}|^{1/2})$ and $db_I = db_{i_1}db_{i_2}\cdots db_{i_k}$ for coordinates $b_j$ on  $V_{b_0}$. In this localization, \eqref{Cauchy 0} is the system of PDEs 
\begin{equation}\label{Cauchy 0 b}
\frac{\partial u_{b,I}}{\partial t} +  P_b^{1/m} u_{b,I} =0, \hskip 1cm u_{b,I}(0,x) = u_{b,I,0}(x),
\end{equation}
on the fibre $M_{b_0}$, where $P_b^{1/m} = \Ps^{1/m}_{|\pi\ii(V_{b_0})}$ in the trivialization over $V_{b_0}$ acts on $\Ci(\R^1\x M_{b_0}, E_{b_0}\otimes |\wedge_{M_{b_0}}|^{1/2})$ and varies smoothly with $b\in V_{b_0}$ (or, rather, its classical Schwartz kernel varies smoothly by assumption that it is a vertical $\pdo$ - for $\Delta_b^{1/2}$ this follows from the choice of smooth vertical metric $g_\tmb$). \eqref{Cauchy 0 b} is  the classical Cauchy problem over the closed manifold $M_b$
with unique solution wave operator $e^{-it P_b^{1/m}}$ with kernel Lagrangian distribution
\begin{equation}\label{Ubt}
    U_{0,b} \in   \If^{-1/4}(\R^1\x M_{b_0})\times M_{b_0}, \End(E_{|M_{b_0}}), C^{U_b}),
\end{equation} as constructed and analysed by Chazarain \cite{C} and Duistermaat - Guillemin \cite{DG}, with respect to the canonical relation \eqref{CUb0t} with $\Phi_b^t$ the flow of  the canonical Hamiltonian vector field $H_b$ on $T^*_{M_{b_0}}$ defined by the leading symbol $p_b^{1/m}$ of $P_b^{1/m}$ in \eqref{p0}. Smooth dependence on $b$  of $\Phi_b^t$ is  ensured by the smooth dependence (by assumption) of $p_b^{1/m}$, and since  $p_b^{1/m} = ({\bf p}^{1/m}_{0})_{|T^*_{M_b}}$ is a section of the vector bundle $\varphi^*(\End E)$ over $T^*_\tmb$ then the Lagrangian submanifolds $C^{U_b}$ pull-back to a smooth trivial vertical Lagrangian subfibration of $(T^*_\tmb \x_B T^*_\tmb)_{|V_{b_0}}$, corresponding to the 
local triviality of a global vertical canonical 
relation $\Cs^\Us$ as a vertical Lagrangian subfibration of  $T^*_\tmb \x_B T^*_\tmb$. Similarly, the trivialization identifies the local family $U_{0,b}^t$ in \eqref{Ubt} with 
\begin{equation}\label{If loc}
 (\Us_0)_{|V_{b_0}} \in\If^0(((\uR\x_B M)\times_B M_{|V_{b_0}})/B E^*\boxtimes E, \Cs^\Us_{|V_{b_0}})_{[0]} 
\end{equation}
Patching together the so constructed solutions over open subsets $V_{b_0}$ of $B$ gives the required global exact solution  in $\If^0(((\uR\x_BM)\times_B M)/B, E^*\boxtimes E, \Cs^\Us)_{[0]}$. \hfill $\Box$  \vskip 3mm

The smooth variation of $U_{0,b}^t$ with $b$ in \eqref{If loc} is essentially the same fact that  the $C^{U_b}$ vary smoothly with b. Precisely, \eqref{If loc} means that $U_{0,b}$ is a finite sum of oscillatory integrals of the form 
\begin{equation}\label{local wave osc integral}
    I(b,t,x,y) = (2\pi)^{-q}\int_{\R^q} e^{i\phi(b,t,x,y,\xi)} a(b,t,x,y,\xi) \ d\xi
\end{equation}
with $\phi$ a non-degenerate phase function in $(x,y,\xi)$ parametrising $C^{U_b}$, insofar as
$d_\xi\phi(b,t,x,y,\xi)=0$ implies $((t, d_t\phi),(x, d_x\phi), (y, -d_y \phi)) \in C^{U_b}$, while $a(b,t,x,y,\xi)$ is an order $0$ classical symbol
$a(b,t,x,y,\xi) \sim \sum_{j\geq 0} a_{-j}(b,t,x,y,\xi)$ with $a_{-j}(b,t,x,y,\xi)$  homogeneous in $\xi$ of degree $-j$ and varying smoothly with $b$. \\

 From \cite{DG} one has that the vertical principal symbol of $\Us_0$ (that is, the principal symbol of $ U_{0,b}$ on $M_b$ for each $b\in B$) is non-vanishing and hence that 
 \begin{equation}\label{sing supp U0}
     {\rm sing \, supp}(\Us) 
 \end{equation}
$$  =\{(t,m,m')\in (\uR\x_B M)\times_B M \ | \ (m,\xi) =  \Phi^t(m',\eta) \}$$
some $\xi \in  T^*_{\tmb,m} = T^*_m (M_{\pi(m)}), \eta\in  T^*_{\tmb,m'} = T^*_{m'} (M_{\pi(m)})$; $(t,m,m')\in (\uR\x_B M)\times_B M $ means $\pi(m)=\pi(m') =b \in B$, corresponding to the classical fact fibrewise. For the vertical Laplacian $\Delta_\tmb$ this is  $$(t,m,m') = (t_b ,(b,x), (b,y))\in (\uR\x_B M)\times_B M$$  such that $x,y\in M_b$ can be joined by a geodesic of length $t_b$ in the fibre $M_b$, which we assume to form a subfibration of $\uR\x_B M\times_B M$, which along the diagonal becomes the set of closed geodesics $\gamma_{k,b}\< M_b$ in the fibres of $\pi$.  We may then take the fibrewise trace, integration over $x=y$ in $M_b$, to obtain the form degree zero families wave trace, immediate from the pointwise formulae of Chazarain \cite{C} and Duistermaat - Guillemin \cite{DG}, which for \eqref{P order 2} is
\begin{eqnarray*}
\Tr(e^{it\Qs^{1/2}})_{[0]} & = &  \Tr(e^{it\Delta_\tmb^{1/2}}) \\[2mm]    
& = &   \int_\tmb U(t,m,m) \\[4mm]
&\sim &  a_{\,0,-n+1}(t+i0)^{-n+1} + a_{\,0,-n+2}(t+i0)^{-n+2} + \cdots
\end{eqnarray*}
\begin{equation}\label{wave trace form zero part}
    +  \sum_{k\in\Z\bs 0} \Biggl( a_{\,L_\kappa,-1}(t-(L_\kappa+i0))^{-1}\biggr.  
 + \ \ \left.\sum_{j\geq 0}  a_{\,L_\kappa,j}(t-(L_\kappa+i0))^j \log(t-(L_\kappa+i0))\right), 
\end{equation} 
\vskip 2mm
\noindent where   $$L_\kappa \in \Ci(B)$$ assigns to $b\in B$ the length of the geodesic $\gamma_{k,b}$ and the $a_{r,l} \in \Ci(B)$  evaluated at $b$ are the wave trace coefficients for $e^{it P_b^{1/m}}$.\\

We turn now to the extension of this to the general case in \propref{Qmthroot} in which the strictly positive form degree component of $\Qs$, and hence that of $\Qs^{1/m}$, is non-zero. We shall treat the variable $t\in \R^1$ as a parameter in the following, rather than as a Schwartz kernel variable. With $t$ as a parameter \eqref{wave kernel} is replaced by 
\begin{equation}\label{wave kernel 00}
    e^{-i t\Ps^{1/m}} \in   \If^0(M\times_B M\, / B, E^*\boxtimes E, \Cs^\Us)_{[0]}.
\end{equation}
\begin{prop}\label{form wave operator}
The solution $e^{-it\Qs^{1/m}}$ of \eqref{Cauchy} exists as a vertical FIO with Schwartz kernel a vertical distribution $\Us$ on M $\times_B M$ whose $d$-form component $\Us_{[d]}$ is a sum of  vertical Lagrangian distribution kernels 
$$\Us^{(\Pp_d)}_{[d]}  \in \If^{\mu_{\Pp_d}} ((M \times_B M)/B, E^*\boxtimes E, \Gs(\Phi^t))_{[d]},$$
where $\mu_{\Pp_d}$ is a positive integer (specified later) depending on a set of values $\Pp_d$ derived from a partition  $d_1 + \cdots + d_k = d$ of $d$ and the orders of the v-$\pdo$s $\Qs^{1/m}_{p_1,\ldots,p_r |d_+r}$.
Each such summand can  be written modulo a smooth family of smoothing operators as a finite sum of vertical oscillatory integrals of the form
\begin{equation}\label{kernel Ff}
    \Ws(t,m,m') = (2\pi)^{-q}\int_{\Vv/(M\x_B M)} e^{i\phi(t,((m,m'),\xi))}\, {\bf w}(t,(m,m'),\xi)_{[d]}
\end{equation}
for $\Vv\to M\x_B M$ a fibration of vector spaces, $(m,m')\in M\x_B M$, so $\pi(m) = \pi(m')$, and $ ((m,m'),\xi) \in \Vv_{(m,m')}$ and $s\in\sigma_k$,  and 
for $\phi$ a non-degenerate vertical phase function parametrising the canonical relation $\Gs(\Phi^t)$ associated to the form degree 0 component $e^{-i t\Ps^{1/m}}$, and ${\bf w}$ a classical v-symbol of order $$\mu_d:=\mu_{\Pp_d}$$ 
with
$${\bf w}(t,(m,m'),\xi)_{[d]} \sim \sum_{j\geq 0} {\bf w}_{\mu_d-j}(t,(m,m'),\xi)_{[d]}$$
and  ${\bf w}_{\mu_d-j}$  homogeneous in $\xi$ of degree $\mu_d-j$ and varying smoothly with $b$ and of differential form degree $d$ on $B$.  In local fibrewise coordinates $m=(b,x), m'=(b,y)$ with $(x,y)\in M_b$ 
\begin{equation}\label{local kernel Ff}
\Ws_t(b,x,y) = (2\pi)^{-q}\int_{\R^r} e^{i\phi(b,t,x,y,\xi)} {\bf w}(b,t,x,y,\xi)_{[d]} \ d\xi
\end{equation}
with $\phi$ a non-degenerate  phase function in $(x,y,\xi)$ parametrising fibrewise  locally an open cone in $\Gs(\Phi_b^t)$, and with  ${\bf w}_{\mu_d-j}$  homogeneous in $\xi$ of degree $\mu_d-j$
$${\bf w}(b,t,x,y,\xi)_{[d]}  \sim \sum_{j\geq 0} {\bf w}_{\mu_d-j}(b,t,x,y,\xi)_{[d]}.$$
\end{prop} 
{\it Proof:}  In Duistermaat - Guillemin \cite{DG} the pointwise wave kernel $U_{0,b}$ of $e^{-itP_b}$ (of \propref{0 form wave eqn}) is constructed from  a parametrix using Duhamel's principle. A second application of this principle effects the extension needed here from the wave kernel for a family of classical operators (of form degree 0) $\Ps^{1/m} \in
\Aa^0(B,\Psi(\Ee)$ to that for
$$\Qs^{1/m} = \Ps^{1/m} + \Qs^{1/m}_+$$
with 
$$\Qs^{1/m}_+  = \Qs^{1/m}_{[1]} + \cdots \Qs^{1/m}_{[\dim B]} \in
\Aa^{>0}(B,\Psi(\Ee))$$ 
of form degree strictly positive. Concretely, the obstruction  to the parametrix $e^{-it\Ps^{1/m}}$ solving \eqref{Cauchy} is the vertical FIO
$$(\partial_t + \Qs^{1/m})e^{-it\Ps^{1/m}} \  \stackrel{\eqref{Cauchy 0}}{=} \ i\Qs^{1/m}_+ e^{-it\Ps^{1/m}}\ \in \ \If^*((\uR\x_BM)\times_B M)/B, E^*\boxtimes E, \Cs^\Us)$$
with the same vertical canonical relation as $e^{-it\Ps}$, see below.
Define the simplex $\sigma_k =\{(t_1, \ldots, T_\kappa) \ | \ 0\leq t_1 \leq \cdots \leq T_\kappa \leq 1 \}$, which can also be parametrised by $s_0 = t_1, s_j = t_{j+1} - t_j, s_k = 1-T_\kappa$ with $s_0 + \cdots + s_l = 1$ and $0\leq s_k\leq 1$. For $k>0$ let
$$  \As^{(k)}_t = \int_{t\sigma_k} i^k\, e^{-i(t-T_\kappa)\Ps^{1/m}}\Qs^{1/m}_+ e^{-i(T_\kappa-t_{k-1})\Ps^{1/m}} \cdots e^{-i(t_2 - t_1)\Ps^{1/m}}  \Qs^{1/m}_+ e^{-it_1\Ps^{1/m}}  \ dt_1\cdots dt_{k-1} $$
\begin{equation}\label{Akt}
= \int_{\sigma_k} (it)^k\, e^{-i s_0 t\Ps^{1/m}}\Qs^{1/m}_+ e^{-is_1 t\Ps^{1/m}} \cdots \Qs^{1/m}_+ e^{-i s_k t \Ps^{1/m}}  \ ds_1\cdots ds_k, 
\end{equation}
where $t\sigma_k :=\{(t_1, \ldots, T_\kappa) \ | \ 0\leq t_1 \leq \cdots \leq T_\kappa \leq t \}$. $\As^{(k)}_t$ is a vertical FIO in $\If^*(((\uR\x_BM)\times_B M)/B, E^*\boxtimes E, \Cs^\Us)$ which we will examine in a moment. Set
$$B^{(k)}(t) = \int_{s\sigma_{k-1}} i^{k-1}\, \Qs^{1/m}_+ e^{-i(t-t_{k-1})\Ps^{1/m}} \cdots \Qs^{1/m}_+ e^{-i(t_2 - t_1)\Ps^{1/m}} \Qs^{1/m}_+ e^{-it_1\Ps^{1/m}}  \ dt_1\cdots dt_{k-1} $$
for $k>1$ and $B^{(1)}(t) =  e^{-it\Ps^{1/m}}$. Then $\As^{(k)}_t = \int^t_0 e^{-i(t -y) \Ps^{1/m}} B^{(k)}(y) \ dy$, from which 
 $$(\partial_t + \Qs^{1/m})\As^{(k)}_t = B^{(k+1)}(t) + B^{(k)}(t)$$
and hence (with $\As^{(0)}_t := e^{-it\Ps^{1/m}}$)
\begin{equation}\label{sum solution}
    (\partial_t + \Qs^{1/m})\l(\sum_{k\geq 0} (-1)^k \As^{(k)}_t\r) = 0.
\end{equation}
The task is to see that the summands $\As^{(k)}_t$ are well-defined vertical FIOs, and to determine their  microlocal class; convergence is not an issue as the sum is  finite by nilpotence of the de Rham algebra of non-zero form degree. Given that these properties hold, and since from \eqref{Akt}
$$ \lim_{t\to 0} \As^{(k)}_t =0,  \ \ k>0,$$
then comparing \eqref{Cauchy}, \eqref{Cauchy 0} and \eqref{sum solution}  gives  
\begin{eqnarray}
     e^{-it\Qs^{1/m}} & = & \sum_{k=0}^{\dim B} (-it)^k  \int_{\sigma_k} e^{-i s_0 t\Ps^{1/m}}\Qs^{1/m}_+ e^{-is_1 t\Ps^{1/m}} \cdots \Qs^{1/m}_+ e^{-i s_k t \Ps^{1/m}}  \ ds  \notag \\
   & =&   e^{-i t\Ps^{1/m}}  -it   \int_0^1 e^{-i u t\Ps^{1/m}}\Qs^{1/m}_+ e^{-i(1-u)t\Ps^{1/m}} \ du \  + \ \cdots \label{Cauchy Qs}
\end{eqnarray}
where $ds:= ds_1\cdots ds_k$. This has $d$-form component   $1\leq d \leq \dim B$
$$\l(e^{-it\Qs^{1/m}}\r)_{[d]} = \sum_{\stackrel{d_1+\ldots +d_k = d}{1\leq d_j \leq d}}  (-it)^k  \int_{\sigma_k} e^{-i s_0 t\Ps^{1/m}}\Qs^{1/m}_{[d_1]} e^{-is_1 t\Ps^{1/m}} \cdots \Qs^{1/m}_{[d_k]} e^{-i s_k t \Ps^{1/m}}  \ ds$$
\begin{equation}\label{Cauchy Qs d} = \sum_{\stackrel{d_1+\ldots +d_k = d}{1\leq d_j \leq d}} \ \sum_{\stackrel{p_{j,1}+\ldots +p_{j,r_j} = d_j}{1\leq p_{j,i} \leq d_j}} (-it)^k  \int_{\sigma_k} e^{-i s_0 t\Ps^{1/m}}\Qs^{1/m}_{p_{1,1},\ldots, p_{1,r_1} |d_1} e^{-is_1 t\Ps^{1/m}} \cdots \Qs^{1/m}_{p_{k,1},\ldots, p_{k,r_k} |d_k} e^{-i s_k t \Ps^{1/m}}  \ ds,  
\end{equation}
where $1\leq k, r_j \leq \dim B$ and where from \propref{Qmthroot} $$\Qs^{1/m}_{p_{j,1},\ldots, p_{j,r_j} |d_j} \in \Psi^{1 -mr_j +\nu_{j,1} + \ldots + \nu_{j,r_j}}(M/B,E)_{[d_j]}$$ 
is a v-$\pdo$ with local vertical polyhomogeneous symbol 
$${\bf q}^{1/m}_{p_{j,1},\ldots, p_{j,r_j} |d_j}(b,x,\xi) \sim \sum_{l\geq 0} {\bf q}^{1/m}_{p_{j,1},\ldots, p_{j,r_j} |d_j,l}(b,x,\xi)$$ of  differential  form order 
$d_j$ with ${\bf q}^{1/m}_{p_{j,1},\ldots, p_{j,r_j} |d_j,l}(b,x,\xi)$  homogeneous in $\xi$ of degree 
 \begin{equation}\label{order Q fractions}
 (1 -mr_j +\nu_{j,1} + \ldots + \nu_{j,r_j}) -l,    
 \end{equation} 
 where $\nu_{j,q}$ is the order of the (by assumption) v-$\pdo$ $\Qs_{[p_{j,q}]}$.
(By an abuse of terminology, the operator is written here instead of its vertical Schwartz kernel in the space of vertical Lagrangian distributions indicated.) \\

\eqref{Cauchy Qs d} is a v-FIO defined by composition of the component v-FIOs $e^{-i s_l t\Ps^{1/m}}$ and $\Qs^{1/m}_{p_{i,1},\ldots, p_{i,r_i} |d_i}.$ The latter as a v-$\pdo$ is an element of
$$\Qs^{1/m}_{p_{i,1},\ldots, p_{i,r_i} |d_i} \in \If^{1 -mr_i +\nu_{i,1} + \ldots + \nu_{i,r_i}} ((M \times_B M)/B, E^*\boxtimes E, \Gs(I_\tmb))_{[d_i]},$$
where $\nu_{i,j}$ is the order of the classical v-$\pdo$ $\Qs_{[p_{i,j}]}$, and hence from \eqref{wave kernel 00}
\begin{equation}\label{wave o pdo}
    e^{-i s_0 t\Ps^{1/m}} \circ \Qs^{1/m}_{p_{1,1},\ldots, p_{1,r_1} |d_1}
    \end{equation}
    $$\in \If^{1 -mr_1 +\nu_{1,1} + \ldots + \nu_{1,r_1}} ((M \times_B M)/B, E^*\boxtimes E, \Gs(\Phi^{s_0 t}))_{[d_1]},$$
and, likewise, 
\begin{equation}\label{wave o pdo o wave o pdo}
e^{-i s_0 t\Ps^{1/m}} \circ \Qs^{1/m}_{p_{1,1},\ldots, p_{1,r_1} |d_1} \circ 
e^{-i s_1 t\Ps^{1/m}} \circ \Qs^{1/m}_{p_{2,1},\ldots, p_{2,r_2} |d_2} 
\end{equation}
$$\in \If^{2 -m (r_1 + r_2) +  \sum_{n=1}^{r_1} \nu_{1,n} + \sum_{n=1}^{r_2} \nu_{2,n}} ((M \times_B M)/B, E^*\boxtimes E, \Gs(\Phi^{(s_0 + s_1) t}))_{[d_1 + d_2]},$$
and iteratively
\begin{equation}\label{form deg d higher wave}
 \Ff_{t,\Pp,d}:= e^{-i s_0 t\Ps^{1/m}}\Qs^{1/m}_{p_{1,1},\ldots, p_{1,r_1} |d_1} e^{-is_1 t\Ps^{1/m}} \cdots \Qs^{1/m}_{p_{k,1},\ldots, p_{k,r_k} |d_k} e^{-i s_k t \Ps^{1/m}}
\end{equation}
$$\in \If^{k -\sum_{i=1}^k (mr_i  +  \sum_{n=1}^{r_i} \nu_{i,n})} ((M \times_B M)/B, E^*\boxtimes E, \Gs(\Phi^t))_{[d]},$$
since $\sum_{i=1}^k  s_i=1 $ and $\sum_{i=1}^k  d_i=d $, and where 
$$\mu_{\Pp_d} = k -\sum_{i=1}^k (mr_i  +  \sum_{n=1}^{r_i} \nu_{i,n})$$
with $\Pp$ denoting the set of data of an unordered non-zero positive integer partition $d_1 + \cdots + d_k = d$ of $d$ along with the sub-partition $p_{i,1} + \cdots + p_{i,r_i} = d_i$ of length $r_i$ (for $1\leq i \leq k$)  of each $d_i$. \\

The vertical canonical relation $\Gs(\Phi^t)$ in \eqref{form deg d higher wave} is independent of $s\in\sigma_k$. This depends, first, on \eqref{comp canon rels}. 
Since $\Qs^{1/m}_{p_{1,1},\ldots, p_{1,r_1} |d_1}$ has vertical canonical relation $\Gs(I\tmb)$ and  the vertical wave operator  $ e^{-i s_0 t\Ps^{1/m}}$ has $\Gs(\Phi^{s_0 t})$, with $\Phi^{s_0 t}$ the fibrewise Hamiltonian flow associated to the principal symbol of $\Ps^{1/m}$, then \eqref{wave o pdo}. Concretely, in local fibrewise coordinates $e^{-i s_0 t\Ps^{1/m}} \circ \Qs^{1/m}_{p_{1,1} }$ is a sum of oscillatory integrals
$ \int_{\R^q} e^{i\phi(b,s_0 t,x,y,\xi)} {\bf w}(b,s_0, t,x,y,\xi) \ d\xi$ with $\phi$ a phase function \eqref{local wave osc integral} parametrising the canonical relation $C_b = G(\Phi_b^{s_0 t})$, and with  ${\bf w}(b,s_0t,x,y,\xi)$ asymptotically as $\xi\to\infty$ equal to
$$  \sum_{\alpha, \beta}\frac{1}{\alpha!\beta!}\,\partial_\xi^{\alpha+\beta} D_x^\alpha {\bf q}^{1/m}(b, y,-\phi_y(x,y,\xi)) \cdot D^\beta_z (e^{i\tilde\phi(x,y,z,\xi)} a(b,s_0t, x,y,\xi))|_{z=y},$$ 
where $\tilde\phi(x,y,z,\xi) = \phi(b,s_0 t,x,z,\xi) - \phi(b,s_0 t,x,y,\xi) + (y-z)\phi_y(b,s_0 t,x,y,\xi)$. Likewise, in  \eqref{wave o pdo o wave o pdo},  for any $t_1, t_2$,  and v-$\pdo$ $A$, the v-FIO $e^{-i  t_1\Ps^{1/m}} A e^{-i  t_2\Ps^{1/m}} $  is in the same  class of vertical Lagrangians as $e^{-i  (t_1 + t_2)\Ps^{1/m}}$,  i.e. associated with $\Phi^{t_1+t_2}$, with symbol and order, and form degree, determined accordingly. \\

Finally, since the phase function is independent of $s\in\sigma_k$, then \eqref{form deg d higher wave} can be written as a sum of vertical oscillatory integrals over the compact set $\sigma_k$ which depend on the parameter $s$ only in the symbol, whilst integration over $\sigma_k$ results in  a symbol of the same classical class. \\

Hence, the conclusion is obtained. \hfill $\Box$  \vskip 3mm

\subsection{The vertical wave trace $\Tr_\tmb(e^{it\Qs^{1/m}})$:}

From \eqref{Cauchy Qs} 
\begin{equation}\label{Cauchy Qs d 2} 
\l(e^{-it\Qs^{-\frac{1}{m}}}\r)_{[d]}  = \sum_{k=0}^{\dim B} \sum_{\Pp_k} \ (-it)^k  \Ff_{t,\Pp_k,d},  
\end{equation}
where $\Pp_k$ is the data comprising an unordered integer partition $d_1, \ldots, d_k$ of $d$, so $d_1+\cdots +d_k = d$, along with sub-partitions $p_{i,1},\ldots, p_{i,r_i}$ of each  $d_j$, so $p_{i,1}+\cdots +p_{i,r_i} = d_i$, for $1\leq d_i, p_{i,i} \leq d_i$. With $\nu_{i,n}={\rm \ord}(\Qs_{[p_{i,n}]})\in\mathbb{N}$,
\begin{equation}\label{Ffd} 
\Ff_{t,\Pp,d} \in \If^{\mu_d} (M \times_B M/B, E^*\boxtimes E, \Gs(\Phi^t))_{[d]}.
\end{equation}
where 
\begin{equation}\label{mu}
\mu_d = k -\sum_{i=1}^k (mr_i  +  \sum_{n=1}^{r_i} \nu_{i,n}).    
\end{equation}
As a vertical Lagrangian distribution, $\Ff_{t,\Pp_k,d}$ may be written as a finite sum of oscillatory integrals of the form \eqref{kernel Ff}, \eqref{local kernel Ff}. \\

At this point, we suppose (as satisfied by the Bismut superconnection curvature)
\begin{equation}\label{max order m}
 v_k = \ord(\Qs_{[k]}) \leq m, \hskip 10mm k=1,\ldots,\dim B,
\end{equation}
so  by \eqref{s homog} the order of each summand of the v-$\pdo$ $(\Qs^{-\frac{1}{m}})_{[d]}$ is strictly negative
\begin{equation}\label{order Qd <0}
  \ord((\Qs^{-\frac{1}{m}})_{[d]}) < 0, \hskip 10mm d=0,\ldots,\dim B.
\end{equation}
 Then as in \cite{DG}, since $ \Us(t,m,m') =   i\Qs^{-\frac{1}{m}} \partial_t \Us(t,m,m'), $ repeated integration in $t$ by parts gives that the functional 
$$\rho \mto \int_{-\infty}^\infty \Us(t,m,m') \rho(t)\ dt$$
is a continuous linear mapping $\Sss(\R)\to \Aa(M \times_B M/B, E^*\boxtimes E)$ to vertical smoothing kernels, which hence may be combined with integration over the fibre 
$$\Tr_\tmb \Us : \rho \mto \Tr_\tmb \int_{-\infty}^\infty \Us(t,m,m') \rho(t)\ dt$$
to define a tempered distribution on $\R$ with values in the de Rham algebra $\Aa(B)$, resulting in the vertical wave trace
\begin{equation}\label{wave trace along diagonal}
    \Tr_\tmb(e^{-it\Qs^{\frac{1}{m}}}) =   \Tr_\tmb \Us = \int_\tmb  \tr(\Us(t,m,m)) 
\end{equation}
$$\ \in \Aa(B)\ox_{\Ci(B)}\Sss^\prime(\R).$$
Here, $\Us(t,m,m) $ is the restriction of $\Us(t, \ ,\ )$ to the fibrewise diagonal in $M\x_BM$, identified with a copy of $M\to B$, with $\Us$ restricting there to a continuous map $\Sss(\R)\to \Aa(M, \End(E))$. Equivalently, 
$$ \Tr_\tmb \Us  = \pi_* \delta^*\Us$$
for $\delta: M\hookrightarrow M \x_B M$ the inclusion map onto the fibrewise diagonal and $\pi_*$  the Gysin (Umkehr) push-forward identified on smooth differential forms with integration over the fibre. \\

That the restriction $\delta^*\Us$ to the diagonal in space variables is defined in this way may also be inferred from  \eqref{Ffd} insofar as   $\Ff_{t,\Pp_k,d}$ is defined with respect to the same canonical relation as the classical vertical wave operator $e^{-itP^{1/m}}$ then the microlocal singularity structure is the same. Concretely, $e^{-it\Qs^{-\frac{1}{m}}}$ may be written as oscillatory integrals with the same vertical phase function as for $e^{-itP^{1/m}}$  and hence with the same singular support. One has  
\begin{equation}\label{Cauchy Qs d trmb} 
\Tr_\tmb\l(e^{-it\Qs^{-\frac{1}{m}}}\r)_{[d]}  = \sum_{k=0}^{\dim B} \sum_{\Pp_k} \ (-it)^k\, \Tr_\tmb\l( \Ff_{t,\Pp_k,d}\r). 
\end{equation}
From \cite{DG} the principal symbol of the vertical Schwartz kernel $\Us_{[0]}(b)$ of $e^{-itP^{1/m}}(b)$ is nowhere zero and hence the same is so for the full vertical  principal symbol of $\Us_{[0]}$ and so along the diagonal, identifying with the canonical relation $\Cs = \Gs(\Phi^t) \< T^*_\tmb \x_B T^*_\tmb$,
 $${\rm sing \,supp}(\delta^*\Us) = {\rm sing \,supp}(\delta^*\Us_{[0]}) $$ $$=  \{(t,(m,\xi))\in \R\x_B T^*_\tmb \ | \ \exists \, \xi\in T^*_{M_{b=\pi(m)}}, \ (m,\xi) = \Phi^t(m,\xi)\}.$$
The fibrewise Hamiltonian flow $\Phi^t$ commutes with positive scalars on $\xi\in \Vv_m$ and  so \cite{DG} the singular support is unchanged on restriction to the vertical cosphere bundle 
$$=  \{(t,(m,\eta))\in \R\x_B S^*_\tmb \ | \ \exists \, \eta\in S^*M_{b=\pi(m)}, \ (m,\eta) = \Phi^t(m,\eta)\}.$$
 $\Tr_\tmb\l( \Ff_{t,\Pp_k,d}\r)$ is, modulo smoothing operators and modulo $\Ci(\R)$, a finite sum of vertical oscillatory integrals of the form 
\begin{equation}\label{kernel Ff diag}
    \int_\tmb\Ws_{\Pp_k}(t, m,m) = (2\pi)^{-q}\int_\tmb \int_{\Vv/M} e^{i\phi(t,((m,m),\xi))}\, {\bf w}_{\Pp_k}(t,(m,m),\xi)_{[d]}
\end{equation}
for $\Vv\to M$ a fibration of rank $q$ vector spaces and $ ((m,m),\xi) \in \Vv_m := \Vv_{(m,m)}$,   
$${\bf w}_{\Pp_k}(t,(m,m'),\xi)_{[d]} \sim \sum_{j\geq 0} \underbrace{{\bf w}_{\mu_d-j}(t,(m,m'),\xi)_{[d]}}_{{\rm homog\  in} \ \xi \ {\rm deg} \ \mu_d-j},$$
$\mu_d$ as in \eqref{mu}, and $\phi$ is a non-degenerate vertical phase function parametrising the vertical graph canonical relation $\Gs(\Phi^t)$ associated to the form degree 0 component $e^{-i t\Ps^{1/m}}$. Evaluated pointwise at $b\in B$ \eqref{kernel Ff diag} becomes
\begin{equation}\label{kernel Ff diag b}
    \l(\int_\tmb\Ws_{\Pp_k}(t, m,m)\r)(b) = (2\pi)^{-q}\int_{M_b} \int_{\R^q} e^{i\phi(b,t,x,x,\xi)}\, {\bf w}_{\Pp_k}(b,t,x,x,\xi)_{[d]} \ d\xi
\end{equation}
with ${\bf w}_{\Pp_k}(b,t,x,x,\xi)_{[d]}  \sim \sum_{j\geq 0} {\bf w}_{\mu_d-j}(b,t,x,x,\xi)_{[d]}.$\\

Let $\{T_\kappa(b) \ | \ \kappa\in\Z\}$ be the discrete set of periodic $H_{q_b}$ solution curves of period $T_b$ on the fibre $M_b$. For each $\kappa\in\Z$ let $L_\kappa :B\to \R^1$ be the function assigning smoothly to $b\in B$ the period $T_\kappa(b)$ on $M_b$. One may think of $L$ as a section of the trivial line bundle $\underline{\R}\to B$. We assume that the period $T_\kappa(b)$ varies smoothly with $b$, for $\Pp=\Delta$ this is just that the geodesic length varies smoothly with the Riemannian metric $g_{M_b}$ on each fibre, and the (graphs of the) sections $L_\kappa$ define in this way a trivial  $\Z$-fold covering space of $B$.  In particular, we may take $L_0$ to be the zero section. 
 Since $B$ is compact we may then consider $t$ in a collar neighbourhood of the graph of $L_\kappa$, so  $t\in (T_\kappa(b)-\epsilon_b, T_\kappa(b)+\epsilon_b)$  with $\epsilon_b$ varying smoothly with $b$ such that no other periods on $M_b$ other than $T_\kappa(b)$ occur in the interval. For such $t$, we shall then say that `$t$ is near to $L_\kappa$'.\\
 
According to the characterization of the singularity structure of $\Us$, Duistermatt-Guillemin \cite{DG} identify the critical manifolds of the vertical wave phase function fibrewise with the submanifold $Z_b \< S^*M_b $ consisting of periodic $H_{q_b}$ solution curves of period $T_\kappa(b)$, coinciding with the the fixed point set of $\Phi_b^{T_\kappa(b)}:S^*M_b \to S^*M_b $, which is assumed to be clean. $Z_b$ is assumed to be a union of even-dimensional connected submanifolds $Z_b = \sqcup_{p=1}^J Z_{b,p} \< S^*M_b$ with $Z_{b,p}$ a clean fixed point set of $\Phi_b^{T_\kappa(b)}$ of codimension $n_p$. We further assume that 
\begin{equation}\label{Zz}
    \Zz :=\cup_b Z_b \too B
\end{equation} 
forms a subfibration of the vertical cosphere bundle (along the fibres) $S^*_\tmb$ - viewed as a fibration with fibre the cosphere bundle over $M_b$, and that the submanifolds $Z_{b,p}$ likewise form disjoint connected subfibrations $\Zz_p \to B$ with fibre codimension $n_p$. \\

Guillemin-Melrose in \cite{GM}, in extending the Poisson summation formula to manifolds with boundary, observe that consequent on \cite{DG} for $t$ near a geodesic length $L_\kappa$ that $\phi_b$ can be chosen  linear in $t$ with  local coordinates on $M_b$ so that
\begin{equation}\label{phi linear in t}
    \phi(b,t,x,x,\xi)  = -|\xi| t + \psi(b,x,\xi),
\end{equation}
with 
\begin{equation}\label{psi on Zb}
    \psi(b,x,\xi) = L_\kappa(b)  \ \ \ {\rm for} \ (x,\xi) \in Z_b.
\end{equation}
Stationary phase approximation methods then give the singularity structure of $\Tr_\tmb \Us$ as in \cite{DG},  \cite{GM}, \cite{C}.  For, for $t$ near to $L_\kappa $ the integral  \eqref{kernel Ff diag b}  is an asymptotic sum of distributions over $j\geq 0$ of terms
$$\int_{M_b}  \int_{\R^q} e^{-i|\xi| t}e^{-i\psi(b,x,\xi)}\, {\bf w}_{\mu_d-j}(b,t,x,x,\xi)_{[d]} \ dx_{g_b} d\xi $$ $$= \int_0^\infty r^{q-1+\mu_d-j}\,e^{-ir t} \
\int_{M_b\x S^{q-1}} e^{-ir\psi(b,x,\eta)}\, {\bf w}_{\mu_d-j}(b,t,x,x,\eta)_{[d]} \ dx_{g_b} d\eta\,dr
$$
which is a localization of 
\begin{equation}\label{polar intmb}
   \int_0^\infty r^{q-1+\mu_d-j}\,e^{-ir t} \
\int_{S^*_\tmb/B} e^{-ir\psi(m,\eta)}\, {\bf w}_{\mu_d-j}(t,m,\eta)_{[d]} \,dr
\end{equation}
which is asymptotically via  v-SPA \eqref{v-SPA0} for $r\to\infty$
\begin{eqnarray}
   & \sim & \int_0^\infty r^{q-1+\mu_d-j}\,e^{-ir t} \ \sum_{p=1}^J \left(\frac{2\pi}{r}\right)^{\frac{n_p}{2}}\, e^{i\frac{\pi}{4}{\rm sgn}(\Zz_p)}\, e^{ir\psi(\Zz_p)} \sum_{l\geq 0} {\bf \alpha}_{\kappa,\mu_d-j,l,p} \, r^{-l} \ dr \notag\\
  &=&   \sum_{p=1}^J  \sum_{l\geq 0}  {\bf \alpha}_{\kappa,\mu_d-j,l,p} \,  \int_0^\infty r^{q+\mu_d-j -\frac{n_p}{2}-l -1} e^{-ir(t-L_\kappa)} \ dr\notag\\
 &\sim&   \sum_{p=1}^J  \sum_{l\geq 0}  {\bf \alpha}_{\kappa,\mu_d-j,l,p}\,(t- L_\kappa + i0)^{-\gamma_{\kappa,d,p,j,l}} \log(t- L_\kappa + i0)^{\mathfrak{H}(\gamma_{\kappa,d,p,j,l})} \label{sing exp basic}
\end{eqnarray}
where
\begin{equation}\label{alphalp}
    {\bf \alpha}_{\kappa,\mu_d-j,l,p} := \int_{\Zz_p/B} \omega_{\mu_d -j,l,p} \ \in\Aa^d(B)
\end{equation}
with $\omega_{\mu_d-j,l,p}$ derivatives along the fibres of $\Zz_p$ of ${\bf w}_{\mu_d-j}(t,m,\eta)_{[d]}$,  where $\mathfrak{H}$ is the Heaviside function $$\mathfrak{H}(y) = \l\{   \begin{array}{ll}
   1  & y  \leq 0 \\
  0 & y > 0 
\end{array}\r.,$$ 
where the constants have been absorbed into $\omega_{l,p}$, and where
\begin{equation}\label{gamma}
\gamma_{\kappa,d,p,j,l} := q+(\mu_d - j)-\frac{n_p}{2}-l.    
\end{equation}
Here, $\mu_d = k -\sum_{i=1}^k (mr_i  +  \sum_{n=1}^{r_i} \nu_{i,n})$ is the order of the summand $\Ff_{t,\Pp_k,d}$ of $\l(e^{-it\Qs^{-\frac{1}{m}}}\r)_{[d]}$ in \eqref{form deg d higher wave}. 
Thus, overall, for $t$ near  $L _k:B\to\R^1$ there is a singularity expansion 
\begin{equation}\label{full vertical wave trace}
    \Tr_\tmb\l(e^{-it\Qs^{-\frac{1}{m}}}\r)  \sim \sum_{k,d=0}^{\dim B}\sum_{\Pp_{k,d}} \sum_{p=1}^J  \sum_{j,l\geq 0}  
a_{\kappa,\mu_d-j,l,p}\,(t- L_\kappa + i0)^{-\gamma_{\kappa,d,p,j,l}} \log(t- L_\kappa + i0)^{\mathfrak{H}(\gamma_{\kappa,d,p,j,l})}
\end{equation}
where $\Pp_{k,d}$ are partitions of $d$, as in \propref{form wave operator}, and where
$$a_{\kappa,\mu_d-j,l,p} \in \Aa^d(B)$$
is a sum of differential form terms of type \eqref{sing exp basic}, \eqref{alphalp}. \\

Here, we have implicitly overlooked dependence of the symbol ${\bf w}_{\mu_d-j}(t,m,\eta)_{[d]}$ on $t$, and likewise the factor $(-it)^k$ in \eqref{Cauchy Qs d 2}. To take account of the first of these,  asymptotically, replace ${\bf w}_{\mu_d-j}(t,m,\eta)_{[d]}$  locally around $t=L_k$ by its Taylor expansion in powers of $t-L_\kappa$ to a given order, then  use
$\partial_r^k e^{-ir (t-L_\kappa)} = (-i(t-L_\kappa))^k e^{-ir (t-L_\kappa)}$  and partial integration in $r$ to thereby adjust the powers of $r$ in \eqref{sing exp basic}  (downwards), and do similarly for the powers of $t$ via $(-it)^k e^{-ir(t-L_\kappa)}  = ( \partial_r^ke^{-irt})e^{irL_\kappa}$. In both cases, this mixes summands slightly but the overall singularity structure \eqref{full vertical wave trace} does not change. \\

\eqref{full vertical wave trace} may be tidied up a little by taking $v=j+l$, $$a_{\kappa,\mu_d,v,p} = \sum_{j+l=v}a_{\kappa,\mu_d-j,l,p}\in \Aa^d(B)$$ and
\begin{equation}\label{gamma 2}
    \gamma_{d,p,v} := q+\mu_d-\frac{n_p}{2}-v
\end{equation} so for $t$ near $L_\kappa$
\begin{equation}\label{full vertical wave trace 2}
    \Tr_\tmb\l(e^{-it\Qs^{-\frac{1}{m}}}\r)  \sim \sum_{k,d=0}^{\dim B}\sum_{\Pp_{k,d}} \sum_{p=1}^J  \sum_{v\geq 0}  
a_{\kappa,\mu_d,v,p}\,(t- L_\kappa + i0)^{-\gamma_{d,p,v}} \log(t- L_\kappa + i0)^{\mathfrak{H}(\gamma_{d,p,v})}.
\end{equation}
There may be  infinitely many terms in \eqref{full vertical wave trace}, \eqref{full vertical wave trace 2} which do not have a log factor $\log(t-L_\kappa+i0)$ even though only finitely many for each $j$ (or $v$) since $q$ is the only  positive integer summand contributing to $\gamma_{\kappa,d,p,v}$, $\mu_d$ is a  negative integer (it is zero if and only if $d=0$) and in view of \eqref{max order m} $0\geq \mu_d \geq k - 2m\sum_{i=1}^k r_i.$
Collecting together coefficients of $(t- L_\kappa + i0)^l$, \eqref{full vertical wave trace 2} may be rewritten in the (schematically simpler) form
\begin{equation}\label{full vertical wave trace 3}
    \Tr_\tmb\l(e^{-it\Qs^{-\frac{1}{m}}}\r)  \ \sim \ \sum_{\kappa\in\Z}  \sum_{j\geq 0}  
{\bf a}_{\kappa,j}\,(t- L_\kappa + i0)^{-q + j} \hskip 25mm 
\end{equation}
$$\hskip 29mm + \ \sum_{\kappa\in\Z}  \sum_{l\geq 0}  
{\bf a}^{\,\prime}_{\kappa,l}\,(t- L_\kappa + i0)^{-q + l} \log(t- L_\kappa + i0),$$
with
$$ {\bf a}_{\kappa,j}  := \sum_d {\bf a}_{\kappa,j, [d]}, \ {\bf a}^{\,\prime}_{\kappa,l}  := \sum_d {\bf a}^{\,\prime}_{\kappa,l, [d]}\ \in\Aa^*(B),$$
which is \eqref{fibrewise Poisson summation}.\\

\subsection{The singularity of $\Us(t,m,m)$ at $t=0$:}

Computations at $L_0:=0$ are facilitated by the `big singularity' $\Zz$ of \eqref{Zz} then consisting of all of $S^*_\tmb$,  with codimension $n=0$ fibres. For $d=0$ matters are simplified further by  $\mu_0 =0.$
Assuming  for $L\neq 0, d=0$  the fibre codimension $n_p \geq 2q-2$, then by \eqref{full vertical wave trace 2}   $\Tr_\tmb(e^{-it\Qs^{-\frac{1}{m}}})_{[0]} $ for $t$ near $0$ is asymptotically 
\begin{equation}\label{wave invariants at 0}
    \sum_{j\geq0}  {\bf a}_j\,(t+ i0)^{-q+j} + \sum_{j\geq 0} {\bf a}^{\,\prime}_j\,(t + i0)^j \log(t + i0)
\end{equation} 
and for $t$ near $L_\kappa\neq 0$ asymptotically 
$${\bf a}_{-1,L}\,(t- L_\kappa + i0)^{-1} + \sum_{v\geq 0} {\bf a}_{v,L}\,(t- L_\kappa + i0)^v \log(t- L_\kappa+ i0).$$
When $\Qs$ is a vertical differential operator there are then no log-terms in the $d=0$ component $({\bf a}_{v,0})_{[0]}=0$ and so we re-obtain \eqref{wave trace form zero part} and \eqref{wave Poisson summation}.\\

We begin with the general identification between the wave invariants \eqref{wave invariants at 0} at $L_0=0$ and the heat trace invariants at $t=0$ for a vertical elliptic differential operator
$\Qs = \Qs_{[0]} + \Qs_{[1]} + \cdots + \Qs_{[\dim B]}$
 with each $\Qs_{[d]}\in \Psi^{\nu_i}(M/B)_{[d]}$ of fixed positive integer order $\nu_d$ and 
\begin{equation}\label{Q0D}
    \Qs_{[0]}=\D_\tmb
\end{equation}
a vertical metric Laplacian (considerations are similar for any vertical positive self-adjoint elliptic differential operator of positive integer order). $\D=\D_\tmb$ is vertically parameter elliptic with Agmon angle $\pi$,  with 
vertical principal   symbol ${\bf q}_2(b, x,\xi) = \langle\xi,\xi\rangle_{g_b} = g^{ij}_b(x)\xi_i\xi_j = |\xi|^2 $ and 
 \begin{equation}\label{q0-la}
     {\bf q}_2 - \la I\in \Ci(T_\tmb\backslash\{0\},p^*(\End(E)))
 \end{equation}
 invertible for $\la\in\C\bs\R_+$.
 As seen in \eqref{res exp d}, 
 $$(\Qs - \la \Is)\ii = $$ 
\begin{equation*}
=   (\D -
\la \Is)\ii  + \sum_{d=1}^{\dim B}\sum_{p_1+\ldots +p_k = d}  (-1)^k (\D -
\la \Is)\ii \Qs_{[p_1]}(\D - \la \Is)\ii \ldots \Qs_{[p_k]}(\D -
\la \Is)\ii,  
\end{equation*} 
which with the assumption $\nu_d \leq 2$
has operator norm  $\| (\Qs - \la \Is)\ii\|^{(l)}_{M/Z} = O(|\la|\ii)$ 
as $\la\to\infty$  in an open sector around the ray $R_\pi$ emanating from the origin which coincides with the negative real axis. Holomorphic functional calculus  constructs the complex power as in \eqref{e:powers 3}
\begin{equation}\label{power zeta}
    \Qs_\pi\si = \int_{C}\la_\pi\si\, (\Qs - \la\Is)\ii \, \dl\la
\end{equation}
 which is trace class for $\re(s)>>0$ and whose trace $\Tr_\tmb(\Qs_\pi\si) \in \Aa^*(B)$ for such $s$ admits a meromorphic continuation
$Z_\Qs(s)$ to all of $\C$ with pole singularity structure 
\begin{equation}\label{zeta form singularity structure}
 \G(s)Z_\Qs(s) \sim
 \end{equation}
 $$-\sum_{j=-\dim B -1}^{-1}
  \frac{\str\left((\Pi\cdot\Qs_+\cdot\Pi)^{-j-1}\right)}{s - j -
  1} \ \ + $$
 $$\sum_{d=0}^{\dim B}\left(\sum_{j\geq 0} \ 
\frac{{\bf b}_{j,d}}{s + \frac{- n + j- \nu_d}{2}}\right. +
\left. \sum_{l\geq 0}\frac{{\bf b}^{'}_{l,d}}{(s + l - 1)^2} +
\frac{{\bf b}^{''}_{l,d}}{s + l - 1}\right) \ ,$$ 
with $\nu_d := \sum_i\nu_i$ and residues $${\bf b}_{j,d}, {\bf b}^{'}_{l,d}, {\bf b}^{''}_{l,d}\in \Aa^{d}(B).$$ 
This holds, in fact, is for any vertical pseudodifferential $\Qs$ (with $\Qs_{[0]}$ elliptic order 2), but  the restriction to differential operators implies the globally determined terms ${\bf b}^{''}_{l,d} =0$ all vanish.  $\Pi_b$ is the $L^2$ orthogonal projection to $\Ker(\D_b) $, defining a smoothing operator for each $b\in B$, and we assume these define a smooth family 
$\Pi\in\Aa(M/B,\Psi^{-\o}(\Ee))$, whose ranges define thus a finite-rank
 vector bundle $\Ker(\D) \to B$.\\

Equivalently, via Mellin transforms, the vertical heat operator 
\begin{equation}\label{heat operator}
e^{-t\Qs} = \int_{\Cc}
e^{-t\la}\,(\Qs - \la\Is)\ii \ \dl\la
\end{equation}
for a contour $\Cc$  coming in on a ray with
argument in $(0,\pi/2)$, encircling the origin, and leaving on a
ray with argument in $(-\pi/2,0)$, has an asymptotic
expansion as $t\to 0+$
\begin{equation}\label{heat trace singularity structure}
\Tr_\tmb(e^{-t\Qs})  \sim  \sum_{d=0}^{\dim B}\left(\sum_{j\geq 0}
 {\bf b}_{j,d} \, t^{\frac{ - n +j -\nu}{2} -1} + \sum_{l\geq
0} ({\bf b}^{'}_{l,d}\log t + {\bf b}^{''}_{l,d}) \,
t^{l-1}\right) \
\end{equation}
with the same  differential form coefficients ${\bf b}_{j,d}, {\bf b}^{'}_{l,d},
{\bf b}^{''}_{l,d}\in \Aa^{d}(B)$ as \eqref{zeta form singularity structure}. The coefficients  ${\bf b}_{j,d}, {\bf b}^{'}_{l,d}$ are locally determined (can be computed from just finitely many homogeneous terms and their derivatives of the symbol of $\Qs\si$), whilst the ${\bf b}^{''}_{l,d}$ are globally determined. If $\Qs$ is a vertical differential operator then there are no log terms -- ${\bf b}^{'}_{l,d}=0$. 
 One may equally use the Duhamel principle to give a Volterra series construction of $e^{-t\Qs}$ \cite{BGV}, as here for the vertical wave operator \eqref{Cauchy Qs}, but \eqref{heat operator} is better suited to computations of the singularity structure \eqref{heat trace singularity structure}. The constructions are equivalent with the  Duhamel principal variational formula for $\partial_s e^{-t\Qs(s)}$ written in terms of the contour integral definition as $\int_{\Cc}
e^{-t\la}\,(\Qs(s) - \la\Is)\ii \cdot \partial_s \Qs(s)\cdot (\Qs(s) - \la\Is)\ii \ \dl\la.$\\

The above statements on the singularity structures are immediate from work of Grubb and Seeley \cite{GS1}, \cite{GS2}, which for the context here were looked into more closely in \cite{Sc}. The identification, here, of the vertical wave trace invariants \eqref{wave invariants at 0} at $t=0$ with the coefficients in \eqref{zeta form singularity structure} and \eqref{heat trace singularity structure} is based on the property that at $t=0$ the wave forms, also, can be computed microlocally via holomorphic functional calculus in a similar way to \eqref{power zeta} and \eqref{heat operator} (not for $L_\kappa \neq 0$). Let $ \Pi_0$ denote the vertical smoothing operator projecting orthogonally onto $\Ker(\sqrt{\D_b})$ at $b\in B$. For small enough $|t|<\delta < |L_{\pm 1}|$ 
\begin{equation}\label{wave operator contour}
e^{-it\sqrt{\Qs}} - \Pi_0 = \Qs\int_{\Lambda} \frac{e^{-it\mu}}{\mu^2} ( \sqrt{\Qs} - \mu)\ii \ \dl\mu
\end{equation}
with contour $\La$  comprising a ray coming in from infinity along the line parallel to the $x$-axis at height $\beta>0$ as far as $x= \varepsilon>0$  with $\varepsilon<\la_1 := \inf_{b\in W\< B}\{\la_1(b)\}$ and $\la_1(b)$ the smallest non-zero eigenvalue of $\sqrt{\D_b}$, then dropping down to the parallel line at height $-\beta<0$ below the $x$-axis and back out to infinity along that ray. The coordinate chart $W$ is chosen to ensure $\varepsilon>0$ and use fibrewise coordinates on $M$ over $W$, since our computations will be local this is enough. The weight $\mu^{-2}$ is there to ensure convergence of the integral and is counterbalanced by the factor $\Qs -\Pi_0= \sqrt{\Qs}^2$ (since $\Qs$ is positive self-adjoint) -- indeed, if  $\Lambda$ is replaced by the finite contour formed by closing off 
$\Lambda$ with a right-hand vertical side at $x=R>\varepsilon> 0$ chosen disjoint from the spectrum of $\sqrt{\Qs}$, to form a  positively oriented rectangle  $\Lambda_R$, then 
\begin{equation}\label{Lambda R}
    \Qs\int_{\Lambda_R} \frac{e^{-it\mu}}{\mu^2} ( \sqrt{\Qs} - \mu)\ii \ \dl\mu = \int_{\Lambda_R} e^{-it\mu} ( \sqrt{\Qs} - \mu)\ii \ \dl\mu. 
\end{equation}
This is precisely what can be done for local symbol computations with $R$ large enough so $\La_R$  contains $|\xi_b|^2$ for $b\in W$ and
\eqref{q0-la} is invertible so the vertical resolvent symbol for $( \sqrt{\Qs} - \mu)\ii $  in reduced $x$-format has $d$-form component a polyhomogeneous symbol
\begin{equation}\label{sqrt r expansion}
    \sqrt{{\bf r}}[\mu](b,x, \xi)_{[d]}  \sim \sum_{j\geq 0}  \sqrt{{\bf r}}_{\, \nu - (k+1) - j}[\mu](b,x,\xi)_{[d]}
\end{equation}
with $\sqrt{{\bf r}}_{\,\nu - (k+1) - j}[\mu](b,x,\xi)_{[d]}$ homogeneous for $t>1$ in $(\mu,\xi)$ of order $\nu - (k+1) - j$, allowing, for small $t$, the wave `amplitude' at form degree $d$ to be constructed microlocally as 
\begin{equation}\label{wave Lambda R}
  {\bf w}(b,x,\xi) = \int_{\Lambda_R} e^{-it\mu} \sqrt{{\bf r}}[\mu](b,x, \xi)_{[d]} \ \dl\mu.
\end{equation}
Then, with $|dx|$ local Lebesgue measure on $M_b$
\begin{equation}\label{local density b}
\int_{\R^q} {\bf w}(b,x,\xi)\,|dx| 
\end{equation}v
defines a (local representative for a) globally defined Riemannian density on $M_b$ which can be integrated along that fibre. \eqref{wave Lambda R} captures not just the symbol/amplitude but also a portion of the phase (being FIO not just $\pdo$ calculus);
from \cite{DG} \S2 it is noted that for such $t$ there is a local coordinate system in which the phase \eqref{kernel Ff} takes the form
\begin{equation}
    \phi(b, t, x,y,\xi) =i\langle x-y, \xi\rangle + t|\xi|_{g_b(x)}.
\end{equation}
and it is the $t$ portion of this that arises from \eqref{wave Lambda R} when computing the local Schwartz kernel of $e^{-it\sqrt{\Qs}}$ from \eqref{wave operator contour}. The essential consequence is that asymptotically, and modulo smoothing operators, for $t$ near 0 and $b\in W$
\begin{equation}\label{v wtr symbol}
    \Tr_\tmb(e^{-it\sqrt{\Qs}} )(b)  \sim \int_{M_b} \int_{\R^q}  \int_{\Lambda_R} e^{-it\mu} \,\sqrt{{\bf r}}[\mu](b,x, \xi)_{[d]} \ \dl\mu.
\end{equation}
 \eqref{wave Lambda R} in this way forms a vertical wave operator analogue of the complex power symbol ${\bf q}_\pi^{1/m}$ in \eqref{1/m powersexpansion}. For the latter, for $\re(s)>0$  the complex power 
\begin{equation}\label{power sq rt zeta}
    (\sqrt{\Qs})_\pi\si = \int_{C}\la_\pi\si\, (\sqrt{\Qs} - \la\Is)\ii \, \dl\la
\end{equation}
 extends to all $s\in\C$ as instructed by Seeley \cite{Se}, and as done for the construction of $\Qs^{1/m}$ above,
 for $\re(s)  > - N $   by  $\sqrt{\Qs}^{\,-s} :=  \sqrt{\Qs}^N \cdot  \sqrt{\Qs}_\pi^{\,-N-s}$ independently of the choice of $N$ -- a similar mechanism is being used in the definition \eqref{wave operator contour} (one could equally well replace $\mu^2$ by $\mu^N$ and then $\sqrt{\Qs}^{\,2}$ by $\sqrt{\Qs}^{\,N}$ which is an important freedom when taking derivatives) -- and $\sqrt{\Qs}^{\,-s}$  so defined is a vertical $\pdo$ with local vertical symbol  for all $s$ modulo smoothing symbols
 \begin{equation}\label{zeta sqrt symbol}
  {\bf z}\si(b,x,\xi) = \int_{\Lambda_R} \la_\pi\si\,  \sqrt{{\bf r}}[\la](b,x, \xi)_{[d]} \ \dl\la. 
\end{equation}
From this formula the local residues  of the singularity structure of the  meromorphically continued zeta trace form $Z_{\sqrt{\Qs}}(s)$
may be computed by integrating the  globally defined Riemannian density on $M_b$ 
\begin{equation}\label{local zeta density b}
\int_{\R^q} {\bf z}\si(b,x,\xi)\ \dl\xi\,|dx|,
\end{equation}
that is, by computing asymptotically
\begin{equation}\label{local zeta sing structure}
Z_{\sqrt{\Qs}}(s)(b) \sim \int_{M_b} \int_{\R^q} {\bf z}\si(b,x,\xi)\ \dl\xi\,|dx|,
\end{equation}
and given $Z_\Qs(s) = Z_{\sqrt{\Qs}}(2s)$ then hence also the residues  $b_{j,d}, b^{'}_{l,d}, b^{''}_{l,d}\in \Aa^{d}(B)$ of $\Qs$ in \eqref{zeta form singularity structure}, as is done in Duistermaat-Guillemin Cor. 2.2 \cite{DG} for the case $d=0$.\\

Before doing that, we take a moment to comment on the choice of contour $\Lambda$ in \eqref{wave operator contour}.  For  the complex powers \eqref{power zeta}, \eqref{power sq rt zeta} one commences with a simple closed key-hole contour $\G_R\sqcup S_R$ with $S_R$ a circle  of radius $R>>0$ (oriented positively and $R$ disjoint from the spectrum of $\sqrt{\D}$) and with $\G_R$ the portion of $\G$ contained inside the disc of radius $R$.
For $\re(s)>0$, crucially,  the integral along $S_R$ tends to $0$ and is convergent along $\G_R$ as $R\to\infty$,  leading via the Cauchy theorem to   \eqref{power zeta}, \eqref{power sq rt zeta}. Similarly, for the heat operator \eqref{heat operator},
noting that it is essential that the contour lie in the right-half plane.
The difficulty for the wave operator  is that one of $e^{\pm it\mu}$ becomes exponentially unbounded as $\mu \to \infty$ along any ray that has non-zero gradient, overwhelming the polynomial $\mu$ convergence there of the resolvent $(\Qs-\mu)\ii$, while along the connecting closing arc ($S_R$ in the above instance) the integral then does not tend to 0 as $R\to\infty$. $\Lambda_R$ is chosen to avoid these limitations, along the right-hand side of the rectangle $e^{\pm it\mu}$ is bounded by $e^{\delta\beta}$, for $\beta$ is constant and $|t|< \delta <|L_{\pm 1}|$, while   $(\Qs-\mu)\ii$ is (only) $O(1)$ along the horizontal rays, and hence the  factor $\mu^{-2}$ in \eqref{wave operator contour}, corrected by composition with $\sqrt{\Qs}^{\,2} = \Qs$, ensures  that as $R\to\infty$ the  integral is convergent  and tends to 0 along the right-hand edge in the limit. \\ 

To compare the locally determined terms in the vertical zeta trace singularity expansion \eqref{zeta form singularity structure}  with those of the wave trace  \eqref{wave invariants at 0} at $t=0$ they may be computed using \eqref{zeta sqrt symbol} and \eqref{wave Lambda R}.\\

We start with the case $d=0$ - that is, the classical pointwise on $B$ wave and zeta traces identifications of \cite{DG}. First up is the expansion $\sqrt{{\bf r}}[\mu](b,x, \xi)_{[0]}  \sim \sum_{j\geq 0}  \sqrt{{\bf r}}_{\, 1 - j}[\mu](b,x,\xi)_{[0]}$, written for brevity in this case 
$$\sqrt{r}[\mu](x, \xi)  \sim \sum_{j\geq 0}  \sqrt{r}_{\, 1 - j}[\mu](x,\xi)$$
of the symbol of the resolvent $(\sqrt{\D} - \mu)\ii$.
\begin{lem}\label{res sqrt symbols}
$$\sqrt{r}_{-1}[\mu](x,\xi) = (|\xi|_{g(x)}- \mu)^{-1}$$
and for $j>0$
    \begin{equation}\label{sqrt r -1-j}
     \sqrt{r}_{ - 1 - j}[\mu](x,\xi) = \sum_{l=[j/2]}^{2j} \frac{(2l-1)!}{2^{2l-1} (l-1)!}\, \sqrt{r}_{-1}[\mu](x,\xi)^{1+2l}  \, r_{l,j}(x,\xi) 
\end{equation}
with $r_{l,j}(x,\xi)$ a homogeneous polynomial in $\xi = (\xi_1,\ldots, \xi_q)$ of degree $2l-j$.
\end{lem}
{\it Proof:}  Consider first the case $\mu=0$: then the homogeneous summands $\sqrt{r}_{\, 1 - j}[0](x,\xi)$ are according to \eqref{cx symbol}, \eqref{e:powersexpansion} computed from the resolvent symbols  $r[\la](x,\xi):= {\bf r}[\la](b,x,\xi)_{[0]} \sim \sum_{j\geq 0} r_{\, -2 - j}[\la](x,\xi)$ of $\D$ by 
$$\sqrt{r}_{ - 1 - j}[0](x,\xi) =  \int_{C_R}\la_\pi^{-\frac{1}{2}}\, r[\la]_{-2-j}(x,\xi) \ \dl\la.$$
Since 
$$r_{-2}[\la](x,\xi) = (|\xi|^2_{g(x)} -\la)\ii$$
is scalar valued, then  its quasi-homogeneous components have the form
\begin{equation}\label{scalar res symbol}
    r[\la]_{-2-j}(x,\xi) = \sum_{l=[j/2]}^{2j} r_{-2}[\la](x,\xi)^{l+1}\, r_{l,j}(x,\xi)
 \end{equation}
 with $r_{l,j}(x,\xi)$ a sum of products of derivatives of the symbol components of $\D$ and a homogeneous polynomial in $\xi = (\xi_1,\ldots, \xi_q)$ of degree $2l-j$. Hence
\begin{eqnarray*}
    \sqrt{r}_{ - 1 - j}[0](x,\xi) & = &  \sum_{l=[j/2]}^{2j}  \int_{C_R}\la_\pi^{-\frac{1}{2}}r_{-2}[\la]^{l+1}\, \dl\la \ r_{l,j}(x,\xi) \\ 
    & = &  \sum_{l=[j/2]}^{2j}  \int_{C_R}\frac{(-1)^l}{l!}\partial_\la^l \l(\la^{-\frac{1}{2}}\r)\cdot r_{-2}[\la]\ \dl\la \ r_{l,j}(x,\xi) \\
     & = &  \sum_{l=[j/2]}^{2j} \frac{1}{l!}\prod_{n=1}^l \l(n-\frac{1}{2}\r) \int_{C_R}\la^{-\frac{1}{2}-l} \, r_{-2}[\la]\ \dl\la \ r_{l,j}(x,\xi) \\
     &=&  \sum_{l=[j/2]}^{2j} \frac{(2l-1)!}{2^{2l-1} l! (l-1)!}\, |\xi|_{g(x)}^{-1-2l}  \, r_{l,j}(x,\xi). 
\end{eqnarray*}
from which the \eqref{sqrt r -1-j} is immediate.                 \hfill $\Box$  
\vskip 3mm 

\begin{prop}\label{local sing wave density}
Asymptotically at form degree $d=0$ the vertical wave trace pointwise 
$\Tr_\tmb(e^{-it\sqrt{\Qs}} )_{[0]}$ is the integral over $M_b$ of  the distribution valued density 
\begin{equation}\label{wave trace d=0 sing str}
\sum_{j\geq 0} \l(\sum_{l=[j/2]}^{2j} \frac{(2l-j + q-1)! }{2^{2l}(l!)^2} \,  c_{l,j}(b,x)\r) d{\rm vol}_{g_b(x)} \cdot (t-i0)^{-q+j},
\end{equation}
where $d{\rm vol}_{g_b}$ is the Riemannian density on $M_b$ and $c_{l,j}(b,x)$ a sphere integral which is zero for $j$ odd.
\end{prop}
{\it Proof:} Mindful of \eqref{local density b}, \eqref{v wtr symbol}, we have that $\int_{\R^q} {\bf w}(x,\xi)_{[0]} \ \dl\xi\, |dx|$ is by \lemref{res sqrt symbols} asymptotically equal to
\begin{equation}\label{wave sing sum}
    \sum_{j\geq 0} \sum_{l=[j/2]}^{2j}  \frac{(2l-1)!}{2^{2l-1} l! (l-1)!}\, \int_{\R^q}  \int_{\Lambda_R} e^{-it\mu} \, \sqrt{r}_{-1}[\mu](x,\xi)^{1+2l}   \ \dl\mu \ \, r_{l,j}(x,\xi) \ \dl\xi \,|dx|.
\end{equation} 
The distribution valued integral in the summands is
$$ \int_{\R^q}  \int_{\Lambda_R} \frac{1}{(2l)!}\,\partial_\mu^{2l} \l( e^{-it\mu}\r) \cdot  \sqrt{r}_{-1}[\mu](x,\xi) \ \dl\mu \ \, r_{l,j}(x,\xi) \ \dl\xi $$
which with $\alpha_l:= (-1)^l/(2l)!$
\begin{eqnarray*}
      &=& \alpha_l t^{2l}\int_{\R^q}   e^{-it|\xi|_{g(x)}} \, r_{l,j}(x,\xi) \ \dl\xi \\
      &=& \alpha_l t^{2l}\int_{\R^q}   e^{-it|\xi|} \, r_{l,j}(x,\xi) \ \dl\xi \,\sqrt{\det g(x)}\\
      &=& \alpha_l t^{2l}\int_0^\infty   e^{-itr} \, r^{2l-j + q-1} \ dr\, \int_{S^{q-1}} r_{l,j}(x,\eta) \ \dl_S\eta\,\sqrt{\det g(x)} \\
      &=&\alpha_l (2l-j + q-1)! \,t^{2l} (t-i0)^{-(2l-j+q)} \, \int_{S^{q-1}} r_{l,j}(x,\eta) \ \dl_S\eta\,\sqrt{\det g(x)}\\
      &=& \alpha_l (2l-j + q-1)! \, \int_{S^{q-1}} r_{l,j}(x,\eta) \ \dl_S\eta\,\sqrt{\det g(x)}\cdot (t-i0)^{-q+j},
\end{eqnarray*}
using the $\xi$ homogeneity of $r_{l,j}$ with $\dl_S\eta$ sphere measure, and a change of variable which  modifies $r_{l,j}(x,\xi)$ by factors $g_{ij}(x)$ but does not affect its $\xi$ homogeneity. Since $r_{l,j}(x,\xi)$ is polynomial of degree $2l-j$,  and thus $l>j/2$, the exponent in $r^{2l-j + q-1}$ is always positive and so there are no log factors.  
Substituting into \eqref{wave sing sum} we obtain  \eqref{wave trace d=0 sing str} with $c_{l,j}(b,x)$ the above sphere integral computed with respect to $g(x) = g_b(x)$ multiplied by (the independent of $l$) power of $i$. When $j$ is odd then so also is the homogeneity of $r_{l,j}$ and hence vanishing of the oscillatory integral, or the sphere integral, then follows by change of variable $\xi \mto -\xi$. \hfill $\Box$  
\vskip 3mm 

\eqref{wave trace d=0 sing str} provides the more precise formulae for the coefficients in the vertical wave singularity expansion \eqref{wave invariants at 0} at $t=0$ of form degree $d=0$
\begin{equation}\label{formulae for ajb d=0}
    {\bf a}_j(b)_{[0]}  =  \sum_{l=[j/2]}^{2j}  \frac{(2l-j + q-1)! }{2^{2l}(l!)^2} \, \int_{S^*M_b} r_{l,j}(x,\eta) \ \dl_S\eta\, \ d{\rm vol}_{g_b(x)},
\end{equation} 
$$ {\bf a}_{2l+1}(b)_{[0]}  =0, \hskip 12mm {\bf a}_{j,0}(b)_{[0]} =0.$$
\vskip 3mm 

On the  other hand, for the zeta trace one has from \eqref{zeta sqrt symbol}, \eqref{sqrt r -1-j}, 
\begin{eqnarray*}
  {\bf z}\si(x,\xi) & \sim & \sum_{j\geq 0} \sum_{l=[j/2]}^{2j}\frac{(2l-1)!}{2^{2l-1} (l-1)!}  \int_{\Lambda_R} \la_\pi\si\, \sqrt{r}_{-1}[\mu](x,\xi)^{1+2l} \   \dl\mu \, r_{l,j}(x,\xi)  \\
  & = &  \sum_{j\geq 0} \sum_{l=[j/2]}^{2j}\frac{(2l-1)!}{2^{2l-1} (l-1)!}  \int_{\Lambda_R} \frac{1}{(2l)!}\,\partial_\mu^{2l} \l(\la_\pi\si\r)\, \sqrt{r}_{-1}[\mu](x,\xi) \   \dl\mu \, r_{l,j}(x,\xi)  \\
  & = &   \sum_{j\geq 0}\sum_{l=[j/2]}^{2j}\frac{s(s+1)\cdots(s+2l-1)}{2^{2l} (l!)^2}  \int_{\Lambda_R} \la_\pi^{-s-2l}\, \sqrt{r}_{-1}[\mu](x,\xi) \   \dl\mu \, r_{l,j}(x,\xi)  \\
   & = &  \sum_{j\geq 0} \sum_{l=[j/2]}^{2j}\frac{s(s+1)\cdots(s+2l-1)}{2^{2l} (l!)^2} \  |\xi|_{g_b(x)}^{-s-2l} \ r_{l,j}(x,\xi) 
\end{eqnarray*}
For simplicity here we shall not consider globally determined residues, so work modulo a smoothing operator error by taking a bump function $\chi(|\xi|) =0$ for $|\xi|< 1/3$, $\chi(|\xi|) =1$ for $|\xi|> 2/3$ and contemplate for $\re(s)>>0$
$$\int_{\R^q} \chi(|\xi|)\,{\bf z}\si(x,\xi) \ \dl\xi \hskip 8cm$$ 
\begin{eqnarray*}
 & \sim &  \sum_{j\geq 0}  \sum_{l=[j/2]}^{2j}\frac{s(s+1)\cdots(s+2l-1)}{2^{2l} (l!)^2} \ \int_{\R^q}\chi(|\xi|)\, |\xi|_{g_b(x)}^{-s-2l} \ r_{l,j}(x,\xi) \ \dl\xi\\
 & = &   \sum_{j\geq 0} \sum_{l=[j/2]}^{2j}\frac{s(s+1)\cdots(s+2l-1)}{2^{2l} (l!)^2} \ \int_{\R^q} \chi(|\xi|)\,|\xi|^{-s-2l} \ r_{l,j}(x,\xi) \ \dl\xi\,\sqrt{\det g(x)}\\
 & \sim &   \sum_{j\geq 0} \sum_{l=[j/2]}^{2j}\frac{s(s+1)\cdots(s+2l-1)}{2^{2l} (l!)^2} \ 
 \int_1^\infty  r^{-s-j+q-1} \ dr \, \int_{S^{q-1}} r_{l,j}(x,\eta) \ \dl_S\eta\,\sqrt{\det g(x)}
\end{eqnarray*}
-- discarding the smoothing operator portion over the unit ball, which may give rise to globally determined residues corresponding to the primed coefficients in \eqref{zeta form singularity structure} ---
\begin{equation*}\label{sing str 1}
  \stackrel{\re(s)> q}{=}    \sum_{j\geq 0}\l( \sum_{l=[j/2]}^{2j}\frac{s(s+1)\cdots(s+2l-1)}{2^{2l} (l!)^2}  \, \int_{S^{q-1}} r_{l,j}(x,\eta) \ \dl_S\eta\,\sqrt{\det g(x)}\r)\cdot
   \frac{1}{s-q+j}
\end{equation*}
which extends meromorphically to $\C$, and setting $s=-q+j$ in the coefficients of the pole at that point gives from \eqref{local zeta sing structure} the singularity structure
$$Z_{\sqrt{\Qs_b}}(s)(b)_{[0]} = Z_{\sqrt{\D_b}}(s)(b)\sim $$
\begin{equation*}\label{local zeta sing structure d=0}
 \sum_{j\geq 0}\underbrace{\l( \sum_{l=[j/2]}^{2j}\frac{(q-j)(q-j+1)\cdots(q-j+2l-1)}{2^{2l} (l!)^2}  \,  \int_{M_b} c_{l,j}(b,x) d{\rm vol}_{g_b(x)}\r)}_{={\bf a}_ j(b)_{[0]}/\G(q-j)} \cdot  \frac{1}{s-q+j}
 \end{equation*}
with $a_ j(b)_{[0]}$ the coefficients \eqref{formulae for ajb d=0} of the wave trace expansion  \eqref{wave invariants at 0} at $t=0$ of form degree $d=0$. Since the sums $\sum_{\la\neq0}\la\si$ and $\sum_{\la\neq0}(\la^{\frac{1}{2}})^{-2s}$ over the real positive eigenvalues $\la$ of $\D$ are equal, then $Z_\D(s) = Z_{\sqrt{\D}}(2s)$, and hence
\begin{equation}\label{sing str Gamma zeta s d=0}
    \G(s) Z_{\D_b}(s)(b) \sim  \sum_{j\geq 0} {\bf a}_j(b)_{[0]}  \cdot  \frac{1}{s-\frac{q-j}{2}}.
\end{equation}
In view of the correspondence between \eqref{zeta form singularity structure} and \eqref{heat trace singularity structure}, one therefore has: 
\begin{prop}\cite{DG}\label{wave heat d=0}
    For $d=0$: the local wave trace invariants \eqref{wave invariants at 0} are the same as the locally determined heat trace invariants 
    $$  {\bf b}_j(b)_{[0]} = \frac{{\bf a}_j(b)_{[0]}}{2}.$$
\end{prop}
This is, of course, no more than the classical result of Duistermatt and Guillemin, but by an analogous analysis the result persists easily to an equality of the higher degree differential form coefficients of the  wave trace around $L_0=0$ with the residues of the higher zeta singularity  expansion. Concretely, from \eqref{local density b}, \eqref{sqrt r expansion} the wave trace invariants at $L_=0$ at differential form degree $d$ can be computed as a finite sum of 
asymptotic sums 
\begin{equation*}
\sum_{j\geq 0} \int_{M_b}\int_{\R^q} \int_{\Lambda_R} e^{-it\mu} \ \sqrt{{\bf r}}_{\, \nu - 2(k+1) - j}[\mu](b,x,\xi)_{[d]}\ \dl\mu |dx|
\end{equation*}
with each such summand deriving from a summand
\begin{equation}\label{R d k summand}
 \Rs_{d_1, \ldots, d_k |d}[\la] := (\Ps -
\la \Is)\ii \Qs_{[d_1]}(\Ps - \la \Is)\ii \ldots \Qs_{[d_k]}(\Ps -
\la \Is)\ii
\end{equation}
of $(\Qs-\la I)\ii$. From \eqref{d parametrix symbols}, the vertical symbol of $\Rs_{d_1, \ldots, d_k |d}[\la]$  has the form
$${\bf r}_{d_1, \ldots, d_k |d}[\la](x,\xi) = \sum_{j\geq 0}{\bf r}_{\nu - 2k -2-j}[\la](b,x,\xi) _{[d]}$$
with ${\bf r}_{\nu - 2k -2-j}[\la](b,x,\xi)_{[d]}$ quasi-homogeneous of degree $\nu - 2k -2-j$ and of differential form degree $d$ on $B$, where $\nu:= \nu_1+\cdots+ \nu_k$ with $\nu_i$ the positive integer order of $\Qs_{d_i}$. Generalising \eqref{scalar res symbol} one has:
\begin{lem}\label{resolvent symbols d}
Each  ${\bf r}_{\nu - 2k -2-j}[\la](b,x,\xi)_{[d]}$ is a finite sum of terms of the form  
    \begin{equation}\label{res degree d summands}
r_{-2}[\la](b,x,\xi)^{l+1}\, {\bf r}_{\nu+2l-2k-j}(b,x,\xi)_{[d]}
\end{equation}
with ${\bf r}_{\nu+2l-2k-j}(b,x,\xi)_{[d]}$ a homogeneous polynomial in $\xi$ of degree $\nu+2l-2k-j$, $l$ a positive integer $\geq j/2$ and, as before, $r_{-2}[\la](b,x,\xi) = (|\xi|^2_{g_b(x)} - \la)\ii$.
\end{lem}
{\it Proof:} From \eqref{d parametrix symbols} and \eqref{scalar res symbol} each homogeous summand is a sum of terms 
\begin{equation*}
 \l(\sum_{m_1=[j_1/2]}^{2j_1} r_{-2}[\la](x,\xi)^{m_1+1}\, r_{m_1,j_1}(x,\xi)\r)
\circ {\bf q}_{\nu_1 - l_1}(b,x,\xi)_{[d_1]}\circ\cdots 
\end{equation*}
$$\hskip 3cm \cdots\,\circ\,{\bf
q}_{\nu_k - l_k}(b,x,\xi)_{[d_k]}\circ  \l(\sum_{m_{k+1}=[j_{k+1}/2]}^{2j_{k+1}} r_{-2}[\la](x,\xi)^{m_k+1}\, r_{m_{k+1},j_{k+1}}(x,\xi)\r).$$
Since $\Qs_{[d_i]}$ is a vertical differential operator then the symbol ${\bf q}_{\nu_i - l_i}(b,x,\xi)_{[d_i]}$ is a homogeneous polynomial in $\xi$ of degree  $\nu_i - l_i$ (zero for $l_i>\nu_i$). Taking $x$-derivatives does not change this whilst $\xi$-derivatives only lowers the order/homogeneity. Similarly, $x$ and $\xi$ derivatives of
 $ r_{-2}[\la](x,\xi)^{m_1+1}r_{m_1,j_1}(x,\xi)$ gives sums of terms of the same type with $l$ increasing and the homogeneity of the, again $\xi$-polynomial, term  $r_{r_1,j_1}(x,\xi)$ decreasing accordingly. Since the symbol composition is a sum of products of $x$ and $\xi$ derivatives of such terms, then bearing in mind the order we reach the conclusion.
\hfill $\Box$  \vskip 3mm

Knowing the general form  \eqref{res degree d summands} is enough for our needs here and so we do not strive for greater accuracy, indeed the rest is a {\it fait accompli} differing from the $d=0$ case only in degrees of (quasi)homogeneity. Briefly, by the same computation as in the proof of \lemref{res sqrt symbols}, each  $\sqrt{{\bf r}}_{\nu - 2k -2-j}[\mu](b,x,\xi)$ is a finite sum of terms 
\begin{equation}\label{sqrt res-symbol degree d}
\sqrt{{\bf r}}_{\nu - 2k -1 - j}[\mu](b,x,\xi)_{[d]}=\sum_l \beta_{l,\nu, k,j} \, (|\xi|_{g_b(x)} - \mu )^{-2l-1}\, \, {\bf r}_{\nu-2k +2l-j}(b,x,\xi)_{[d]}
\end{equation}
with $\beta_{l,\nu,k,j}$ a rational function of $l,\nu,k,j$. Next, carrying through the evaluation in the proof of \propref{local sing wave density},  the contribution of the component \eqref{R d k summand} to the vertical wave trace  $\Tr_\tmb(e^{-it\sqrt{\Qs}} )_{[d]}(b)$ at $t=0$ evaluated at $b\in B$ is asymptotically the distribution valued density 
$$\sum_{j\geq 0} \sum_l \frac{\beta_{l,\nu,k,j}}{(2l)!}(2l-j+\nu-2k+ q-1)! \,  \int_{M_b} c_{l,\nu,k,j}(b,x)_{[d]} \ d{\rm vol}_{g_b(x)} $$
\begin{equation}\label{wave trace d sing str}
\x \  (t-i0)^{-q+j + 2k - \nu}
\end{equation}
\vskip 3mm
with $c_{l,\nu,k,j}(b,x) =\int_{S^{q-1}} {\bf r}_{\nu+2l-2k-j}(b,x,\eta)_{[d]}\ d_S\eta$. 
On the other hand, from \eqref{zeta sqrt symbol}, \eqref{local zeta sing structure} and using \eqref{sqrt res-symbol degree d} we have that $Z_{\sqrt{\Qs}}(s)(b)$ is asymptotically 
$$\sim \int_{M_b} \int_{\R^q} \int_{\Lambda_R} \la_\pi\si\, \sum_l \beta_{l,\nu, k,j} \, (|\xi|_{g_b(x)} - \mu )^{-l-1}\, \, {\bf r}_{\nu+2l-2k-j}(b,x,\xi)_{[d]} \ \dl\la \ \dl\xi\,|dx|$$
and by following through the same computation  \eqref{local zeta sing structure d=0} we obtain at form degree $d$ that the contribution coming from  the component \eqref{R d k summand} to the singularity structure of $Z_{\sqrt{\Qs_b}}(s)(b)_{[d]}$ is:
\begin{equation*}
 \sum_{j\geq 0}\l( \sum_l \frac{\beta_{l,\nu,k,j}}{(2l)!} 
(q-j+ \nu -2k)\cdots(q-j+ \nu -2k+2l-1)\,  \int_{M_b} c_{l,\nu,k,j}(b,x)_{[d]} \  d{\rm vol}_{g_b(x)}\r)\cdot 
 \end{equation*}
 $$\x \  \frac{1}{s-q+j +2k -\nu}$$
 
 \begin{equation}\label{local zeta sing structure d}
=\  \sum_{j\geq 0}\frac{1}{\G(q-j+\nu-2k)}\l( \sum_l \frac{\beta_{l,\nu,k,j}}{(2l)!} 
(q-j+ \nu -2k+2l-1))! \,  \int_{M_b} c_{l,\nu,k,j}(b,x)_{[d]} \  d{\rm vol}_{g_b(x)}\r)
 \end{equation}
 $$\x \  \frac{1}{s-q+j+2k-\nu}.$$
 \vskip 3mm
 Crucially, the factor $\G(q-j+\nu-2k)$ is independent of $l$. Since $l$ parametrises the summands contributing to $\sqrt{{\bf r}}_{\nu - 2k -1 - j}[\mu]$, the factor is thus `universal' in relating corresponding coefficients of the singularity expansions. Further contributions to the coefficient of the distribution $(t-i0)^{-q+j + 2k - \nu}$ coming from other summands \eqref{R d k summand} with corresponding `$j,k,\nu$'s which give the same numerical value $j+ 2k - \nu$, and with varying form degree $d$, leave unchanged the factor $1/\G(q-j+\nu-2k)$  by which the coefficient of $(t-i0)^{-q+j + 2k - \nu}$ in \eqref{wave trace d sing str}
 differs from the residue of $1/(s-q+j+2k-\nu)$ in \eqref{local zeta sing structure d}. We may thus better write this in a less specific manner. Noting in view of the assumption \eqref{max order m} that  $i := j +2k-\nu \geq j\geq 0$, then by  summing over all coefficients of terms contributed by the various components \eqref{R d k summand}, we may rewrite \eqref{wave trace d sing str} as 
\begin{equation}\label{wave trace d sing str 2}
\Tr_\tmb(e^{-it\sqrt{\Qs}}) \ \sim \ \sum_{i\geq 0} {\bf a}_i \cdot  (t-i0)^{-q+i } \hskip 1cm {\rm t\ near \ 0}
\end{equation}
with 
$${\bf a}_i = \sum_{d=0}^{\dim B} ({\bf a}_i)_{[d]} \in \Aa^*(B),$$
as expected from \eqref{wave invariants at 0} (there are no log-terms since $\sqrt{\Qs}^2$ is a vertical differential operator), and infer from  \eqref{local zeta sing structure d} that this is related to the zeta function singularity structure of $\sqrt{\Qs}$ by 
\begin{equation}\label{local zeta sing structure d 2}
Z_{\sqrt{\Qs}}(s) \ \sim \ \sum_{i\geq 0} \frac{{\bf a}_i}{\G(q-i)} \cdot   \frac{1}{s-q+i}.
\end{equation}
or, rather
\begin{equation}\label{local zeta sing structure d 3}
\G(s)\,Z_{\sqrt{\Qs}}(s) \ \sim \ \sum_{i\geq 0} {\bf a}_i \cdot   \frac{1}{s-q+i}.
\end{equation}
Given $ Z_\Qs(s) = Z_{\sqrt{\Qs}}(2s)$ one therefore has: 
 \begin{prop}\label{wave zeta heat d}
With the above assumptions, the vertical wave trace invariants \eqref{wave invariants at 0} at $t=0$ coincide with the locally determined vertical zeta trace coefficients \eqref{zeta form singularity structure} and the vertical heat trace invariants \eqref{heat trace singularity structure} - one has
\begin{equation}\label{wave zeta heat}
    {\bf b}_i = \frac{{\bf a}_i}{2} \hskip 5mm {\rm in} \ \Aa^*(B).
\end{equation}
    where ${\bf b}_i := \sum_d {\bf b}_{i,d}$.   
\end{prop}

\subsection{Cohomology classes of wave supertrace characteristic forms:}

Let $\Ee=\Ee^+\oplus\Ee^- \to M$ be a $\mathbb{Z}_2$-graded vector bundle over the total space of the fibration $\pi$. Then there is an induced $\Z_2$-grading 
\begin{equation}\label{grading}
\Aa(M/B,\Ee)  = \Aa^+(M/B,\Ee) \oplus \Aa^-(M/B,\Ee)
\end{equation}
on $\Aa(M/E,\Ee) := \G(M,\pi^*(\wedge T^*_B)\otimes \Ee\otimes|\wedge_{\pi}|^{1/2})$ 
with 
$$\Aa^+(M/B,\Ee) = \sum_l  \Aa^{2l}(M/B,\Ee^+)\oplus \Aa^{2l-1}(M/B,\Ee^-)$$ 
and 
$$ \Aa^-(M/B,\Ee) = \sum_l  \Aa^{2l}(M/B,\Ee^-)\oplus \Aa^{2l-1}(M/B,\Ee^+).$$ 
A superconnection \cite{Q,B,BGV} on $\pi_*(\Ee)$ adapted to a 1st-order 
vertical self-adjoint elliptic differential operator 
\begin{equation}\label{graded Dirac}
\Ds = \begin{bmatrix}
  0 & \Ds^- \\
  \Ds^+ & 0 \\
\end{bmatrix} \ \in\Aa^0(M/B,\Psi^r(\Ee))    
\end{equation}
is a differential  operator $\Bf$ on $\Aa(B,\pi_*(\Ee))  = \Ci(M,\pi^*(\wedge T^*
B)\otimes \Ee\otimes|\wedge_{\pi}|^{1/2})$ of odd-parity  for \eqref{grading} and satisfying
\begin{equation}\label{e:superconnection}
    \Bf( \pi^*\omega\wedge\psi) = \pi^* d\omega \wedge \psi +
    (-1)^{|\omega|}\pi^*\omega\wedge\Bf(\psi)
\end{equation}
for $\omega\in\Aa(B)$ and $\psi\in \Aa(B,\pi_*(\Ee))$ and such
that $\Bf_{[0]} = \Ds \ ,$ where 
$$\Bf = \Bf_{[0]} + \Bf_{[1]} + \cdots + \Bf_{[\dim B]}, \hskip 10mm \Bf_{[i]}: \Aa^d(B,\pi_*(\Ee))\to
\Aa^{d+i}(B,\pi_*(\Ee)).$$ 
From \eqref{e:superconnection}, 
$\Bf_{[1]}$ is a classical  connection
while $\Bf_{[i]}( \pi^*\omega\wedge\psi) = (-1)^{|\omega|}\pi^*\omega\wedge\Bf_{[i]}\psi$ for $i\neq 1$ is a vertical (smooth family of) differential operators
$\Bf_{[i]}\in\Aa^i(M/B,\Psi(\Ee))$. In a local weak
trivialization $\Aa(U/B,\pi_*(\Ee)_{|U}) \cong \Aa(U)\otimes
\Ci(M_{b_0},\Ee^{b_0})$ around $b_0\in U$ the superconnection $\Bf$ takes the local
coordinate form $\Bf_{|U} = d_{U} + \sum_I D_I db_I,$ where $db_I
= db_{i_1}\ldots db_{i_m}$ and $D_I$ a classical differential operator on
$\Ci(M_{b_0},\Ee^{b_0})$.\\

The curvature of $\Bf$ is the vertical differential operator 
 $\Bf^2 \in \Aa(M/B,\pi_*(\Ee))$ with
 \begin{equation}\label{graded Laplacian}
     \Bf^2_{[0]} = \Ps^2 =  \begin{bmatrix}
  \Ds^-\Ds^+ &  0\\
  0 & \Ds^+\Ds^- \\
\end{bmatrix} \in\Aa^0(M/B,\Psi^r(\Ee)).
 \end{equation} 
$\Bf^2$ is thus parameter elliptic (with Agmon
angle $\pi$). The Bismut superconnection \eqref{B-superconnection} 
 provides a canonical analytic representative $\Str_\tmb(e^{-t\Bf^2})$ (for $t>0$) for the Chern character of the index bundle \cite{B}, \cite{BGV}. Here, we are interested rather in the wave supertrace characteristic form (for $t\in\R^1$)
\begin{equation}\label{superwavetrace}
    \Str_\tmb\l(e^{-it\sqrt{\Bf^2}}\r) \in \Aa(B)\ox_{\Ci(B)}\Sss^\prime(\R).
\end{equation} 
The proof of the analytic families index formula uses a $t$-rescaling  of $\Bf$ in which $\Bf_{[j]}$ is multiplied by $t^{(-j+1)/2}$,
yielding a transgression formula for the Chern character form with limits that exist as $t\to\infty$ and as $t\to 0+$. Here, no such transgression formula is sought (or apparent), but it is relevant to consider \eqref{superwavetrace} thus rescaled. We use the equivalent but modified rescaling  (replacing $\sqrt{t}$ by $t$) - define
$$\delta_t:\Aa^j(B)\to \Aa^j(B), \ \ \ \delta_t(\omega_{[j]}) = t^{-j}\omega_{[j]},$$
and then on $\Bs_{[k]}\in \If^m(M\times_B N/B, E^*\boxtimes F, \Lambda^\tmnb)_{[k]}$ by
$t\delta_t \cdot \Bs_{[k]}\cdot \delta_t\ii = t^{-k+1}\Bs_{[k]}$
and on a superconnection
$\Bf_t = t\delta_t \cdot \Bf \cdot \delta_t\ii = t\Bf_{[0]} + \Bf_{[1]} + t^{-1}\Bf_{[2]}+\cdots + t^{-\dim B}\Bf_{[\dim B]}.$
\begin{lem}\label{rescaled wave trace}
One has
 \begin{equation}\label{rescaled delta t}
      e^{-i\,{\rm sgn}(t)\sqrt{\Bf_t^2}} =  \delta_t \cdot e^{-it\sqrt{\Bf^2}}\cdot \delta_t\ii
 \end{equation}  
 and hence 
 \begin{equation}\label{str rescaled delta t}
      \Str(e^{-i\,{\rm sgn}(t)\sqrt{\Bf_t^2}}) = \delta_t \l(\Str_\tmb(e^{-it\sqrt{\Bf^2}})\r),
 \end{equation} 
i.e.
$\Str(e^{-i\,{\rm sgn}(t)\sqrt{\Bf_t^2}})(\phi) = \delta_t \l(\Str_\tmb(e^{-it\sqrt{\Bf^2}})(\phi)\r)$ in $\Aa(B)$ for each  $\phi\in\Sss(\R^1).$
\end{lem}
{\it Proof:}  Since $(\delta_t \cdot \Bf^2 \cdot \delta_t\ii -
\la\Is)\ii = \delta_t \cdot (\Bf^2  -
\la\Is)\ii \cdot \delta_t\ii$, then mindful of \eqref{1/m}
\begin{eqnarray*}
\sqrt{\Bf_t^2}   &:=& \Bf_t^2 (\Bf_t^2)^{-1+\frac{1}{2}}_\pi \\
&=& t^2 \delta_t \cdot \Bf^2\cdot \delta_t\ii  \int_{C}\la_\pi^{-\frac{1}{2}}\, (t^2 \delta_t \cdot \Bf^2\cdot \delta_t\ii  -
\la\Is)\ii \dl\la \\  
&=&   \delta_t \cdot t^2 \Bf^2  \int_{C}\la_\pi^{-\frac{1}{2}}\, (t^2  \Bf^2 -
\la\Is)\ii \, \dl\la \cdot \delta_t\ii  \\  
&=&   \delta_t \cdot t^2 |t|^{-1} \Bf^2  \int_{C}\la_\pi^{-\frac{1}{2}}\, (\Bf^2 -
\la\Is)\ii \, \dl\la \cdot \delta_t\ii \\
&= & \delta_t \cdot |t| \sqrt{\Bf^2} \cdot \delta_t\ii
\end{eqnarray*}
for $t\neq 0$, and hence
\begin{equation}\label{sqrt delta t}
 {\rm sgn}(t)  \sqrt{\Bf_t^2}  =  \delta_t \cdot \l(t \sqrt{\Bf^2}\r) \cdot \delta_t\ii.
\end{equation}
Since $\delta_t\cdot \sqrt{\Bf^2} \cdot \delta_t\ii = \sqrt{\delta_t \cdot\Bf^2\cdot \delta_t\ii}$ then  
$\delta_t \cdot t\sqrt{\Bf^2} \cdot \delta_t\ii = t(\sqrt{\Ds^2} +  \delta_t \cdot (\sqrt{\Bf^2})_+\cdot \delta_t\ii)  $ and so 
\begin{eqnarray*}
e^{-i\,{\rm sgn}(t)\sqrt{\Bf_t^2}} &\stackrel{\eqref{sqrt delta t}}{=} & e^{-i\, t(\sqrt{\Ds^2} + \delta_t \cdot(\sqrt{\Bf^2})_+ \cdot \delta_t\ii)} \\
  & \stackrel{\eqref{Cauchy Qs}}{=} & \sum_{k=0}^{\dim B} (-it)^k  \int_{\sigma_k} e^{-i s_0 t\sqrt{\Ds^2}}  (\delta_t \cdot\sqrt{\Bf^2}_+ \cdot \delta_t\ii) \cdots (\delta_t \cdot\sqrt{\Bf^2}_+ \cdot \delta_t\ii) \,e^{-i s_k t \sqrt{\Ds^2} }  \ ds \\
& \stackrel{\eqref{Cauchy Qs}}{=} & \delta_t \cdot e^{-it\sqrt{\Bf^2}}\cdot \delta_t\ii.
\end{eqnarray*}
\eqref{str rescaled delta t} then follows since for any v-FIO $\Bs$ 
$$\Str(\delta_t \cdot \Bs \cdot \delta_t\ii)(\phi) 
=\sum_{k\geq 0} \Str(\delta_t \cdot \Bs_{[k]} \cdot \delta_t\ii)(\phi) 
=\sum_{k\geq 0} t^{-k} \Str(\Bs_{[k]})(\phi) 
=\delta_t\l(\Str(\Bs)(\phi)\r).$$ 
\hfill $\Box$  \vskip 3mm 

From \eqref{full vertical wave trace 3} and \lemref{rescaled wave trace} there is a singularity expansion
\begin{eqnarray}
    \Str_\tmb\l(e^{-i\,{\rm sgn}(t)\sqrt{\Bf_t^2}} \r)  &\sim&   \sum_{d= 0}^{\dim B} \sum_{j\geq 0}
{\bf a}_{j,d}\,(t + i0)^{-q + j -d} \label{wave trace sconn} \\ 
&& +  \sum_{d= 0}^{\dim B}\sum_{\kappa\neq 0}  \sum_{j\geq 0}  
{\bf a}_{\,\kappa, j,d}\,(t- L_\kappa + i0)^{-q + j-d} \log(t-L_\kappa+i0)^{\mathfrak{H}(q-j)} \notag
\end{eqnarray}
with differential form coefficients with $ {\bf a}_{j,d}, {\bf a}_{\,\kappa,j,d} \in\Aa^d(B)$; there are no log terms at $L_0=0$ since $\Bf^2$ is a vertical differential operator.
\begin{prop}\label{closed coeeficients}
   The coefficients ${\bf a}_{j,d}$ at the $t=0$ singularity are closed differential forms, representing cohomology classes in $H^d(B,\R^1)$.  If $L_\kappa : B\to \R^1$ is a constant function, then same is true of the coefficients ${\bf a}_{\,\kappa,j,d}$. 
\end{prop}
{\it Proof:} By \lemref{rescaled wave trace}  it is enough to prove this for the un-scaled wave trace $\Str_\tmb(e^{-it\sqrt{\Bf^2}})$. 
Let $\Us(t) = e^{-it\sqrt{\Bf^2}}$ and
let $\Vs(t)$ be the (super)commutator $[\Bf, U(t)]$ (noting that $\Us(t)$ has even parity). Then 
\begin{eqnarray*}
    \partial_t\Vs(t) = [\Bf,\Us^{\,\prime}(t))]  &= &  [\Bf,-i\sqrt{\Bf^2}\,\Us(t)] \\ 
    &= &  [\Bf,-i\sqrt{\Bf^2}] \,\Us(t)) - i\sqrt{\Bf^2} \,[\Bf,\Us(t))] \\ 
    &= &  - i\sqrt{\Bf^2}\, [\Bf,\Us(t)]
\end{eqnarray*}
since $[\Bf,-i\sqrt{\Bf^2}]  =0$ from $\sqrt{\Bf^2} = \Bf^2 \int_{C}\la_\pi^{-1+\frac{1}{m}}\, (\Bf^2 -
\la\Is)\ii \, \dl\la$ and $[\Bf,(\Bf^2 - \la\Is)\ii]=0$. Thus, $\Vs(t)$ for small $t$ is the unique solution to the ODE on $M$
$$\Vs^{\,\prime}(t) = - i\sqrt{\Bf^2}\,\Vs(t), \ \ \ \ \Vs(0) =0,$$
and hence must be identically zero, i.e. 
\begin{equation}\label{BU=0}
    [\Bf,\Us(t))]=0.
\end{equation}
(The solution $\Vs(t)$ exists for all $t$ via $\Vs(t+t') =\Vs(t)e^{-it'\sqrt{\Bf^2}} + e^{-it\sqrt{\Bf^2}} \Vs(t')$.) 
With $\Us$  the vertical wave operator with $t$ considered as a variable, then for a test function for $\rho\in\Sss(\R^1)$  $$\Us(\rho) = \int_{-\infty}^\infty \Us(t)\rho(t)\ dt \in \Aa((M\x_BM)/B,E^*\boxtimes E, \Lambda_\pi^{1/2}\boxtimes \Lambda_\pi^{1/2})$$ 
 is a vertical smoothing kernel, as is seen by repeated partial integrations using $\partial_t\Us(t) = -i\sqrt{\Bf^2}\Us(t)$, with vertical trace $$\Str(U(\rho)) = \int_{-\infty}^\infty  \Str_\tmb(e^{-it\sqrt{\Bf^2}})(\rho(t)) \ dt\in\Aa(B).$$ 
As such, by Lemma 9.15 of \cite{BGV} 
\begin{equation}\label{dstr}
    d_B\Str_\tmb\l(\Us(\rho)\r) = \Str_\tmb\l([\Bf,\Us(\rho)]\r).
\end{equation}
The vertical kernel of $\Us(t) = e^{-it\sqrt{\Bf^2}}$ can be written with respect to a local identification $\Aa(U_B, \pi_*(E)) \cong \Aa(U_B)\otimes \Ci(M_0,E_0\ox |\Lambda_{M_0}|^{1/2})$ defined by fibrewise coordinates over a chart $U_B\< B$ as 
$\sum_I \Us_I(t,b,m,m') \, db_I$ for $m,m'\in M_0$. Over $U_B$ the superconnection looks like
$$\Bf_{|U_B} = d_{U_B} + \sum_J D_J \, db_J$$
with $D_J$ a family of differential operators on $\Ci(M_0,E_0\ox |\Lambda_{M_0}|^{1/2})$ parametrised by $b\in U_B$, as observed in \cite{BGV} {\it loc. cit.} and as in that argument 
    $$[\Bf,\Us(\rho)]_{|U_B} = \int_{M_0}  d_{U_B} \int_{-\infty}^\infty \Us(t,b,m,m)_I \, db_I \, \rho(t) \ dt  +  \l[\sum_J D_J\, db_J, \Us(\rho)_{|U_B}\r].$$
But the supertrace of the commutator term vanishes, so the differential form \eqref{dstr} is given in coordinates over $U_B$ by 
$$ \int_{M_0} \sum_{i,I}\partial_{b_i}\int_{-\infty}^\infty \Us(t,b,m,m)_I \ db_i\wedge db_I \ \rho(t) \ dt$$
which since $ \Us(t,b,m,m)_I\rho(t)$ is smooth in $b$ and Schwartz class in $t$
$$  = \int_{M_0} \int_{-\infty}^\infty\sum_{i,I}\partial_{b_i}\Us(t,b,m,m)_I \ db_i\wedge db_I \, \rho(t) \ dt $$
$$= \Str_\tmb\l([\Bf,\Us(t)]_{|U_B}(\rho)\r)$$
Hence \eqref{BU=0} and \eqref{dstr} prove that $\Str_\tmb\l(\Us_t(\rho)\r)$ is a closed differential form on $B$ for any $\rho$. On the other hand, as in \eqref{sing exp basic}, \eqref{wave trace sconn} one has that $d_B\Str_\tmb\l(\Us_t(\rho)\r) $ is equal asymptotically to 
$$\sum_{d=0}^{\dim B}\sum_\kappa \sum_{j\geq 0} d_B{\bf a}_{\kappa,j,d } \,  \l(\int_0^\infty r^{q-j +d -1} e^{-ir(t-L_\kappa)} \ dr\r)(\rho)\hskip 10mm $$
$$\ \ \ \ \ \hskip 25mm + \ 
{\bf a}_{\kappa,j,d }\,d_BL_\kappa\,\l(\int_0^\infty ir^{q-j +d} e^{-ir(t-L_\kappa)} \ dr\r)(\rho)$$
with ${\bf a}_{0,j,d }= {\bf a}_{j,d }$ in \eqref{wave trace sconn}. If $L_\kappa$ is constant on $B$ (as for $L_0=0$) then  $d_BL_\kappa=0$, leaving only the sum over the $d_B{\bf a}_{\kappa,j,d } $ which since the test function $\rho$ is arbitrary implies these must all be zero.  For, taking $\rho$ with compact support in a sufficiently small interval around $L_\kappa$ we may restrict attention to a fixed $\kappa$, and moreover choosing $\rho$ analytic in $s= t-L_\kappa$, this leaves a sum 
$$ \sum_{d= 0}^{\dim B}\sum_{j = 0}^{k+d}  
d_B{\bf a}_{\,\kappa, j,d}\,(s + i0)^{-q + j-d} \log(s+i0)^{\mathfrak{H}(q-j)}(\rho) + o(1) =0, $$
and replacing $\rho(s)$ by its dilation $\rho_{\epsilon\ii}(s):=\rho(\epsilon\ii s)$ then this sum over log-homogeneous distributions becomes a sum of the form 
$$ \sum_{d= 0}^{\dim B}\sum_{j = 0}^{k+d}  
d_B{\bf a}_{\,\kappa, j,d}\,C_{\,\kappa, j,d}(\rho) \epsilon^{-q + j-d} (\log\epsilon)^{\mathfrak{H}(q-j)} + o(1) =0$$
as $\epsilon\to 0$ for constants $C_{\,\kappa, j,d}(\rho)$, which by the argument of \cite{H1} Lem. 3.2.1. immediately implies that $d_B{\bf a}_{\,\kappa, j,d}=0$ for $j \leq k+d$, thus leaving a sum 
$$ \sum_{d= 0}^{\dim B}\sum_{j > k+d}  
d_B{\bf a}_{\,\kappa, j,d}\,(s + i0)^{-q + j-d} \log(s+i0)^{\mathfrak{H}(q-j)}(\rho) + o(1) =0, $$
and by replacing $\rho(s)$ now by the reciprocal dilation $\rho_{\epsilon}(s):=\rho(\epsilon s)$ then again the argument of \cite{H1} Lem. 3.2.1. implies $d_B{\bf a}_{\,\kappa, j,d}=0$ for $j > k+d$. 
\hfill $\Box$  \vskip 3mm 

Generically, then, the coefficients ${\bf a}_{\,\kappa,j,d}$ do not represent cohomology classes. In cases where they are closed, such as for constant $L_\kappa$, one may conjecture that they will be exact forms due to the canceling of periodic geodesics of the Laplacians $D^*D$ and $DD^*$, and so zero in cohomology (but may, nevertheless, define secondary classes in cases where the wave trace coefficients are known to vanish). \\

 \begin{prop}\label{wave A hat}
 For the Bismut  superconnection
 \begin{equation}
     \frac{1}{2}\,{\bf a}_{q+ d,d} =   (2\pi )^{-\frac{q}{2}}\l(\int_{M/B}\widehat{A}(\smb)
    \, \ch(E)\r)_{[d]},
 \end{equation}
and the vertical wave trace  $\Str_\tmb(e^{-i\,{\rm sgn}(t)\sqrt{\Bf_t^2}} )$ expansion \eqref{wave trace sconn} then has a limit as $t\to 0$. 
 \end{prop}
{\it Proof:}  The heat vertical operator $e^{-(\Bf_{\sqrt{t}})^2}$  can be constructed via the Duhamel method \cite{BGV}, as mirrored here for the vertical wave operator, or, since it has advantages for computations for $t$ near to 0, the equivalent functional calculus construction for which one considers a contour $\Cc$ coming in on a ray with argument in $(0,\pi/2)$, encircling the origin, and leaving on a
ray with argument in $(-\pi/2,0)$. Then for $t>0$ 
\begin{equation}\label{v-heat operator}
e^{-t\Bf^2} = \int_{\Cc}
e^{-t\la}(\Bf^2 - \la\Is)\ii \
  \dl\la
\end{equation}
$$ = e^{-t\Ds^2} + \sum_{d=1}^{\dim B}\sum_{p_1+\cdots +p_k = d}  (-1)^k \int_{\Cc}
e^{-t\la}(\Ds -
\la \Is)\ii \Bf^2_{[p_1]}(\Ds - \la \Is)\ii \cdots \Bf^2_{[p_k]}(\Ds -
\la \Is)\ii \
  \dl\la.$$
The partial integration identity 
$e^{-t\Bf^2} = \int_{\Cc}
t^{-m} e^{-t\la}\partial_\la^m (\Bf^2 - \la\Is)\ii \
  \dl\la$ shows that $e^{-t\Bf^2}$ is a vertical smoothing operator and hence has a vertical heat supertrace
\begin{equation*}
\Str_\tmb(e^{-t\Bf^2}) = \int_\tmb \tr_m(\Hs_t(m,m)) \  \in \Aa(B),
\end{equation*}
where in local fibrewise coordinates $m= (b,x), m'=(b,y)$ the vertical heat kernel is computed using the  vertical resolvent symbol \eqref{d parametrix symbols} by
$$ \Hs_t(m,m')=\int_{\R^q} \int_{\Cc}
e^{-t\la} e^{i\langle\xi. x-y\rangle} \
  {\bf r}[\la](b,x,y,\xi) \ \dl\la\ \dl\xi.$$
By similar computations to those above  we hence infer an asymptotic expansion
$\Str_\tmb(e^{-t\Bf^2}) =  \sum_{d=0}^{\dim B} \sum_{j\geq 0} {\bf a}_{j,d}\, t^{\frac{- q +j }{2}}$ as $t\to 0$ with  ${\bf b}_{j,d} = \int_\tmb  \Omega_{j,d} \in\Aa^d(B)$
 and $\Omega_{j,d}$ built from a sum of integrals
$\int_{\R^q} \int_{\Cc} e^{-\la}\,{\bf r}_{\nu - (k+1)m - j}[\la](b,x,\xi)_{[d]} \ \dl\la\,\dl\xi |dx|$, and hence since  $\Str_\tmb(e^{-(\Bf_{\sqrt{t}})^2}) = t^{1/2}\delta_{t^{1/2}} \Str_\tmb(e^{-t\Bf^2})\delta_{t^{1/2}}\ii$ an asymptotic expansion
\begin{equation*}
 \Str_\tmb(e^{-(\Bf_{\sqrt{t}})^2}) = \sum_{d=0}^{\dim B} \sum_{j\geq 0} {\bf b}_{j,d}\,t^{\frac{-q+j -d }{2}}.
\end{equation*}
There are no $\log t$ terms because $\Bf$ is defined by differential operators. Choosing $\Bf_{\sqrt{t}}$ to specifically be the rescaled Bismut superconnection $t^{1/2}\Ds + \nabla^{\pi_*(\Ee)} + \frac{1}{4t^{1/2}}\,\cc(T)$ then  the local families index formula \cite{B,BGV} states that in this case the singular terms vanish and
\begin{equation}\label{local ind formula}
\lim_{t\to 0+} \Str_\tmb(e^{-(\Bf_{\sqrt{t}})^2}) =  (2\pi )^{-\frac{q}{2}}\int_{M/B}\widehat{A}(\smb)\,\ch(E)
\end{equation}
in $\Aa^*(B)$.
 Thus, 
${\bf b}_{j,d} = 0$ for $j <q +d$ while $\sum_{d=0}^{\dim B} {\bf b}_{d,d}$ is equal to the right-hand side of \eqref{local ind formula}. From \eqref{wave zeta heat}, the vertical wave invariants in \eqref{wave trace sconn} therefore satisfy  
$${\bf a}_{j,d} = 0 \ \ {\rm for} \ j < q+d$$
while the coefficient ${\bf a}_{d,d}$ of $(t-i0)^0$ is twice the $d$-form component of the right-hand side of \eqref{local ind formula}. 
\hfill $\Box$  \vskip 3mm 

{\it Comment}: \ In view of this, one may ask if the local families index formula is visible through the vertical wave trace $\Str_\tmb(e^{-(\Bf_{\sqrt{t}})^2})$ around $t=0$.  For, it is immediate that the differential form  order  $d=0$ component of $\Str_\tmb(e^{-(\Bf_{\sqrt{t}})^2})$ gives the pointwise index of the vertical Dirac-type  operator \eqref{graded Dirac} -- from \eqref{graded Laplacian} and \eqref{wave kernel 0}
\begin{equation}\label{str pointwise}
\Str_\tmb(e^{-it\sqrt{\Bf^2}})_{[0]}(b) = \Tr_\tmb(e^{-it\sqrt{\Ds_b^-\Ds_b^+}}) - \Tr_\tmb(e^{-it\sqrt{\Ds_b^+\Ds_b^-}})
\end{equation}
and, as in \cite{C,DG}, $\Tr_\tmb(e^{-it\sqrt{\Ds_b^-\Ds_b^+}})=\sum_{j\geq 0} e^{-it\la_j}$ with $\la_j$ the 
eigenvalue spectrum of $\Ds_b^-\Ds_b^+$, and by an elementary well-known argument all the non-zero eigenvalues cancel in the difference reducing \eqref{str pointwise} to $ \index (\Ds_b^+)$, and constant in $b$.
This is the same argument as for the heat trace, though the latter is an exact smooth object whilst the distributional wave trace is not.  For this reason one may expect the index not to be visible via the wave trace expansion; that is, the index is not visible from {\em asymptotic} wave trace formulae.  Concretely, the contour integral representation $I_{t,\Bf^2}$ in \eqref{wave operator contour} differs from $e^{-it\sqrt{\Bf^2}}$ by a smoothing operator term that carries the index. For the single Laplacian $\Ds_b^-\Ds_b^+$ one has
$$I_{t,\Ds_b^-\Ds_b^+} = e^{it\sqrt{\Ds_b^-\Ds_b^+}} - \Pi_{0,b} $$
with $\Pi_{0,b}$ a the smoothing-operator projection with range the kernel of $\Ds_b^-\Ds_b^+$, while
$$I_{t,\Bf^2} = e^{-it\sqrt{\Bf^2}} - \Pi_0 e^{-it\sqrt{(\Bf^2)_+}}\Pi_0  $$
with $\Pi_0 $ a corresponding smooth family of smoothing projections, noting that the finite sum $e^{-it\sqrt{(\Bf^2)_+}} = \sum_{k=0}\frac{(-it)^k}{k!}(\Bf^2)_+^k$ is a smooth function of $t$. The Chern character is not seen by the distributional wave trace, even though the local family index density is. \\[5mm]

\appendix

{\Large {\bf Appendix:}} \\

\section{Contextualising with classical Poisson summation}
 
 
 The wave trace formula \eqref{wave Poisson summation} due to \cite{C} and \cite{DG} is a generalisation to closed Riemannian manifolds of the distributional classical Poisson summation formula
\begin{equation}\label{classical p sum}
    \sum_{n\in\Z} e^{int} = 2\pi  \sum_{n\in\Z} \d(t-2\pi n).
\end{equation}
For this, on the left-hand side the integers $n^2$ are understood as the eigenvalues of the Laplacian $\D= -d^2/dx^2$ over the circle $S^1 = \R/2\pi\Z$ and the sum as the trace of the (full) wave operator $\cos(t\sqrt{\D})$, while the right-hand side is summed over the lengths $L_n =2\pi n$ of closed geodesics on $S^1$. \\

The equality \eqref{classical p sum} is an instance ($E=\R$, $G=2\pi\Z$, $X=S^1$) of computing a trace by pushing-down a Schwartz kernel from the total space $E$ of a smooth $G$-fibration $\rho : E\to X$ to the base manifold $X$. For this one assumes a (smooth, free and proper) Lie group $G$ action $E\x G \to E$ and a $G$-invariant density $d_g\mu$ on $E$. The mantra then is that objects on $X = E/G$ are the same thing as $G$-invariant objects on $E$.
So, to give a continuous operator $A: \Ci(X)\to \Ci(X)$  is in this mantra the same thing as to give the corresponding continuous operator $\Aa :  \Ci(E)^G \to \Ci(E)^G$
on the space of $G$-invariant functions $f = g^*f$ on $E$. $\Aa$ is thus an operator on $\Ci(E)$ which commutes with the $G$-action   $  g^*(\Aa\psi) = \Aa(g^*\psi)$, 
which is equivalent to $G$-invariance of the Schwartz kernel $K_\Aa$ of $\Aa$
\begin{equation}\label{G inv kernel}
    K_\Aa(xg,yg) = K_\Aa(x,y)
\end{equation}
for $(x,y)\in E\x E$, 
and hence that $K_\Aa$ pushes-down to $X\x X$ as
\begin{equation}\label{push down K}
    K_A( \ol x, \ol y) := \int_{G} K_\Aa(x,yg) \ d_g\mu
\end{equation}
which we take, by the assumed correspondence, to be the Schwartz kernel of $A$. Here, $(x, y)$ can be any element in the product fibre $\rho\ii(\ol x) \x \rho\ii(\ol y)$ and if $T(x,y)$ denotes the right-hand side then $T(xg_1, yg_2) = T(x,y)$ for all $g_1,g_2\in G$. The operator $A$ on $X$ is by construction the push-down $\rho_*(\Aa)$  
\begin{equation}\label{push down operator}
   A\phi = \rho_*(\Aa)(\phi):= \Aa(\rho^*\phi),
\end{equation}
or, equivalently, via the Schwartz kernel theorem, given \eqref{G inv kernel} holds then  $\rho$ induces by duality with pull-back on functions a push-forward on distributions and $K_A = (1\x\rho_*) K_\Aa$.
There is the question of which $\Aa$ on $E$ pushes-down to a given $A$ on $X$; for example, on a Riemannian fibration whether the Laplacian on $E$ pushes-down to the Laplacian on $X$, see \cite{GLP}. \\

Given continuity along the diagonal, setting $x=y$ \eqref{push down K} becomes the trace formula 
\begin{equation}\label{G inv kernel diag}
    K_A( \ol x, \ol x) = \int_{G} K_\Aa(x,xg) \ d_g\mu
\end{equation} and 
\begin{equation}\label{TrA G inv kernel}   
\Tr(A) = \int_X\int_{G} K_\Aa(x,xg) \ d_g\mu.
\end{equation}
The choice of $x\in\rho\ii(\ol x)$ does not matter locally, it is a choice of trivialization of the fibre $E_{\ol x} = \rho\ii(\ol x)$ defining a diffeomorphism   $E_{\ol x}\cong G.$ 
Globally, without such a choice, it is the identity
$ \ol\tau^* K_A = \int_{E/X} \tau_g^* K_\Aa$, with $\tau, \ol\tau_g$ the inclusion maps of the diagonal into $E\x E, X\x X$, respectively.   For examples of \eqref{TrA G inv kernel} in the case $\dim G >0$ see  \cite{Be} \S2 and \cite{BGV} Prop. 5.7. If $\dim G = 0$ and $\rho$ is a principal bundle, then $E$ is a covering space, $G=\G$  a discrete group,  and  \eqref{G inv kernel diag} becomes 
\begin{equation}\label{G=Gamma}   
K_A( \ol x, \ol x) = \sum_{\g\in\G} K_\Aa(x,x\g).
\end{equation}
Concretely, if $E$ is the universal cover of Riemannian $X$ then $\G = \pi_1(X)$ and integrating over $X$ gives the generic Poisson/ Arthur-Selberg summation formula
\begin{equation}\label{PAS}   
\Tr(A) = \sum_{l\geq 0}\sum_{\g} \int_D K_\Aa(x,x\g^l)
\end{equation}
in which the inner sum is over primitive geodesics $\g$ and $D\<E$ is a fundamental domain for the action of $\G$ on $E$. Taking $M=\R, \G=2\pi\Z, A= e^{\pm it\sqrt{d^2/dx^2}}$ (or appropriate linear combinations of it) one obtains back \eqref{classical p sum} by taking $K_\Aa$ to be the FIO kernel $K_\Aa(x,y) = F_{\xi \to x-y} \ii(\cos(|\xi| t))$, where $F$ is Fourier transform.\\

A vertical families Poisson summation formula considers this framework applied to a  fibration $\pi: M\to B$ of manifolds, as in \thmref{thm1} except here allowing non-compact $M$, but now with a fibrewise $G$-action on $M$; thus, an action which restricts to a $G=G_b$ action $M_b \x G \to M_b$ on each fibre of $\pi$, yielding a quotient fibration
$$\pi_G : M\sslash G \to B$$
with fibre $M_b/G = M_b/G_b$ at $b\in B$. The subscript on $G_b$ is to emphasize the specific action of $G$ on $M_b$. The objective is to push-down a vertical Schwartz kernel on $\Ks^M$ on the fibre product fibration $M\x_B M \to B$ to a vertical Schwartz kernel $\Ks^{M\sslash G}$ on $(M\sslash G) \x_B (M\sslash G) \to B$ -- that is, a family parametrised by $b\in B$ of kernels $\Ks^M_b$ on $M_b \x M_b$ to a family of kernels $K^{M\sslash G}_b$ on $M_b/G \x M_b/G = X_b \x X_b$ -- associated to the commutative diagrams of fibrations
\begin{equation*}
\begin{tikzcd}
M\   \arrow[r,  "\tau" ] \arrow[dr, "\pi"']
& \ M\sslash G \arrow[d, "\mu"]\\
& B
\end{tikzcd} \hskip 5mm {\rm and} \hskip 5mm 
\begin{tikzcd}
M\x_B M \   \arrow[r] \arrow[dr]
& \  (M\sslash G) \x_B (M\sslash G)  \arrow[d]\\
& B
\end{tikzcd}.
\end{equation*}
Given that the symmetry property $g^*\Ks^M = \Ks^M$ for the form valued kernels on $M$ hold, and so fibrewise \eqref{G inv kernel} holds for each $\Ks_b, $
\begin{equation}\label{G inv kernel vert}
    \Ks^M_b(mg,m'g) = \Ks^M_b(m,m'),    \ \  \pi(m)=\pi(m') =b,
\end{equation}
then the vertical family of push-down kernels $\Ks^{M\sslash G}$ on $ (M\sslash G) \x_B (M\sslash G) $  is 
\begin{equation}\label{push down K vert}
    \Ks_b^{M\sslash G} ( \ol m, \ol m') := \int_{G_b}  \Ks^M_b(m,m'g), 
\end{equation}
where it is assumed $\Ks_b^M$ includes a $G_b$-invariant density, and where $$\pi(m) = \pi(m') =b  \ \ {\rm  and} \ \  \tau(m) = \ol m \in M_b/G, \ \tau(m')=\ol m' \in M_b/G.$$ 
Let $\As^{M\sslash G}$ be the vertical family of operators along the fibres of $M\sslash G$ defined by $\Ks^{M\sslash G}$, so $\As^{M\sslash G} = \tau_*(\As^M)$ is the fibrewise push-forward. As in the single operator case,  $M\sslash G$ is the space of interest whilst it is envisioned that the fibres $M_b$ of the fibration $\pi$ be geometrically flat non-compact spaces simple enough to know the fibrewise kernel $K^M$ with some precision. If the  quotient fibres $M_b/G$ are compact and $\Ks_b^{M\sslash G}$ trace class, this yields the vertical trace formula on $M\sslash G$
\begin{equation*}\label{TrMG}
    \Tr_{(M\sslash G)/B} (\As^{M\sslash G}) = \int_{(M\sslash G)/B} \int_G \Ks^M(m,mg) \ \in\Aa^*(B),
\end{equation*}
so at  $b= \pi(m=m_b)$
\begin{equation*}\label{TrMG2}
    \Tr_{(M\sslash G)/B} (\As^{M\sslash G})(b) = \int_{M_b/G_b} \int_{G_b} \Ks^M(m_b,m_b g_b),
\end{equation*}
which for discrete $G=\G$ reads
\begin{equation}\label{TrMG3}
    \Tr_{(M\sslash \G)/B} (\As^{M\sslash \G})(b) = \sum_{\g_b\in\Gamma_b} \int_{M_b/\G_b}  \Ks^M(m_b,m_b \g_b).
\end{equation}
If  $M$ is itself compact then there is  the vertical trace $\Trmb(\As^M) =\int_{M/B} \Ks^M (m,m)$ also, but this will in general differ from $ \Tr_{(M\sslash G)/B} (\As^{M\sslash G}) $. \\

An elementary
example of vertical families Poisson summation is the fibrewise $2\Z$ action $\xi \cdot 2n =\xi + 4\pi n$ on the universal cover  fibration $\pi:\R^1\to S^1:= \R^1/\Z $  which quotients $\pi$ down to the  non-trivial double cover (spin) $\Z_2$-fibration $\mu: 2S^1:=\R^1\sslash 2 \Z \to S^1$. The fibre of $2S^1$ above $b\in S^1$ comprises two points $2S^1_b = \{\sqrt{b},-\sqrt{b}\}$ (with the fibres of $\pi$ over $\pm \sqrt{b}$ respectively a copy of the even and odd integers); if $S^1$ is identified with the geometric unit circle then $\sqrt{b} = e^{ip/2}$ and  $-\sqrt{b} = e^{i(p+2\pi)/2}$ some $p$  or any of its $4\pi\Z$ translates, and one has
\begin{equation*}
\begin{tikzcd}
\R^1 \ \  \arrow[r,  "\tau(p) = e^{ip/2}" ] \arrow[dr, "\pi(p) = e^{ip}"']
& \ \ 2S^1 \arrow[d, "\mu(\theta) = \theta^2"]\\
& S^1
\end{tikzcd}.
\end{equation*}
To a smooth family of Schwartz class functions $f: S^1\to \Sss(\R^1), \, b\mto f_b$, we associate the vertical (fibrewise) FIO kernel $K_{\Aa_b}(p,q) = \int^\infty_{-\infty} e^{i(p-q)\zeta} f_b(\zeta) \ d\zeta$ for $p,q\in\pi\ii(b)$.
$K_{\Aa_b}(p,q)$ is fibrewise $4\pi\Z$ invariant \eqref{G inv kernel} and so pushes-down to $2S^1$ as the vertical family $\Ks^f$ of kernels $K^{f_b}(\ol p,\ol q) = \sum_n K_{\Aa_b}(p,q +4\pi n)$ of a vertical FIO $\As^{2S^1}$ on $2S^1$, where $\tau(p) = \ol p, \tau(q)=\ol q\in 2S^1_b.$
 Each $K^{f_b}$ has just four possible values on $2S^1_b\x 2S^1_b$ 
  \begin{equation}\label{4 kernels}
      K^{f_b}(\sqrt{b},\sqrt{b}), \  K^{f_b}(\sqrt{b},-\sqrt{b}), \ K^{f_b}(-\sqrt{b},\sqrt{b}), \ K^{f_b}(-\sqrt{b},-\sqrt{b}),
  \end{equation} 
and is $\Z_2$-invariant so $K^{f_b}(\sqrt{b},\sqrt{b}) =  K^{f_b}(-\sqrt{b},-\sqrt{b})$ and equal to $\sum_n  F(f_b)(4\pi n)$ while $K^{f_b}(\sqrt{b},-\sqrt{b}) =  K^{f_b}(-\sqrt{b},\sqrt{b})=\sum_n  F(f_b)(2\pi (2n +1))$. 
Thus, integrating over the fibre, \eqref{TrMG3} is in this case the vertical trace formula 
$$\Tr_{2S^1/S^1} (\As^{2S^1})(b) = K^{f_b}(\sqrt{b},\sqrt{b}) +  K^{f_b}(-\sqrt{b},-\sqrt{b}) = 2\sum_n  F(f_b)(4\pi n).$$
On the other hand, in view of $\Z_2$ invariance, $\Ks^f$  in turn pushes-down  from $2S^1$ to the null-fibration $id: S^1\to S^1$  by \eqref{G=Gamma} as 
$$K_{S^1}(b) = K^{f_b}(\sqrt{b},\sqrt{b}) +  K^{f_b}(\sqrt{b},-\sqrt{b}) = \sum_n  F(f_b)(2\pi n)$$
which is the usual Poisson sum in \eqref{classical p sum} (note comment following \eqref{TrMG3}).\\

The push-down to the $k$-fold cover $kS^1 \to S^1$ is entirely similar, with the above sums replaced by a sum over the ${\rm k}^{th}$ roots of unity. Likewise, discrete quotients of higher dimensional rank vector bundles yield trace formulae over fibrations of tori. \\

\thmref{thm1} and \thmref{thm2} provide generalizations of such vertical families Poisson summation formulae to arbitrary compact Riemannian fibrations, coinciding with the generalization of \eqref{wave Poisson summation} to \eqref{classical p sum} in the case that $B$ is a single point manifold. \\

\section{Proof of \propref{tmxn}}

We shall show the dual identification for tangent bundles since the proof is perhaps a little more transparent, though the argument is entirely dualizable to give directly the result for cotangent bundles.  There is a  commutative diagram
\begin{equation*}
\begin{tikzcd}
M\x_B N \arrow{r}{\mu} \arrow[swap]{dr}[above]{\hskip 5mm \   s}\arrow{d}{\nu} &  M\arrow{d}{\pi} \\
N \arrow{r}{\sigma}   & B
\end{tikzcd}
\end{equation*}
from which one has (since $s=\pi \circ \mu = \sigma\circ \nu$) the exact sequence 
\begin{equation}\label{fibre exact}
0\to T_{M \x_B N}  \stackrel{\mu^* d\mu + \nu^* d\nu  }{\too} \mu^* T_M + \nu^* T_N \stackrel{ s^*d\pi \mu^\#  - s^*d \sigma \nu^\#}{\too}   s^* T_B \to 0,
\end{equation} 
where $\beta^\# :\beta^* E \to E$ is the covering bundle map for the pull-back of a vector bundle $E\to Y$ and smooth map $\beta: X\to Y$. Using \eqref{cotangent splitting} to write this as
\begin{equation}\label{two exact}
    0\to T_\tmnb + s^*T_B \stackrel{\mu^* d\mu + \nu^* d\nu  }{\too} \mu^* T_\tmb + s^*T_B + \nu^* T_\tnb + s^*T_B  \stackrel{ s^*d\pi \mu^\#  - s^*d \sigma \nu^\#}{\too}   s^* T_B \to 0,
\end{equation} 
\eqref{fibre exact}  is seen to split into the direct sum of the two short exact sequences 
\begin{equation}\label{1exact}
0\to T_\tmnb \stackrel{\mu^* d\mu + \nu^* d\nu  }{\too} \mu^* T_\tmb + \nu^* T_\tnb   \stackrel{ s^*d\pi \mu^\#  - s^*d \sigma \nu^\#}{\too}  0,
\end{equation} 
and 
\begin{equation}\label{2exact} 
0\to s^*T_B \stackrel{\mu^* d\mu + \nu^* d\nu  }{\too}  s^*T_B + s^*T_B  \stackrel{ s^*d\pi \mu^\#  - s^*d \sigma \nu^\#}{\too}   s^* T_B \to 0
\end{equation} 
with $s^*d\pi \mu^\#$ restricting on $\mu^* T_\tmb$ to the zero map, as $T_\tmb = \ker(d\pi)$, and likewise $s^*d \sigma \nu^\#$ to the zero map on $\nu^* T_\tnb$, hence the final arrow of \eqref{1exact}. On the other hand,  $\mu^* d\mu$ and $\nu^* d\nu$ are complementary linear projection operators to the vertical tangent spaces over $M$ and $N$ respectively, defining a bijection for the non-zero arrow of \eqref{1exact}. \eqref{1exact} is thus a vector bundle isomorphism 
\begin{equation}\label{1exact2}
T_\tmnb \stackrel{\mu^* d\mu + \nu^* d\nu  }{\too} \mu^* T_\tmb + \nu^* T_\tnb.  
\end{equation} 
The fibrewise projections $\mu^* d\mu$ and $\nu^* d\nu$ both reduce on the pull-back bundle $s^* T_B$ to the identity operator, so the first non-zero arrow of \eqref{2exact} is the diagonal map $\xi\mapsto (\xi,\xi)$ and the second is $(\xi,\eta) \mapsto \xi-\eta$.  Taking duals of \eqref{1exact2} gives the first identity of \eqref{fibre symp isom}. \\

 $T^*_\tmb = \cup_b T^*_{M_b}$ has a canonical vertical(fibrewise)-symplectic form $\omega_\tmb$ --
restricting to the standard symplectic two form $\omega_{M_b}$  on each $T^*_{M_b}$, but which is degenerate on the manifold $T^*_\tmb$.  Globally, to the projection 
$p:T^*_\tmb\to M$ one assigns the 1-form $\alpha$ on $T^*_\tmb$ by $\a_{m,\xi} = \xi\circ dp_{(m,\xi)}$, with $\xi\in T^*_m M_b$ (i.e $\pi(m) =b$) considered as in $T^*_M$ via inclusion (extending by zero)  \eqref{cotangent splitting}. Then the (degenerate) 2-form $\omega_\tmb = d\alpha$  coincides on restriction to each $T^*_{M_b}$ with the standard non-degenerate symplectic form $\omega_{M_b}$. From the first identity of \eqref{fibre symp isom} we obtain the vertical-symplectic form $\mu^* \omega_\tmb + \nu^*\omega_\tnb$ on $\mu^*T^*_\tmb + \nu^*T^*_\tnb $ which by convention is usually written as $ \omega_\tmb + \omega_\tnb$ on $T^*_\tmnb$. \\

The final identification of \eqref{fibre symp isom} is the twist fibrewise diffeomorphism defined by the map which acts by $-I$ on the fibres of $T^*_\tnb$, given in local fibrewise coordinates by $(b, x,y, \xi,\eta) \mto (b, x,\xi, y, -\eta)$, taking the vertical-symplectic form $ \mu^*\omega_\tmb + \nu^*\omega_\tnb$ on $T^*_\tmnb$ to the vertical-symplectic form $ \omega_\tmb - \omega_\tnb$ on $T^*_\tmb \x_B T^*_\tnb$; that is, to the degenerate 2-form which  coincides on restriction to the fibre $T^*_{M_b} \x T^*N_b$ with the (non-degenerate) symplectic 2-form  $\omega_{M_b} - \omega_{N_b}$. 
\hfill $\Box$ \\[3mm]

\end{document}